\documentclass[times]{nmeauth}

\usepackage{amsthm,amsmath}
\usepackage[utf8]{inputenc} 
\usepackage{bm}
\usepackage{dsfont}
\usepackage{amssymb}
\usepackage{mathtools}
\usepackage{graphicx}
\graphicspath{ {Figures/} }
\usepackage{subcaption}
\usepackage{lscape}
\usepackage{calrsfs}
\usepackage{algorithm2e}
\usepackage{hyperref}
\def\cref#1{{Eq.~(\ref{#1})}}
\def\Cref#1{{Eq.~(\ref{#1})}}
\usepackage{algorithmic}
\usepackage{varwidth}
\usepackage[toc,page]{appendix}
\usepackage[acronym,smallcaps]{glossaries}
\usepackage[parfill]{parskip}

\makeatletter
\def\myfnt{\ifx\protect\@typeset@protect\expandafter\footnote\else\expandafter\@gobble\fi}
\makeatother

\newcommand{\Fa}{\mathcal{F}}
\newcommand{\La}{\mathcal{L}}

\newcommand{\Pa}{\mathcal{P}}

\newcommand{\dpar}[2]{\frac{\partial{#1}}{\partial{#2}}}

\newacronym{pod}{POD}{Proper Orthogonal Decomposition}
\newacronym{fsi}{FSI}{Fluid-Structure Interaction}
\newacronym{ale}{ALE}{Arbitrary-Lagrangian-Eulerian}
\newacronym{svd}{SVD}{Singular Value Decomposition}
\newacronym{vms}{VMS}{Variational Multiscale}
\newacronym{fe}{FE}{Finite Element}
\newacronym{oss}{OSS}{Orthogonal Subgrid-Scales}
\newacronym{sgs}{SGS}{Subgrid-Scales}
\newacronym{bdf}{BDF}{Backward Differentiation Formula}

\DeclareMathAlphabet{\pazocal}{OMS}{zplm}{m}{n}
\DeclareMathAlphabet\mathbfcal{OMS}{cmsy}{b}{n}
\allowdisplaybreaks

\begin{document}


\title{Three-Field Fluid-Structure Interaction by Means of the Variational Multiscale Method}


\author{Alexis Tello\affil{1}, Ramon Codina\affil{1,2}}

\address{
\affilnum{1}Universitat Polit\`ecnica de Catalunya (UPC), Jordi Girona 1-3, 08034 Barcelona, Spain\\
\affilnum{2}International Centre for Numerical Methods in Engineering (CIMNE),
C/ Gran Capit\`a S/N, 08034 Barcelona, Spain}
\corraddr{Corresponding author ramon.codina@upc.edu (R. Codina)}
\date{}

\begin{abstract} 
Three-field \gls{fsi} formulations for fluid and solid are applied and compared to the standard two field-one field formulation for fluid and solid, respectively. Both formulations are applied in a non linear setting for a Newtonian fluid and a neo-Hookean solid in an updated Lagrangian form, both approximated using finite elements and stabilized by means of the Variational Multiscale (VMS) Method to permit the use of arbitrary interpolations. It is shown that this type of coupling leads to a more stable solution. Even though the new formulation poses the necessity of additional degrees of freedom, it is possible to achieve the same degree of accuracy as standard FSI by means of coarser meshes, thus making the method competitive. We enhance the stability of the formulation by assuming that the sub-grid scales of the model evolve through time. Benchmarking of the formulation is carried out. Numerical results are presented for semi-stationary and a fully transient cases for well known benchmarks for 2D and 3D scenarios.
\end{abstract}


\maketitle



\section{Introduction}

The term `mixed methods' in mechanics is generally applied to a formulation that approximates separately different variables, for example stress and displacement fields. In our case we apply the term to deal with the splitting of the stress tensor into its corresponding deviatoric and spherical parts, leading to a displacement--stress--pressure formulation approximated by the finite element (FE) method. This kind of stress splitting techniques are by no means new and have been shown to work properly for both solid mechanics and fluid dynamics. In \cite{Codina2009} it was shown that it is possible to approximate successfully the Stokes problem by means of a three-field splitting and a stabilization of the Galerkin formulation using a Variational Multiscale (VMS) approach, in particular assuming that the sub-grid scales of the model belong to a space orthogonal to the space of the FE scale. Later, \cite{Chiumenti2015} applied this same three-field splitting technique to a linear solid mechanics setting, comparing it with a displacement--pressure splitting. This is the basis of our work in terms of solid mechanics; however there are other types of approximations to the elastic problem, for example \cite{Scovazzi2017} explores a velocity--stress splitting that proves to be convenient and robust for time evolution settings. The assumptions and approximations made to develop stress--displacement and strain--displacement formulations are detailed in \cite{Cervera2010} and applied in \cite{Cervera2015} to approximate compressible and incompressible plasticity by means of a VMS method, producing a method with enhanced stability and convergence properties in comparison to the displacement based --or irreducible-- formulation. In the context of geometrically nonlinear solid mechanics, a formulation accounting for the incompressible limit using a total Lagrangian approach is proposed in \cite{castanar19}.

In the field of fluid dynamics the spectrum is wider. In \cite{Castillo2014a,Castillo2014b} a FE velocity--stress--pressure formulation is applied to flows with non-linear viscosity, stabilized again by means of a VMS method. In \cite{castillo2018} the same method is shown to behave better than velocity--pressure approximations for high Weissenberg number flows of viscoelastic fluids. Later, \cite{Moreno2019a} would expand this work, reformulating into a logarithmic version of the problem able to cope with higher elastic effects.

Both in solids and in fluids, a major reason for using stabilized FE methods is that the Galerkin method is only stable for certain choices of the interpolating spaces for the unknowns, which turn out to be very restrictive. There are two inf-sup conditions to be met (see e.g. \cite{ruas-1985-1}), one between the displacement and the pressure space to yield stable pressures and another one between the stresses and the displacements in order to have control on the displacement gradients. In the case of fluids, there is also the need of using stabilized FE methods when convection dominates. The inf-sup conditions require complex interpolations, one of the first being the element proposed in \cite{marchal-crochet} and analyzed in \cite{fortin-pierre-1989-1}; see also \cite{ruas-1996-1,sandri-1993-1}. Stabilized FE methods allow one to use arbitrary interpolations \cite{Codina2009}, thus simplifying enormously the implementation.

In terms of the interaction problem, research can be broadly grouped into two categories based on how the FE mesh is treated, namely, conforming and non conforming methods. Essentially, conforming mesh methods consider interface conditions as physical boundary conditions, thus treating the interface as part of the solution. In this approach, the mesh reproduces or conforms to the interface; when the interface is moved it is also necessary to displace the mesh, which carries on all related problems of mesh recalculation and inherent instabilities of the method, be it staggered or monolithic, see \cite{Badia2008,Farhat2014,Farhat2006,Bord2013,Bazilevs2006,Bazilevs2008,LeTallec2001}.
On the other hand, non-conforming methods treat the interface and boundary as constraints imposed on the governing equations, which makes it possible to use meshes that do not reproduce the interface. The main problems in this case is the treatment of the interface conditions and the complexity of the formulation, see for example \cite{Baiges2017,Glowinski2017} for further reading. For a general review of significant Fluid-Structure Interaction (FSI) advances and developments, see \cite{Hou2012}. 

As far as we know, no attempt has been made to treat FSI using a three-field approach for both the solid and the fluid. The benefits of mixed formulations applied to either fluid or solid are evident, from the extension to different applications to better convergence properties. Our assumption is that the interaction problem can `inherit' these properties. Moreover, taking advantage of an equal matching of the unknowns (velocity/stress/pressure-displacement/stress/pressure), or field to field coupling, may lead to a more accurate solution. We show in this paper that this is in fact the case, and the gain in accuracy and the better behavior of coupling schemes may compensate the increase in the number of degrees of freedom (DOFs) with respect to irreducible formulations. 

The paper is organized as follows. In Section~2 we present the numerical approximation employed for the solid. Being this a non-standard three-field formulation, allowing to reach the incompressible limit, we describe it in some detail. The point of departure is an updated Lagrangian statement of the finite strain solid mechanics problems, assuming the solid to behave as a neo-Hookean material. We describe the splitting of the stress tensor that leads to the three-field formulation employed, the linearization, time integration and FE approximation using a VMS formulation. The next ingredient is the three-field formulation for the fluid, assumed to be Newtonian and incompressible, described in Section~3. The FSI problem is then described in Section~4, using as building blocks the solid and the fluid numerical models. No attempt has been made to design a particularly robust and efficient iterative scheme, and in this work we have restricted ourselves to the classical Dirichlet-Neumann coupling (the fluid uses the motion of the interface determined by the solid and the solid is computed with the stresses on the interface provided by the fluid). Numerical results are then presented in Section~5, and concluding remarks close the paper in Section~6. 

\section{Three-field elasto-dynamic solid equations}\label{sec:sldsection}

In this section we give an overview of the three-field elasto-dynamic solid equations, their manipulation to obtain the governing equations we use in this work and the way to approximate them using a stabilized FE formulation.

The conservation of momentum for a solid can be written as:
\begin{align}\label{gov_eqn}
    \rho_{\textrm{sl}}\partial_{tt} d_{i} -\partial_{j} \sigma_{ij} &= \rho_{\textrm{sl}} f_{\textrm{sl},i} \quad
    \textrm{in} \; \Omega_{\rm sl}, ~t \in ]0, t_{\rm f}[, 
\end{align}
where $\Omega_{\rm sl}$ is the domain of $\mathbb{R}^{n_{\textrm{d}}}$ where the solid moves during the time interval $[0,t_{\rm f}]$, $n_{\rm d} = 2$ or 3 being the number of space dimensions, $\partial_t$ denotes the partial time derivative (and thus $\partial_{tt}$ is the second derivative with respect to time), and $\partial_i$ the derivative with respect to the $i$th Cartesian coordinate $x_i$, $i = 1,\dots, n_{\rm d}$. The unknowns of the problem are the displacement field and the Cauchy stress tensor, with Cartesian components $d_i$ and $\sigma_{ij}$, respectively, $i,j = 1,\dots, n_{\rm d}$. Here and in what follows, vectors and tensors are assumed to be represented by their Cartesian components. In particular, $f_{\textrm{sl},i}$ is the vector of external forces acting on the solid. Its density is given by $\rho_{\textrm{sl}}$. Finally, in Eq.~(\ref{gov_eqn}) and below repeated indexes imply summation over the number of space dimensions.

Eq.~(\ref{gov_eqn}) will be expressed in an updated Lagrangian reference system, and therefore we will need the constitutive law to be given in terms of the Cauchy stress $\sigma_{ij}$. If the expression of the stress in terms of the displacement is inserted into Eq.~(\ref{gov_eqn}), we will call {\em irreducible} the resulting problem posed for $d_i$ alone. Initial and boundary conditions have to be appended to this equation.

\subsection{Split of the stress tensor and field equations}\label{sec:tensor_split}

We can define the pressure and the deviatoric component of the stress tensor as:  
\begin{subequations}
\begin{align}
        p_{\textrm{sl}}    & = \frac{1}{n_{\textrm{d}}}\sigma_{kk}, \label{volumetric}\\
        s_{\textrm{sl},ij}  & = \sigma_{ij} - p_{\textrm{sl}} \delta_{ij}, \label{deviatoric}
\end{align}
\end{subequations}
where $\delta_{ij}$ is Kronecker's delta. This split can be used both in an updated Lagrangian and in a total Lagrangian framework (see \cite{castanar19} for the latter). Since we will restrict ourselves to the first option, the range of constitutive equations is more restrictive. In particular, we will consider the solid to behave as a neo-Hookean material, with constitutive law:
\begin{align}\label{neohook}
    \sigma_{ij} &= \frac{1}{J}\left[(\lambda \textrm{ln}(J)-\mu_{\textrm{sl}}) \delta_{ij} + \mu_{\textrm{sl}} b_{ij}\right],
\end{align}
where $J$ is the determinant of the displacement gradient $F_{iJ}$, $\lambda$ and $\mu_{\textrm{sl}}$ are Lame's parameters, $ b_{ij}$ is the left Cauchy tensor and $ b_{ii}$ its trace. These are defined as follows:
\begin{align}\label{nonli_terms}
    J = \textrm{det}(F_{iJ}),  \quad
    F_{iJ} = \dpar{x_i}{X_J},  \quad
     b_{ij} = F_{iK}F_{jK},
\end{align}
where $x_i$ are the coordinates defined in the spatial frame of reference and $X_I$ are the coordinates defined in the material frame of reference. Note also that here and below we will use lower case letters for spatial indexes and upper case for material ones.

\Cref{volumetric} can be expressed as:
\begin{align}\label{vol_press}
    p_{\textrm{sl}} &=\frac{1}{J}\left[ (\lambda \textrm{ln}(J)-\mu_{\textrm{sl}}) + \frac{\mu_{\textrm{sl}}\, b_{ll}}{n_{\rm d}}\right] ,
\end{align}
and the deviatoric part of the stress, \cref{deviatoric}, as:
\begin{align}\label{dev_stress}
    s_{\textrm{sl},ij} &= \frac{\mu_{\textrm{sl}}}{J}\left( b_{ij}-\frac{ b_{ll}}{n_{\textrm{d}}} \delta_{ij}\right).
\end{align}
Finally we can write the resulting system of equations as an alternative to \cref{gov_eqn} as:
\begin{align}\label{trifil_prob}
    \rho_{\textrm{sl}}\partial_{tt} d_{i} -\partial_{j} s_{\textrm{sl},ij} - \partial_{i}p_{\textrm{sl}}&= \rho_{\textrm{sl}} f_{\textrm{sl},i}, \nonumber\\
     \frac{J}{2\mu_{\textrm{sl}}}s_{\textrm{sl},ij}       &=\frac{1}{2}\left( b_{ij}-\frac{ b_{ll}}{n_{\textrm{d}}}\delta_{ij}\right) ,  \\
    \frac{J}{\lambda}p_{\textrm{sl}} &= \textrm{ln}(J) + \frac{\mu_{\textrm{sl}}}{\lambda}\left(\frac{b_{ll}}{n_{\rm d}} -1 \right). \nonumber
\end{align}
This is the three-field form of the solid mechanics equations we shall consider, the unknowns being $d_i$, $s_{\textrm{sl},ij} $ and $p_{\textrm{sl}}$. We shall write these unknowns as $\bm{y}_{\textrm{sl}} = [\bm{d},\bm{s}_{\textrm{sl}},p_{\textrm{sl}}]$, and denote by $\mathcal{A}(\bm{y}_{\textrm{sl}})$ the spatial nonlinear operator associated to \cref{trifil_prob}, so that these equations can be written as
\begin{align}
\mathcal{M} \partial_{tt} \bm{y}_{\textrm{sl}} + \mathcal{A}(\bm{y}_{\textrm{sl}}) 
= \bm{F}_{\rm sl},\label{eq:sol-compact}
\end{align}
where $\mathcal{M} := {\rm diag} [\rho_{\textrm{sl}} \bm{I}_{n_{\rm d}} ,  {\bf 0} , 0]$, $\bm{I}_{n_{\rm d}}$ being the identity on vectors and $\bm{F}_{\rm sl} := [ \rho_{\textrm{sl}} \bm{f}_{\textrm{sl}} ,  {\bf 0}, 0]$.

\subsection{Linearization}\label{sec:solid_linearization}

In order to linearize the problem, we can re-write \cref{trifil_prob} by approximating our variables  in terms of their increment $\delta \bm{y}_{\textrm{sl}}=[\delta \bm{d},\delta \bm{s}_{\textrm{sl}}, \delta p_{\textrm{sl}}]$ in the following form:
\begin{equation}\label{increment}
    \bm{y}_{\textrm{sl}} = \tilde{\bm{y}}_{\textrm{sl}} + \delta \bm{y}_{\textrm{sl}},
\end{equation}
where $\tilde{\bm{y}}_{\textrm{sl}} = [\tilde{\bm{d}},\tilde{\bm{s}}_{\textrm{sl}},\tilde{p}_{\textrm{sl}}]$ is the vector consisting of previously known values of $\bm{y}_{\textrm{sl}}$. The Newton-Raphson linearization is obtained inserting this split in the equations to be solved and neglecting quadratic terms of the increments. \Cref{trifil_prob} can be re-written as (see the Appendix for the linearization of $J$ and of its logarithm):
\begin{subequations}\label{eqn:trifil_prob_LU}
\begin{flalign}
    &\quad \rho_{\textrm{sl}}\partial_{tt}\delta d_{i}-\partial_{j} \delta s_{\textrm{sl},ij} - \partial_{i} \delta p_{\textrm{sl}}= \rho_{\textrm{sl}} f_{\textrm{sl},i} +\partial_{j} \tilde{s}_{\textrm{sl},ij} + \partial_{i} \tilde{p}_{\textrm{sl}} - \rho_{\textrm{sl}}\partial_{tt} \tilde{d}_{i}, &\\
    &\quad \frac{\tilde{J}}{2\mu_{\textrm{sl}}}\delta s_{\textrm{sl},ij} - \frac{1}{2}\left(\dpar{\delta d_i}{X_K}\tilde{F}_{jK}+\tilde{F}_{iK}\dpar{\delta d_j}{X_K}\right) + \frac{1}{n_{\textrm{d}}}\left( \tilde{F}_{lK}\dpar{\delta d_l}{X_K}\right)\delta_{ij} &\\ 
    &\qquad\qquad\qquad\qquad + \frac{\tilde{J}}{2\mu_{\textrm{sl}}}\left(\tilde{F}_{Kl}^{-1}\dpar{\delta d_l}{X_K}\right)\tilde{s}_{\textrm{sl},ij}=   \frac{1}{2}\left(\tilde{b}_{ij}-\frac{\tilde{ b}_{ll}}{n_{\textrm{d}}} \delta_{ij}\right)-\frac{\tilde{J}}{2\mu_{\textrm{sl}}}\tilde{s}_{\textrm{sl},ij},\nonumber&\\
    &\quad \frac{\tilde{J}}{\lambda}\delta p_{\textrm{sl}} + \left[\left(\frac{\tilde{J}\tilde{p}_{\textrm{sl}}}{\lambda}-1\right)\tilde{F}_{Kl}^{-1}-\frac{2\mu_{\textrm{sl}}}{\lambda n_{\textrm{d}}}\tilde{F}_{lK}\right]\dpar{\delta d_l}{X_K} = -\frac{\tilde{J}}{\lambda}\tilde{p}_{\textrm{sl}} + \textrm{ln}(\tilde{J})+\frac{\mu_{\textrm{sl}}}{\lambda}\left(\frac{\tilde{b}_{ii}}{n_{\textrm{d}}} -1 \right). &
\end{flalign}
\end{subequations}
Let us denote by $\La_{\textrm{sl}}(\tilde{\bm{y}}_{\textrm{sl}} ; \delta\bm{y}_{\textrm{sl}})$ the spatial linear operator on $\delta\bm{y}_{\textrm{sl}}$, for given $\tilde{\bm{y}}_{\textrm{sl}}$, appearing in the left-hand-side (LHS) of these equations. Using the notation in \cref{eq:sol-compact}, we may write
\begin{align}
\mathcal{M} \partial_{tt} \delta \bm{y}_{\textrm{sl}} + \La_{\textrm{sl}}(\tilde{\bm{y}}_{\textrm{sl}} ; \bm{\delta y}_{\textrm{sl}})
=  \bm{F}_{\rm sl}  - \mathcal{M} \partial_{tt} \tilde{\bm{y}}_{\textrm{sl}} - \mathcal{A} (\tilde{\bm{y}}_{\textrm{sl}}).
\label{eq:sol-lin-ibvp}
\end{align}
In this way, upon convergence $\delta\bm{y}_{\textrm{sl}} \approx \bm{0}$ and $\tilde{\bm{y}}_{\textrm{sl}}$ is the solution of \cref{eq:sol-compact}.

\subsection{Initial and boundary-value problem}

The problem to be solved consists in finding 
$\bm{y}_{\rm sl} : \Omega_{\rm sl} \times ]0,t_{\rm f}[ \longrightarrow \mathbb{R}^{n_{\rm d}}\times  \mathbb{R}^{n_{\rm d}}\otimes \mathbb{R}^{n_{\rm d}}\times \mathbb{R}$ as the solution to \cref{eq:sol-compact} and such that
\begin{align*}
    &\quad d_{i}                  = d_{i,D}            &   &\text{on}\;\Gamma_{\textrm{sl},D},   ~~t\in\;]0,t_{\textrm{f}}[,  \\ 
    &\quad n_{\textrm{sl},j}\sigma_{ij}= t_{{\rm sl}, i}              &   &\text{on}\;\Gamma_{\textrm{sl},N},  ~~t\in\;]0,t_{\textrm{f}}[,  \\
    &\quad d_{i}               = d_{i}^{0}          &   &\text{in}\;\Omega_\textrm{sl},      ~~t=0, \\
    &\quad \dot{d}_{i}         = \dot{d}_{i}^{0}    &   &\text{in}\;\Omega_\textrm{sl},      ~t=0, 
\end{align*}
where $\bm{d}^{0}$ is a prescribed initial displacement, $\dot{\bm{d}}^{0}$ is a prescribed initial velocity, $\bm{d}_D$ is a prescribed displacement on the boundary $\Gamma_{\textrm{sl},D}$, $\bm{t}_{\rm sl}$ is a prescribed traction on the boundary $\Gamma_{\textrm{sl},N}$, and $\bm{n}_{\textrm{sl}}$ is the normal to the solid domain. Here it is assumed that $\Gamma_{\textrm{sl},D}$ and $\Gamma_{\textrm{sl},N}$ are a partition of $\partial\Omega_{\rm sl}$, although later we will also introduce the interface boundary with the fluid.

\subsection{Weak form }\label{sec:weak_trifil}

Let us denote by $\langle\cdot,\cdot\rangle_\omega$ the integral of the product of two functions in a domain $\omega$, with the subscript omitted when $\omega = \Omega_{\rm sl}$. When the two functions belong to $L^2$, we will replace this symbol by $(\cdot,\cdot)_\omega$. Let also $Y_{\rm sl}$ be the space where the unknown $\bm{y}_{\rm sl}$ must belong for each time $t$, satisfying the Dirichlet boundary conditions, and let $Y_{{\rm sl},0}$ be the space of functions satisfying the homogeneous counterpart of these Dirichlet conditions. 

The weak form of the three-field elasto-dynamic solid equations consists in finding $\bm{y}_{\textrm{sl}} = [\bm{d},\bm{s} _{\textrm{sl}},p_{\textrm{sl}}] : [0,t_{\rm f}] \longrightarrow Y_{\rm sl}$ such that 
\begin{subequations}\label{sys:weaktrifil} 
\begin{align}
   &\quad (\rho_{\textrm{sl}}\partial_{tt}\delta d_{i},e_{i}) +  \left(\dpar{e_{i}}{x_j},\delta s_{\textrm{sl},ij}\right) + \left(\dpar{e_{i}}{x_i},\delta p_{\textrm{sl}}\right)= \langle e_{i},\rho_{\textrm{sl}} f_{\textrm{sl},i}\rangle + \langle e_{i},t_{{\rm sl},i}\rangle_{\Gamma_{\textrm{sl},N}} - (\rho_{\textrm{sl}}\partial_{tt} \tilde{d}_{i},e_{i}) \nonumber \\
   &\qquad\qquad -\left(\dpar{e_{i}}{x_j},\tilde{s}_{\textrm{sl},ij}\right) - \left(\dpar{e_{i}}{x_i},\tilde{p}_{\textrm{sl}}\right) , \label{eq:weak1}\\
    &\quad \left(\xi_{\textrm{sl},ij},\frac{\tilde{J}}{2\mu_{\textrm{sl}}}\delta s_{\textrm{sl},ij}\right) - \left(\xi_{\textrm{sl},ij},\frac{1}{2}\left(\dpar{\delta d_i}{X_K}\tilde{F}_{jK} +\tilde{F}_{iK}\dpar{\delta d_j}{X_K}\right)\right)   \nonumber\\
    &\qquad\qquad + \left(\xi_{\textrm{sl},ij},\frac{1}{n_{\textrm{d}}}\left( \tilde{F}_{lK}\dpar{\delta d_l}{X_K}\right)\delta_{ij}\right) + \left(\xi_{\textrm{sl},ij},\frac{\tilde{J}}{2\mu_{\textrm{sl}}}\left(\tilde{F}_{Kl}^{-1}\dpar{\delta d_l}{X_K}\right)  \tilde{s}_{\textrm{sl},ij}\right) 
    \nonumber\\
    &\qquad\qquad = \left(\xi_{\textrm{sl},ij},\frac{1}{2}\left(\tilde{ b}_{ij}-\frac{\tilde{ b}_{ll}}{n_{\textrm{d}}} \delta_{ij}\right)\right)-\left(\xi_{\textrm{sl},ij},\frac{\tilde{J}}{2\mu_{\textrm{sl}}}\tilde{s}_{\textrm{sl},ij}\right) , \label{eq:weak2}\\
    &\quad \left(q_{\textrm{sl}},\frac{\tilde{J}}{\lambda}\delta p_{\textrm{sl}}\right) + \left(q_{\textrm{sl}},\left[\left(\frac{\tilde{J}\tilde{p}_{\textrm{sl}}}{\lambda}-1\right)\tilde{F}_{Kl}^{-1}-\frac{2\mu_{\textrm{sl}}}{\lambda n_{\textrm{d}}}\tilde{F}_{lK}\right]\dpar{\delta d_l}{X_K}\right) \nonumber\\
    &\qquad\qquad= \left(q_{\textrm{sl}},-\frac{\tilde{J}}{\lambda}\tilde{p}_{\textrm{sl}} + \textrm{ln}(\tilde{J})+\frac{\mu_{\textrm{sl}}}{\lambda}\left(\frac{\tilde{b}_{ii}}{n_{\textrm{d}}} -1 \right)\right) ,  \label{eq:weak3} 
\end{align}
\end{subequations}
for all $\bm{z}_{\textrm{sl}}=[\bm{e},\bm{\xi}_{\textrm{sl}},q_{\textrm{sl}}]\in Y_{{\rm sl},0}$, $t\in ]0,t_{\rm f}[$, and satisfying initial conditions in a weak sense. 

Let us define the form $B_{\textrm{sl}}$ as:
\begin{flalign}
&\quad B_{\textrm{sl}}(\tilde{\bm{y}}_{\textrm{sl}};\delta \bm{y}_{\textrm{sl}},\bm{z}_{\textrm{sl}}) =    \left(\dpar{e_{i}}{x_j},\delta s_{\textrm{sl},ij}\right) + \left(\dpar{e_{i}}{x_i},\delta p_{\textrm{sl}}\right) &\nonumber \\ 
&\qquad\qquad+ \left(\xi_{\textrm{sl},ij},\frac{\tilde{J}}{2\mu_{\textrm{sl}}}\delta s_{\textrm{sl},ij}\right) -\left(\xi_{\textrm{sl},ij},\frac{1}{2}\left(\dpar{\delta d_i}{X_K}\tilde{F}_{jK}+\tilde{F}_{iK}\dpar{\delta d_j}{X_K}\right)\right)& \nonumber\\ 
&\qquad\qquad+ \left(\xi_{\textrm{sl},ij},\frac{1}{n_{\textrm{d}}}\left( \tilde{F}_{lK}\dpar{\delta d_l}{X_K}\right)\delta_{ij}\right) \nonumber + \left(\xi_{\textrm{sl},ij},\frac{\tilde{J}}{2\mu_{\textrm{sl}}}\left(\tilde{F}_{lK}^{-1}\dpar{\delta d_l}{X_K}\right)  \tilde{s}_{\textrm{sl},ij}\right) & \nonumber\\
&\qquad\qquad + \left(q_{\textrm{sl}},\frac{\tilde{J}}{\lambda}\delta p_{\textrm{sl}}\right) + \left(q_{\textrm{sl}},\left[\left(\frac{\tilde{J}\tilde{p}_{\textrm{sl}}}{\lambda}-1\right)\tilde{F}_{Kl}^{-1}-\frac{2\mu_{\textrm{sl}}}{\lambda n_{\textrm{d}}}\tilde{F}_{lK}\right]\dpar{\delta d_l}{X_K}\right),& \label{trifil_bilinear}
\end{flalign}
and a form $L_{\textrm{sl}}$ as,
\begin{flalign}\label{trifil_linear}  
&\quad L_{\textrm{sl}}(\tilde{\bm{y}}_{\textrm{sl}};\bm{z}_{\textrm{sl}}) 
=  \langle e_{i},\rho_{\textrm{sl}} f_{\textrm{sl},i}\rangle 
+ \langle e_{i},t_{\textrm{sl},i} \rangle_{\Gamma_{{\rm sl},N}} -\left(\dpar{e_{i}}{x_j},\tilde{s}_{\textrm{sl},ij}\right) 
- \left(\dpar{e_{i}}{x_i},\tilde{p}_{\textrm{sl}}\right) \nonumber &\\
&\qquad\qquad + \left(\xi_{\textrm{sl},ij},\frac{1}{2}\left(\tilde{ b}_{ij}-\frac{\tilde{ b}_{ll}}{n_{\textrm{d}}} \delta_{ij}\right)\right)-\left(\xi_{\textrm{sl},ij},\frac{\tilde{J}}{2\mu_{\textrm{sl}}}\tilde{s}_{\textrm{sl},ij}\right) & \\
&\qquad\qquad + \left(q_{\textrm{sl}},-\frac{\tilde{J}}{\lambda}\tilde{p}_{\textrm{sl}} + \textrm{ln}(\tilde{J})+\frac{\mu_{\textrm{sl}}}{\lambda}\left(\frac{\tilde{b}_{ii}}{n_{\textrm{d}}} -1 \right)\right),\nonumber&
\end{flalign}
which enable us to write \cref{sys:weaktrifil} in the following simplified form:
\begin{equation} \label{eqn:trifil_bilinear_compressed}  
(\rho_{\textrm{sl}}\partial_{tt} \delta d_{i},e_{i}) +  B_{\textrm{sl}}(\tilde{\bm{y}}_{\textrm{sl}};\delta \bm{y}_{\textrm{sl}},\bm{z}_{\textrm{sl}}) = L_{\textrm{sl}}(\tilde{\bm{y}}_{\textrm{sl}};\bm{z}_{\textrm{sl}})
- (\rho_{\textrm{sl}}\partial_{tt} \tilde{d}_{i},e_{i}), 
\end{equation}
for all $\bm{z}_{\textrm{sl}}:=[\bm{e},\bm{\xi}_{\textrm{sl}},q_{\textrm{sl}}] \in Y_{{\rm sl},0}$. Initial conditions have to added to this variational equation.

\subsection{Time discretization}

Let us consider a uniform partition of the time interval $]0,t_{\textrm{f}}[$ of size $\delta t$, and let us denote with superscript $n$ the time level. For the temporal discretization the following second order Backward Differences scheme (BDF2) will used:
\begin{equation*} 
    \bm{a}^{n+1} = \frac{1}{\delta t^{2}}\left(2\bm{d}^{n+1} -
    5\bm{d}^{n} + 4\bm{d}^{n-1} - \bm{d}^{n-2}\right), 
\end{equation*}
where $\bm{d}^{n+1}$ and $\bm{a}^{n+1}$ are approximations to the position and acceleration ($\partial_{tt}\bm{d}$) vectors at time $t^{n+1} = (n+1)\delta t$. Note that it is possible to use any other time integration scheme.

\subsection{Galerkin spatial discretization}\label{sldsectiondiscrete}

Let $\Pa_h$ denote a FE partition of the solid domain $\Omega_{\textrm{sl}}$. The diameter of an element domain $K \in \Pa_h$ is denoted by $h_K$ and the diameter of the FE partition by $h=\textrm{max}\{h_K|K \in \Pa_h\}$. As usual, a subscript $h$ is used to refer to FE functions and spaces. In particular, we can construct the approximating space for the unknown and the test functions, $Y_{{\rm sl},h}$ and $Y_{{\rm sl},h,0}$, respectively, in the usual manner. We consider here conforming approximations. 

The Galerkin FE approximation to the problem can be written as: find  $\bm{y}_{\textrm{sl},h} = [\bm{d}_h,
\bm{s}_{{\rm sl},h} ,p_{\textrm{sl},h}] : [0,t_{\rm f}] \longrightarrow Y_{{\rm sl},h}$ such that 
\begin{align}\label{disc_SLD} 
    (\rho_{\textrm{sl}}\partial_{tt}\delta \bm{d}_h,\bm{e}_h) + B_{\textrm{sl}}(\tilde{\bm{y}}_{\textrm{sl},h};\delta \bm{y}_{\textrm{sl},h},\bm{z}_{\textrm{sl},h}) & = L_{\textrm{sl},h}(\tilde{\bm{y}}_{\textrm{sl}};\bm{z}_{\textrm{sl},h}) - (\rho_{\textrm{sl}}\partial_{tt}\tilde{\bm{d}}_h,\bm{e}_h),
\end{align}
in $]0,t_{\rm f}[$ for all $\bm{z}_{\textrm{sl},h} \in  Y_{{\rm sl},h,0}$, and satisfying the initial conditions weakly. Considering the problem discretized in time, at each time step the problem to be solved until convergence is
\begin{align}
    (\rho_{\textrm{sl}}\delta \bm{a}^{n+1}_h,\bm{e}_h) + B_{\textrm{sl}}(\tilde{\bm{y}}^{n+1}_{\textrm{sl},h};\delta \bm{y}^{n+1}_{\textrm{sl},h},\bm{z}_{\textrm{sl},h}) & = L_{\textrm{sl},h}(\tilde{\bm{y}}^{n+1}_{\textrm{sl},h};\bm{z}_{\textrm{sl},h}) - (\rho_{\textrm{sl}}\tilde{\bm{a}}^{n+1}_h,\bm{e}_h),\label{disc_SLDn} 
\end{align}
for $n = 1,2,\dots$ and for all $\bm{z}_{\textrm{sl},h} \in  Y_{{\rm sl},h,0}$, $\bm{d}_h^0$ and $\bm{d}_h^1$ being given by the initial conditions and $ \bm{a}^{2}$ properly initialized (for example using a first order BDF scheme).
 
\subsection{A VMS method for an abstract stationary linear problem}\label{sec:VMS}

Problem (\ref{disc_SLDn}) is unstable, unless stringent requirements are met for the interpolating spaces of the variables in play (see, e.g., \cite{ruas-1985-1,Codina2009}). In order to be able to use arbitrary interpolations, a stabilized FE method is required. Here we describe the VMS approach we follow, first for an abstract linear stationary problem. This summary is required to extend it to {\rm second} order problems in time, a class of problems not considered previously.

Let us consider a generic bilinear form $B(\bm{y},\bm{z})$ and a linear form $L(\bm{z})$, where $\bm{y}$ is a vector of $n$ unknowns and $\bm{z}$ is the vector of test functions of the problem $B(\bm{y},\bm{z}) = L(\bm{z}) $ for all $\bm{z}$. 
Let $Z$ be the functional space where the continuous problem is posed (for simplicity, the same for unknowns and test functions) and $Z_{h} \subset Z$ the FE approximation. The idea of VMS methods is to split our unknowns into a FE part and a sub-grid scale (or simply sub-scale) that needs to be modelled \cite{Hughes1998,Codina2018}. Thus, let  $Z = Z_h \oplus \breve{Z}$, where $\breve{Z}$ is the complement of $Z_h$ in $Z$. This split will cause associated splittings $\bm{y} = \bm{y}_{h} + \breve{\bm{y}}$ and $\bm{z} = \bm{z}_{h} + \breve{\bm{z}}$ and, due to the linearity of $B$, we can write the problem as:
\begin{subequations}\label{eq:generic}
\begin{align}
B(\bm{y}_{h},\bm{z}_{h}) + B(\breve{\bm{y}},{\bm{z}}_h) &= L(\bm{z}_{h}) && \forall \bm{z}_{h} \in Z_{h}, \label{fem_generic}\\
B(\bm{y}_{h},\breve{\bm{z}}) + B(\breve{\bm{y}},\breve{\bm{z}}) &= L(\breve{\bm{z}}) && \forall \breve{\bm{z}} \in \breve{Z}.\label{sgs_generic}
\end{align}
\end{subequations}
It is possible to express the generic form $B(\bm{y},\bm{z})$ as:
\begin{subequations} \label{Bil_gen_vol_surf}
\begin{align}
    B(\bm{y},\bm{z}) &= \sum_{K}\left[\langle\La(\bm{y}),\bm{z}\rangle_{K} + \langle\Fa(\bm{y}),\bm{z}\rangle_{\partial K} \right]
    =   \sum_{K}\left[\langle\bm{y},{\cal L}^\ast \bm{z}\rangle_{K} + \langle\bm{y},{\cal F}^\ast(\bm{z})\rangle_{\partial K} \right],
\end{align}
\end{subequations}
where $\La(\cdot)$ is now the linear operator of the problem being solved and $\Fa(\cdot)$ is the associated flux operator acting on the inter-elemental boundaries $\partial K$, whereas ${\cal L}^\ast$ and ${\cal F}^\ast$ are the corresponding formal adjoints (see \cite{Codina2018} for details). 

Using \cref{Bil_gen_vol_surf} we can avoid computing derivatives of the sub-scale in \cref{fem_generic} and write problem (\ref{eq:generic}) as
\begin{subequations}\label{eqn:bilinear_red_gen}
\begin{align}
B(\bm{y}_{h},\bm{z}_{h}) + \sum_{K}\left[\langle\breve{\bm{y}},\La^{*}(\bm{z}_{h})\rangle_{K} +          
 \langle\breve{\bm{y}},\Fa^{*}(\bm{z}_{h})\rangle_{\partial K}\right] & = L(\bm{z}_{h}) &&\forall \bm{z}_{h} \in Z_{h},        \label{bilinear_u_redef_gen} \\ 
\sum_K \left[\langle\breve{\bm{z}},\La(\breve{\bm{y}}) + \La({\bm{y}_h})\rangle_{K} +  
\langle\breve{\bm{z}},\Fa(\breve{\bm{y}}) + \Fa({\bm{y}_h})\rangle_{\partial K} \right] & = L(\breve{\bm{z}}) && \forall \breve{\bm{z}} \in \breve{Z}. \label{bilinear_u_redef2_gen} 
\end{align}
\end{subequations}
If we assume that the inter-elemental fluxes across the element are continuous, the second term on the LHS of equation \cref{bilinear_u_redef2_gen} vanishes. We can express this equation as:
\begin{align}\label{sgs_gen} 
\La(\breve{\bm{y}}) = \breve{\Pi} (\bm{r}), 
\end{align}
in the space of sub-scales, where we have defined the residual $\bm{r}=\bm{f} - \La(\bm{y}_{h})$, $\La(\bm{y}) = \bm{f}$ is the linear equation being approximated and $\breve{\Pi}$ is the projection onto the space of sub-grid scales. We favor the option of taking this space as $L^2$ orthogonal to the FE space (what we call {\em orthogonal sub-grid scales}), which yields $\breve{\Pi} = {\Pi}^\bot_h = I - {\Pi}_h$, $I$ being the identity operator and  ${\Pi}_h$ the $L^2$ projection onto the FE space. However, in order to simplify the exposition, in this paper we will consider the classical option in residual-based stabilized formulations $\breve{\Pi} = I$. See \cite{Codina2018} and references therein for further discussion.

Assuming that $ \breve{\Pi} = I$, one can now use different arguments (see \cite{Codina2018}) to approximate the solution of \cref{sgs_gen} within each element $K$ as 
\begin{equation} 
\breve{\bm{y}} = \bm{\tau}_{K}\bm{r}, \label{var_sgs_gen} 
\end{equation} 
where $\bm{\tau}_{K}$ is the matrix of stabilization parameters. It is problem dependent, the particular expression we use for the FSI problem is presented later. 

For simplicity, in the following we shall omit the element boundary term in \cref{bilinear_u_redef_gen}, although it would be necessary if discontinuous pressure and stress interpolations are used (see \cite{Codina2009}). Thus, the final problem is obtained inserting \cref{var_sgs_gen} in \cref{bilinear_u_redef_gen}; this problem is posed in terms of the FE unknown only.

\subsection{Dynamic sub-grid scales for second order equations in time}\label{sec:dyn_sgs_trifil}

Let us consider now the time evolution version of the previous linear problem, considering second order time derivatives, i.e., the problem to be solved is now 
$$\partial_{tt}\bm{y} + \La(\bm{y}) = \bm{f},$$ 
with adequate initial and boundary conditions. In our particular FSI problem, we only have time derivatives of the displacement; the final equations to be solved are described in the next subsection.

The time evolution counterpart of  \cref{bilinear_u_redef_gen}, neglecting the element boundary term, is 
\begin{align}
(\partial_{tt}\bm{y}_h , \bm{z}_{h} ) + \sum_K\langle\partial_{tt}\breve{\bm{y}} , \bm{z}_{h} \rangle_K
+ B(\bm{y}_{h},\bm{z}_{h}) + \sum_{K}\langle\breve{\bm{y}},\La^{*}(\bm{z}_{h})\rangle_{K} & = L(\bm{z}_{h}) &&\forall \bm{z}_{h} \in Z_{h}, \label{eq:DSGS}
\end{align}
whereas the counterpart of approximation (\ref{var_sgs_gen}) is 
\begin{equation} 
\partial_{tt}\breve{\bm{y}} + \bm{\tau}_{K}^{-1}\breve{\bm{y}} =\bm{r} - \partial_{tt}\bm{y}_h. \label{eqn:sgs_dyn_disp} 
\end{equation} 
This equation is approximated using finite differences in time. Suppose that the finite difference scheme employed leads to the following approximation at time $t^{n+1}$:
\begin{equation}\label{eqn:sgs_accel} 
\partial_{tt}\breve{\bm{y}} \vert_{t^{n+1}}\approx \frac{\gamma_{1}\breve{\bm{y}}^{n+1}}{\delta t^{2}} -\frac{\breve{\bm{y}}^{\theta}}{\delta t^{2}},
\end{equation} 
where $\breve{\bm{y}}^{\theta}$ is the collection of terms belonging to the previous time steps in our time integration scheme and $\gamma_1$ is the coefficient multiplying the variable at step $n+1$. In this way, for a BDF2 scheme $\gamma_{1} = 2$ and $\breve{\bm{y}}^{\theta} = 5\bm{y}^{n} - 4\bm{y}^{n-1} +\bm{y}^{n-1}$. 
Substituting \cref{eqn:sgs_accel} into \cref{eqn:sgs_dyn_disp} and regrouping all unknowns from the current time step $n+1$ on the LHS and all known terms of time step $n$ into the right-hand-side (RHS) we obtain:
\begin{equation}\label{eqn:sgs_dyn_disp_with_accel} 
\frac{\gamma_{1}\breve{\bm{y}}^{n+1}}{\delta t^{2}}  + \bm{\tau}_{K}^{-1}\breve{\bm{y}}^{n+1} = \bm{R}^{n+1} + \frac{\breve{\bm{y}}^{\theta}}{\delta t^{2}},
\end{equation} 
where $\bm{R}^{n+1}$ includes all time discretized terms of the RHS of \cref{eqn:sgs_dyn_disp}. We thus have:
\begin{equation}\label{eqn:sgs_dyn} 
\breve{\bm{y}}^{n+1}  = \bm{\tau}_{t}\left(\bm{R}^{n+1} + \frac{\breve{\bm{y}}^{\theta}}{\delta t^{2}} \right),\qquad
\bm{\tau}_{t} :=\left(\frac{\gamma_{1}}{\delta t^{2}}\bm{I}_{n}  + \bm{\tau}_{K}^{-1}\right)^{-1},
\end{equation} 
where $\bm{I}_n$ is the $n\times n$ identity matrix. From \cref{eqn:sgs_dyn_disp} we also have that
\begin{equation}\label{eqn:sgs_fem} 
\langle\partial_{tt}\breve{\bm{y}},\bm{z}_{h}\rangle_K  = \langle\bm{R} - \bm{\tau}_{K}^{-1}\breve{\bm{y}} ,\bm{z}_{h}\rangle_K,
\end{equation} 
Then, after time discretization and substituting \cref{eqn:sgs_dyn} into \cref{eqn:sgs_fem} we have:
\begin{equation}\label{eqn:sgs_dyn_final} 
\langle\partial_{tt}\breve{\bm{y}} \vert_{t^{n+1}},\bm{z}_{h}\rangle_K \approx \langle(\bm{I}_n - \bm{\tau}_{K}^{-1}\bm{\tau}_{t})\bm{R}^{n+1},\bm{z}_{h}\rangle_K - \left\langle\bm{\tau}_{K}^{-1}\bm{\tau}_{t}\frac{\breve{\bm{y}}^{\theta}}{\delta t^{2}},\bm{z}_{h}\right\rangle_K.
\end{equation} 
The final time discrete problem for the FE unknown $\bm{y}_h^{n+1}$ is obtained by replacing \cref{eqn:sgs_dyn_final} and \cref{eqn:sgs_dyn} into the time discrete form of \cref{eq:DSGS}. Note that both in \cref{eqn:sgs_dyn_final} and in \cref{eqn:sgs_dyn} we need $\breve{\bm{y}}^{\theta}$, i.e., the sub-scales of previous time steps need to be stored at the numerical integration points. They will act as internal variables in a solid mechanics problem.

The idea presented is the same as for first order problems in time (see \cite{Codina2007}), the main difference being that now the sub-scales need to be stored in more previous time steps to allow computing an approximation to the second derivative. Even though it is not the purpose of this paper to study dynamic sub-scales for second order problems in detail, we have found them crucial to improve the behavior of iterative schemes. In particular, most of our FSI cases shown later would not converge without the use of dynamic sub-grid scales, their effect has been found critical in these cases.
 
\subsection{Stabilization of the non-linear three-field solid problem}\label{sec:stab_trifil}

Following the development from the previous section we can express in particular the forms and operators for our stabilized three-field solid problem. If we consider the splitting  of the unknown in the FE component and the sub-scale, instead of  \cref{disc_SLDn} we obtain:
\begin{subequations}\label{eqn:split_nonli_trifil}
\begin{align}
    & (\rho_{\textrm{sl}}\delta \bm{a}_h,\bm{e}_{h}) + \sum_K \langle \rho_{\rm sl} \breve{\bm{a}} ,  \bm{e}_{h}\rangle_K
    + B_{\textrm{sl}}(\tilde{\bm{y}}_{\textrm{sl},h};\delta \bm{y}_{\textrm{sl},h},[\bm{e}_{h},{\bf 0},0]) \nonumber\\
    & \qquad +  \sum_K \left\langle\dpar{e_{h,i}}{x_j}, \breve{s}_{{\rm sl},ij}\right\rangle_K 
    + \sum_K \left\langle\dpar{e_{h,i}}{x_i}, \breve{p}_{\rm sl}\right\rangle_K 
    = L_{\textrm{sl}}(\tilde{\bm{y}}_{\textrm{sl},h};[\bm{e}_{h},{\bf 0},0]) - (\rho_{\textrm{sl}}\tilde{\bm{a}}_h,\bm{e}_h), \label{eq:d1-a}\\
    & B_{\textrm{sl}}(\tilde{\bm{y}}_{\textrm{sl},h};\delta \bm{y}_{\textrm{sl},h},[0,\bm{\xi}_{\textrm{sl},h},0]) 
    + \sum_K\left\langle\xi_{\textrm{sl},h,ij}, \frac{\tilde{J}}{2\mu_{\textrm{sl}}}\breve{s}_{ij}\right\rangle_K
    - \sum_K\left\langle\xi_{\textrm{sl},h,ij},\dpar{\breve{d}_{i}}{X_K}\tilde{F}_{jK}\right\rangle_K \nonumber\\ 
    &\qquad+ \sum_K\left\langle\xi_{\textrm{sl},h,ij},\frac{1}{n_{\textrm{d}}}\tilde{F}_{lK}\dpar{\breve{d}_{l}}{X_K}\delta_{ij}\right\rangle_K 
    + \sum_K\left\langle\xi_{\textrm{sl},ij},\frac{\tilde{J}}{2\mu_{\textrm{sl}}}\tilde{F}_{Kl}^{-1}\dpar{\breve{d}_{l}}{X_K}\tilde{s}_{\textrm{sl},ij}\right\rangle_K  \nonumber\\
    &\qquad =   L_{\textrm{sl}}(\tilde{\bm{y}}_{\textrm{sl},h};[0,\bm{\xi}_{\textrm{sl},h},0]) , \\
    & B_{\textrm{sl}}(\tilde{\bm{y}}_{\textrm{sl},h};\delta \bm{y}_{\textrm{sl},h},[0,0,q_{\textrm{sl},h}]) 
    + \sum_K\left\langle q_{\textrm{sl},h},\frac{\tilde{J}}{\lambda}\breve{p}_{\rm sl}\right\rangle_K \nonumber\\
    & \qquad + \sum_K\left\langle q_{\textrm{sl},h},\left[\left(\frac{\tilde{J}\tilde{p}_{\textrm{sl},h}}{\lambda}-1\right)\tilde{F}_{Kl}^{-1}-\frac{2\mu_{\textrm{sl}}}{\lambda n_{\textrm{d}}}\tilde{F}_{lK}\right]\dpar{\breve{d}_{l}}{X_K}\right\rangle_K  
     = L_{\textrm{sl}}(\tilde{\bm{y}}_{\textrm{sl},h};[0,0,q_{\textrm{sl},h}]),
\end{align}
\end{subequations}
where the sub-scale at the current time step is $\breve{\bm{y}}_{\rm sl} = [ \breve{\bm{d}} ,  \breve{\bm{s}}_{\rm sl} ,  \breve{p}_{\rm sl}]$ and $\breve{\bm{a}}=\partial_{tt}\breve{\bm{d}}$. In these equations, we have made use of the fact that tensors $\bm{s}_h$ and $\bm{\xi}_{\textrm{sl},h}$ are symmetric, and it is understood that all terms are evaluated at $t^{n+1}$, the variables with a tilde being guesses to the unknowns (from a previous iteration step, for example). Observe that the problem is linear in $\delta \bm{y}_{\textrm{sl},h}$.

This version of the problem is unfeasible as it deals with gradients of the sub-scales for which we do not have an approximation. These terms are integrated by parts and derivatives transferred to the test functions. As explained for the abstract problem, the stabilization terms can then be written as (see \cref{bilinear_u_redef_gen}):
\begin{align}
\sum_K \langle \breve{\bm{y}}_{\textrm{sl},h} , \La^{*}_{\rm sl}(\tilde{\bm{y}}_{\textrm{sl}};\bm{z}_{\textrm{sl},h}) \rangle_K
+ \sum_K \langle \breve{\bm{y}}_{\textrm{sl}} , \Fa^{*}_{\rm sl}(\tilde{\bm{y}}_{\textrm{sl}};\bm{z}_{\textrm{sl},h}) \rangle_{\partial K}.
\label{eq:stab-terms}
\end{align} 
As explained before, we shall neglect the second term, although this approximation can be relaxed (see \cite{codina09}) and, in fact, is needed if stresses or pressures are discontinuous. Concerning the first term, operator $\La^{*}_{\rm sl}(\tilde{\bm{y}}_{\textrm{sl}};\bm{z}_{\textrm{sl},h}) $ has three components:
$$\La^{*}_{\rm sl}(\tilde{\bm{y}}_{\textrm{sl}};\bm{z}_{\textrm{sl},h}) = [\La_{{\rm sl},1}^{*}(\tilde{\bm{y}}_{\textrm{sl}};\bm{z}_{\textrm{sl},h}),\La_{{\rm sl},2}^{*}(\tilde{\bm{y}}_{\textrm{sl}};\bm{z}_{\textrm{sl},h}),\La_{{\rm sl},3}^{*}(\tilde{\bm{y}}_{\textrm{sl}};\bm{z}_{\textrm{sl},h})]^{T},$$ 
the first being a vector, the second a tensor and the third a scalar, given by:
\begin{subequations} \label{tr_adjoint}
\begin{align}
\La_{{\rm sl},1}^{*}(\tilde{\bm{y}}_{\textrm{sl}};\bm{z}_{\textrm{sl},h})\vert_{i} & = \dpar{}{X_K}\left[\xi_{\textrm{sl},h,ij}\tilde{F}_{jK}
     -\frac{1}{n_{\textrm{d}}}\xi_{\textrm{sl},h,ll}\tilde{F}_{iK}\right]  -\dpar{}{X_K}\left(\frac{\tilde{J}}{2\mu_{\textrm{sl}}}\xi_{h,lm}\tilde{s}_{{\rm sl},ml}\tilde{F}_{Ki}^{-1}\right) \nonumber\\
     &-\dpar{}{X_K}\left(\left(\frac{\tilde{J}\tilde{p}_{\textrm{sl},h}}{\lambda}-1\right)q_{\textrm{sl},h}\tilde{F}_{Ki}^{-1}\right)+\left(\frac{2\mu_{\textrm{sl}}}{\lambda n_{\textrm{d}}}\dpar{}{X_K}\left(q_{\textrm{sl},h}\tilde{F}_{iK}\right)\right), 
     \label{tr_momentum_adjoint}\\
\La_{{\rm sl},2}^{*}(\tilde{\bm{y}}_{\textrm{sl}};\bm{z}_{\textrm{sl},h}) \vert_{ij}
& = \dpar{e_{h,i}}{x_j} + \frac{\tilde{J}}{2\mu_{\textrm{sl}}}\xi_{\textrm{sl},h,ij}, \label{tr_stress_adjoint}\\
\La_{{\rm sl},3}^{*}(\tilde{\bm{y}}_{\textrm{sl}};\bm{z}_{\textrm{sl},h}) & = \dpar{e_{h,i}}{x_i}  + \frac{\tilde{J}}{\lambda}q_{\textrm{sl}}.\label{tr_press_adjoint}
\end{align}
\end{subequations}

To obtain the final expression of the stabilized problem, we need the expression of the sub-scales to be introduced in \cref{eq:stab-terms}, as well as the expression of $\langle \rho_{\rm sl} \breve{\bm{a}} ,  \bm{e}_{h}\rangle_K$ in \cref{eq:d1-a}. For that, we just have to apply the general expressions obtained in the previous section, taking into account that we only have time derivatives of the displacements. 

Using the same arguments as in \cite{Codina2009}, we will take the matrix of stabilization parameters within each element $K$ as
\begin{align}
& \bm{\tau}_{\textrm{sl},K} = {\rm diag} \left[
\tau_{{\rm sl},1,K} \bm{I}_{n_{\rm d}},
\tau_{{\rm sl},2}  \bm{I}_{n_{\rm d} \times n_{\rm d} } ,
\tau_{{\rm sl},3}  
\right], \label{eq:tau-sl1}\\
& \tau_{{\rm sl},1,K} := \left(  {c}_{\textrm{sl},1}\frac{\mu_{\textrm{sl}}}{h_{K}^{2}}  \right)^{-1},\quad
\tau_{{\rm sl},2} := {c}_{\textrm{sl},2}2\mu_{\textrm{sl}}, \quad
\tau_{{\rm sl},3} := {c}_{\textrm{sl},3}2\mu_{\textrm{sl}}, \label{eq:tau-sl2}
\end{align}
where $\bm{I}_{n_{\rm d} \times n_{\rm d} }$ is the identity on second order tensors and
${c}_{\textrm{sl},1} = 4.0$, ${c}_{\textrm{sl},2}=0.1$ and ${c}_{\textrm{sl},3}=0.1$ are numerical constants defined in the same way as in \cite{Castillo2014a} for linear elements, $h_K$ being divided by the polynomial order for higher order interpolations. 

Using \cref{eq:tau-sl1}, from expression (\ref{eqn:sgs_dyn}) we now have, within each element $K$:
\begin{subequations}\label{eq:sgs-sl}
\begin{align}
\breve{\bm{d}}^{n+1} & = \tau_{{\rm sl},1,t} 
\left( \bm{R}_{{\rm sl},1}^{n+1} + \frac{\rho_{\rm sl}\breve{\bm{d}}^{\theta}}{\delta t^2}\right),\qquad 
\tau_{{\rm sl},1,t} := \left( \frac{\rho_{\rm sl}\gamma_1}{\delta t^2} + \tau_{{\rm sl},1,K}^{-1}\right)^{-1}, 
\label{eq:sgs-d}\\
\breve{\bm{s}}_{\rm sl}^{n+1} & = \tau_{{\rm sl},2} \bm{R}_{{\rm sl},2}^{n+1},
\label{eq:sgs-ssl}\\
\breve{p}_{\rm sl}^{n+1} & = \tau_{{\rm sl},3} {R}_{{\rm sl},3}^{n+1},
\label{eq:sgs-psl}
\end{align}
\end{subequations}
where 
\begin{align}
\bm{R}^{n+1}_{\rm sl} = [ \bm{R}_{{\rm sl},1}^{n+1} , \bm{R}_{{\rm sl},2}^{n+1} , {R}_{{\rm sl},3}^{n+1} ]^T 
= \bm{F}_{\rm sl}^{n+1}  - [\rho_{\rm sl}\tilde{\bm{a}}_{\textrm{sl}}^{n+1} , {\bf 0} , 0]^T - \mathcal{A} (\tilde{\bm{y}}^{n+1}_{\textrm{sl}})
\end{align}
is the residual of the equation being solved at time step $n+1$, as it appears in \cref{eq:sol-lin-ibvp}.

From expression  (\ref{eqn:sgs_dyn_final}) we have:
\begin{align}
\langle \rho_{\rm sl} \breve{\bm{a}}^{n+1} , \bm{e}_h \rangle_K  = 
 \left\langle (1 - \tau_{{\rm sl},1,K}^{-1} \tau_{{\rm sl},1,t} )  \bm{R}_{{\rm sl},1}^{n+1} , \bm{e}_h \right\rangle_K 
-  \left\langle \tau_{{\rm sl},1,K}^{-1} \tau_{{\rm sl},1,t}   \frac{\rho_{\rm sl} \breve{\bm{d}}^\theta}{\delta t^2} , \bm{e}_h \right\rangle_K.
\label{eq:sgs-sl-ac}
\end{align}
It is seen from this expression that for the solid we have decided to approximate the acceleration of the sub-scale with the same kind of time integrator as for the FE scale, this is, a BDF2 scheme.

The fully discrete and stabilized problem is now completely defined. It is given by equations (\ref{eqn:split_nonli_trifil}), with the terms involving the subgrid scales given by the first term in \cref{eq:stab-terms}, the adjoint operator given in Eqs.~(\ref{tr_adjoint}), the matrix of stabilization parameters in \cref{eq:tau-sl1}, the sub-scales in Eqs.~(\ref{eq:sgs-sl}) and the acceleration of the sub-scales appearing in (\ref{eqn:split_nonli_trifil}) given in \cref{eq:sgs-sl-ac}.

\section{Three-field Navier-Stokes equations}\label{sec:nssection}

In this section we present a short review of the three-field Navier-Stokes equations and the FE approximation we use, as it is explained for example in \cite{Castillo2014a}.

\subsection{Governing equations}

Let $\Omega_{\textrm{fl}} \subset \mathbb{R}^{n_{\rm d}}$ be the domain where the fluid moves in the time interval $]0,t_{\textrm{f}}[$. We consider this fluid as incompressible and Newtonian, with viscosity $\mu_{\rm fl}$. If $\bm{u}$ is the fluid velocity, $p_{\rm fl}$ the pressure and $\bm{s}_{\rm fl}$ the deviatoric part of the stress tensor, the problem to be solved consists in finding 
$\bm{y}_{\rm fl} := [\bm{u} , \bm{s}_{\rm fl} , p_{\rm fl}]: \Omega_{\rm fl} \times ]0,t_{\rm f}[ \longrightarrow \mathbb{R}^{n_{\rm d}}\times  \mathbb{R}^{n_{\rm d}}\otimes \mathbb{R}^{n_{\rm d}}\times \mathbb{R}$  such that:
\begin{subequations}\label{nstrifilstrong}
\begin{flalign}
    &\quad & \rho_{\textrm{fl}}\partial_{t}\bm{u} - \nabla \cdot
    \bm{s}_{\textrm{fl}} + \rho_{\textrm{fl}}\bm{u}\cdot\nabla\bm{u} + \nabla p_{\textrm{fl}}      &= \rho_{\textrm{fl}}\bm{f}_{\textrm{fl}} & &\; \text{in} \; \Omega_{\textrm{fl}},    & t\in\;]0,t_{\textrm{f}}[,&\\ 
    &\quad &\frac{1}{2\mu_{\textrm{fl}}}\bm{s}_{\textrm{fl}} - \nabla^{\rm s}\bm{u} &= \bm{0}& &\;\text{in} \; \Omega_{\textrm{fl}},    & t\in\;]0,t_{\textrm{f}}[,&\\ 
    &\quad & \nabla\cdot\bm{u}  &= \bm{0}  & &\;\text{in} \; \Omega_{\textrm{fl}},    & t\in\;]0,t_{\textrm{f}}[,&\\ 
    &\quad & \bm{u}                   &= \bm{u}_{D}  &   &\;\text{on} \; \Gamma_{\textrm{fl},D},& t\in\;]0,t_{\textrm{f}}[,&\nonumber\\ 
    &\quad & \bm{n}_{\textrm{fl}}\cdot\bm{\sigma} &= \bm{t}_{\rm fl}    &    &\;\text{on} \; \Gamma_{\textrm{fl},N},& t\in\;]0,t_{\textrm{f}}[,&  \nonumber\\
    &\quad & \bm{u}   &= \bm{u}^{0}     &   &\;\text{in} \; \Omega_{\textrm{fl}},    & t=0,\;\;\;\;\;\; & \nonumber
\end{flalign}
\end{subequations}
where $ \nabla^{\rm s}\bm{u} $ is the symmetrical part of the velocity gradient, $\rho_{\textrm{fl}}$ is the fluid density, $\bm{f}_{\textrm{fl}}$ is the force vector, $\bm{u}^{0}$ is a prescribed initial velocity, $\bm{u}_D$ is a prescribed velocity on the boundary $\Gamma_{\textrm{fl},D}$, $\bm{t}_{\rm fl}$ is a prescribed traction on the boundary $\Gamma_{\textrm{fl},N}$, and $\bm{n}_{\textrm{fl}}$ is the normal to the boundary of the fluid domain. As for the solid, $\Gamma_{\textrm{fl},D}$ and $\Gamma_{\textrm{fl},N}$ are assumed to be a partition of $\partial\Omega_{\rm fl}$, but later we will introduce the interface with the solid. Note that the sign of $p_{\rm fl}$ is positive in compression, whereas the sign of $p_{\rm sl}$ is positive in traction.

\subsection{Weak form} 

We will use in the following a notation analogous to that of the solid. Thus, let $Y_{\rm fl}$ be the space where the unknown $\bm{y}_{\rm fl}$ must belong for each time $t$, satisfying the Dirichlet boundary conditions, and let $Y_{{\rm fl},0}$ be the space of functions satisfying the homogeneous counterpart of these conditions. The weak form of the three-field Navier-Stokes problem consists in finding $\bm{y}_{\rm fl} : [0,t_{\rm f}] \longrightarrow Y_{\rm fl}$ solution of the variational problem:
\begin{subequations}\label{sys:nstrifilweak} 
\begin{align}
    & (\rho_{\textrm{fl}}\partial_{t}\bm{u},\bm{v}) 
    + (\bm{s}_{\textrm{fl}},\nabla^{\rm s}\bm{v}) 
    + \langle \rho_{\textrm{fl}}\bm{u}\cdot\nabla\bm{u},\bm{v}\rangle \nonumber\\ 
    & \qquad - (p_{\textrm{fl}},\nabla\cdot\bm{v})  = \langle \rho_{\textrm{fl}}\bm{f}_{\textrm{fl}},\bm{v}\rangle + \langle \bm{t}_{\textrm{fl}} ,\bm{v}\rangle_{\Gamma_{\textrm{fl},N}} &  &\text{in} \; \Omega_{\textrm{fl}},~ t\in\;]0,t_{\textrm{f}}[,\label{eq:ns-1}\\ 
    & \frac{1}{2\mu_{\textrm{fl}}} ( \bm{s}_{\textrm{fl}},\bm{\xi}_{\textrm{fl}}) - (\nabla^{\rm s}\bm{u},\bm{\xi}_{\textrm{fl}}) = {0} & &\text{in} \; \Omega_{\textrm{fl}},~t\in\;]0,t_{\textrm{f}}[,\label{eq:ns-2}\\ 
    & (\nabla\cdot\bm{u},q_{\textrm{fl}})  = {0} &  &\text{in} \; \Omega_{\textrm{fl}},t\in\;]0,t_{\textrm{f}}[,\label{eq:ns-3}
\end{align}
\end{subequations}
for all ${\bm{z}_{\rm fl}} := [\bm{v},\bm{\xi}_{\textrm{fl}},q_{\textrm{fl}}] \in Y_{{\rm fl},0}$, and satisfying the initial conditions in a weak sense. We may now define a form $B_{\textrm{fl}}$ as:
\begin{align}
B_{\textrm{fl}}(\hat{\bm{u}} ; \bm{y}_{\textrm{fl}},\bm{z}_{\textrm{fl}}) 
& =  (\bm{s}_{\textrm{fl}},\nabla^{\rm s}\bm{v}) + \langle\rho_{\textrm{fl}}\hat{\bm{u}}\cdot\nabla\bm{u},\bm{v}\rangle - (p_{\textrm{fl}},\nabla\cdot\bm{v})   \nonumber\\ 
& + (\nabla\cdot\bm{u},q_{\textrm{fl}}) + \frac{1}{2\mu_{\textrm{fl}}} (\bm{s}_{\textrm{fl}},\bm{\xi}_{\textrm{fl}}) - (\nabla^{\rm s}\bm{u},\bm{\xi}_{\textrm{fl}}),\label{eq:ns_bilinear}  
\end{align}
and a form,
\begin{equation}\label{eq:ns_linear}  
\begin{aligned}
    L_{\textrm{fl}}(\bm{z}_{\textrm{fl}}) =   \langle \rho_{\textrm{fl}}\bm{f}_{\textrm{fl}},\bm{v}\rangle + \langle \bm{t}_{\textrm{fl}} ,\bm{v}\rangle_{\Gamma_{\textrm{fl},N}},
\end{aligned}
\end{equation}
which enable us to write \cref{sys:nstrifilweak} as:
\begin{equation} \label{eqn:ns_trifil_bilinear_compressed}  
(\rho_{\textrm{fl}}\partial_{t}\bm{u},\bm{v}) +  B_{\textrm{fl}}({\bm{u}} ; \bm{y}_{\textrm{fl}},\bm{z}_{\textrm{fl}}) = L_{\textrm{fl}}(\bm{z}_{\textrm{fl}})\qquad \forall \bm{z}_{\textrm{fl}} \in Y_{{\rm fl},0}, 
\end{equation}
with the initial conditions holding in a weak sense.

\subsection{Time discretization}

For the temporal discretization, we consider classical finite difference schemes. In particular, we have used the second order Backward Difference (BDF2) scheme in the applications, which has the following form:
\begin{equation*} 
    \left.\frac{\partial \bm{u}_{h}}{\partial t} \right\vert_{t^{n+1}} \approx \frac{3 \bm{u}_{h}^{n+1}- 4 \bm{u}_{h}^{n} + \bm{u}_{h}^{n-1}}{\delta t} =: \delta_{t} \bm{u}_{h}^{n+1},
\end{equation*}
with the notation inherited from the problem for the solid.

\subsection{Galerkin spatial discretization}\label{nssectiondiscrete}

Let us consider now a FE partition of the fluid domain, for which we will use the same notation as for the solid domain. From this we may construct FE spaces $Y_{{\rm fl},h} \subset Y_{{\rm fl}}$ and $Y_{{\rm fl},0,h} \subset Y_{{\rm fl},0}$ in the usual manner. Again, only conforming approximations will be considered. 

The spatial discretization for the Navier-Stokes problem can be written as: find $\bm{y}_{\textrm{fl},h}  : [0,t_{\rm f}] \longrightarrow Y_{{\rm fl},h} $ as the solution to the problem:
\begin{align} 
    (\rho_{\textrm{fl}}\partial_{t}\bm{u}_h,\bm{v}_h) + B_{\textrm{fl}}(\bm{u}_h;\bm{y}_{\textrm{fl},h},\bm{z}_{\textrm{fl},h}) & = L_{\textrm{fl}}(\bm{z}_{\textrm{fl},h}),\label{disc_NS} 
\end{align}
for all $\bm{z}_{\textrm{fl},h} \in Y_{{\rm fl},0,h}$, and satisfying the initial conditions weakly. As it is well known, this problem is unstable for two main reasons: the need of having compatible velocity-pressure and stress-velocity interpolations, similar to those of the solid problem, and the possibility to encounter convection dominated flows. Both are overcome using the following stabilized formulation.

\subsection{VMS stabilization}\label{sec:ns_stab_trifil}

For the fluid problem we employ the same VMS strategy as for the solid. Now, instead of linearizing the equations using a Newton-Raphson strategy and solving for the incremental unknowns at each iteration, we employ a fixed-point strategy, assuming known the transport velocity of the convective term, $\tilde{\bm{u}}_h$. Thus, let us introduce the linearized Navier-Stokes operator and its formal adjoint:
\begin{align}
{\cal L}_{\rm fl} (\tilde{\bm{u}}_h ; \bm{y}_{\textrm{fl},h}) & 
:= \left[ \begin{matrix} 
{\cal L}_{{\rm fl},1} (\tilde{\bm{u}}_h ; \bm{y}_{\textrm{fl},h}) \\
{\cal L}_{{\rm fl},2} (\bm{y}_{\textrm{fl},h}) \\
{\cal L}_{{\rm fl},3} (\bm{y}_{\textrm{fl},h}) \end{matrix}\right]
= \left[ \begin{matrix}
- \nabla\cdot\bm{s}_{{\rm fl},h} + \nabla p_{\textrm{fl},h} + \rho_{\textrm{fl}}\tilde{\bm{u}}_h\cdot\nabla\bm{u}_{h} \\ 
 \frac{1}{2\mu_{\textrm{fl}}}\bm{s}_{h} - \nabla^{\rm s}\bm{v}_{h} \\
 \nabla\cdot\bm{u}_{h}
\end{matrix}\right], \nonumber\\
{\cal L}_{\rm fl}^\ast (\tilde{\bm{u}}_h ; \bm{z}_{\textrm{fl},h}) & 
:= \left[ \begin{matrix} 
{\cal L}_{{\rm fl},1}^\ast (\tilde{\bm{u}}_h ; \bm{z}_{\textrm{fl},h}) \\
{\cal L}_{{\rm fl},2}^\ast (\bm{z}_{\textrm{fl},h}) \\
{\cal L}_{{\rm fl},3}^\ast (\bm{z}_{\textrm{fl},h}) \end{matrix}\right]
= \left[ \begin{matrix}
 \nabla\cdot\bm{\xi}_{{\rm fl},h} - \nabla q_{\textrm{fl},h} - \rho_{\textrm{fl}}\tilde{\bm{u}}_h\cdot\nabla\bm{v}_{h} \\ 
 \frac{1}{2\mu_{\textrm{fl}}}\bm{\xi}_{h} + \nabla^{\rm s}\bm{v}_{h} \\
 -\nabla\cdot\bm{v}_{h}
\end{matrix}\right]. \nonumber
\end{align}
We define the matrix of stabilization parameters within each element $\bm{\tau}_{\textrm{fl},K}$ as:
\begin{align}
& \bm{\tau}_{\textrm{fl},K} = {\rm diag} \left[
\tau_{{\rm fl},1,K} \bm{I}_{n_{\rm d}},
\tau_{{\rm fl},2}  \bm{I}_{n_{\rm d} \times n_{\rm d} } ,
\tau_{{\rm fl},3}  
\right], \label{eq:tau-fl1}\\
& \tau_{{\rm fl},1,K} := \left(  {c}_{\textrm{fl},1}\frac{\mu_{\textrm{fl}}}{h_{K}^{2}}  
+ {c}_{\textrm{fl},2}\frac{\rho_{\textrm{fl}\vert \tilde{\bm{u}}_{h}\vert_K }}{h_{K}} \right)^{-1},\quad
\tau_{{\rm fl},2} := {c}_{\textrm{fl},3}2\mu_{\textrm{sl}}, \quad
\tau_{{\rm fl},3} := {c}_{\textrm{fl},4}2\mu_{\textrm{sl}}, \label{eq:tau-fl2}
\end{align}
where, when using linear elements, ${c}_{\textrm{fl},1} = 4.0$, ${c}_{\textrm{fl},2} = 1.0$, ${c}_{\textrm{fl},3} = 0.1$ and ${c}_{\textrm{fl},4} = 0.1$ are numerical constants, as defined in \cite{Castillo2014a}. For higher order interpolations the element size $h_K$ needs to be divided by the polynomial order. In \cref{eq:tau-fl2}, $\vert \tilde{\bm{u}}_{h}\vert_K $ can be taken as the maximum of the Euclidian norm of $\bm{u}_{h}$ in element $K$.

With these ingredients we can already write down the stabilized VMS-based formulation we employ for the fluid. Considering all terms evaluated at time $t^{n+1}$, it consists of finding $\bm{y}_{{\rm fl},h} =  [{\bm u}_h, {\bm s}_{{\rm fl},h},{p}_{{\rm fl}}]\in Y_{{\rm fl},h}$ such that:
\begin{align}
& (\rho_{\textrm{fl}}\delta_t {\bm u}_h , {\bm v}_h) + \sum_K \langle \rho_{\textrm{fl}} \delta_t \breve{\bm u} , {\bm v}_h \rangle_K
+ B_{\textrm{fl}}(\bm{u}_h;\bm{y}_{\textrm{fl},h},\bm{z}_{\textrm{fl},h}) \nonumber\\
& \qquad + \sum_K \langle \breve{\bm{y}}_{\textrm{fl}} , {\cal L}_{\rm fl}^\ast ({\bm{u}}_h ; \bm{z}_{\textrm{fl},h})\rangle_K
 = L_{\textrm{fl}}(\bm{z}_{\textrm{fl},h}),
\end{align}
for all $\bm{z}_{\textrm{fl},h}\in Y_{{\rm fl},0,h}$, where $\breve{\bm{y}}_{\textrm{fl}} = [\breve{\bm u}, \breve{\bm s}_{\rm fl},\breve{p}_{\rm fl}]$ is obtained from:
\begin{align*}
\rho_{\textrm{fl}} \delta_t \breve{\bm u}  + \tau_{{\rm fl},1,K}^{-1} \breve{\bm u} 
& = \rho_{\textrm{fl}}\bm{f}_{\textrm{fl}} - {\cal L}_{{\rm fl},1} ({\bm{u}}_h ; \bm{y}_{\textrm{fl},h}), \\
\breve{\bm s}_{\rm fl} & = - \tau_{{\rm fl},2} {\cal L}_{{\rm fl},2} (\bm{y}_{\textrm{fl},h}), \\
\breve{p}_{\rm fl} & = - \tau_{{\rm fl},3}{\cal L}_{{\rm fl},3} (\bm{y}_{\textrm{fl},h}).
\end{align*}

\section{Three-field Fluid-Structure Interaction}\label{fsisection}

Once the problems for solid and for the fluid and their approximation have been described, we may proceed to write the FSI problem. In what follows, we shall assume that there is common moving boundary $\Gamma_I(t)$ between the solid and the fluid, so that the boundary of the solid domain $\Omega_{\rm sl}(t)$ is $\partial\Omega_{\rm sl} = \Gamma_{{\rm sl},D} \cup \Gamma_{{\rm sl},N}  \cup \Gamma_I$, whereas the boundary of the fluid domain $\Omega_{\rm fl}(t)$ is $\partial\Omega_{\rm fl} = \Gamma_{{\rm fl},D} \cup \Gamma_{{\rm fl},N}  \cup \Gamma_I$, with void intersection between boundary components. Again, subscript $D$ refers to Dirichlet boundary conditions and subscript $N$ to Neumann boundary conditions, the same as described in previous sections.

\subsection{Problem setting}

To write the equations to be solved, we consider the problem discretized in time, but still continuous in space, to simplify the writing. The FE approximation can be done as explained in the previous sections. We will comment on specific aspects of the FE approximation in the next subsection.

The solid equations are those introduced in Section~\ref{sec:sldsection}, written in an updated Lagrangian reference. However, to cope with the time dependency of the fluid domain we need to slightly modify the equations for the fluid introduced in Section~\ref{sec:nssection}. To this end we use the Arbitrary Lagrangian Eulerian (ALE) approach (see for example \cite{Donea1999}). Let $\bm{u}_{\rm dom}$ be the velocity assigned to the points of the fluid domain, which needs to match the velocity of $\partial\Omega_{\rm fl}$, i.e., to match the velocity of the moving boundary $\Gamma_I$ and vanish on the rest of $\partial\Omega_{\rm fl}$. Using the ALE reference, the only modification with respect to the purely Eulerian formulation is to replace the transport velocity $\bm{u}$ of the advective term by $\bm{c} := \bm{u} - \bm{u}_{\rm dom}$, so that this advective term becomes $\bm{c}\cdot\nabla \bm{u}$. If $\bm{c} = {\bf 0}$ we would obtain an updated Lagrangian formulation for the fluid as well.

The problem to be solved is the following:
\begin{itemize}
\item[] Loop over the number of time steps:
\begin{itemize}
\item[] At each time step, iterate until convergence, $(k)$ being the iteration counter:
\begin{itemize}
\item[$\bullet$] Solve the equations for the fluid, i.e., find $\bm{y}^{(k)}_{\textrm{fl}}\in Y_{\rm fl}$ such that:
\begin{align*} 
(\rho_{\textrm{fl}}\delta_{t}\bm{u}^{(k)},\bm{v}) +  B_{\textrm{fl}}({{\bm{u}}}^{(k-1)} - \bm{u}_{\rm dom} ; \bm{y}_{\textrm{fl}}^{(k)},\bm{z}_{\textrm{fl}}) = L_{\textrm{fl}}(\bm{z}_{\textrm{fl}})\qquad \forall \bm{z}_{\textrm{fl}} \in Y_{{\rm fl},0}.
\end{align*}
\item[$\bullet$] Solve the equations for the solid (written in terms of the iterative increments of the unknown), i.e., find $\delta\bm{y}^{(k)}_{\textrm{sl}} \in Y_{{\rm sl},0}$ such that:
\begin{align*} 
& (\rho_{\textrm{sl}}\delta \bm{a}^{(k)},\bm{e}) +  B_{\textrm{sl}}({\bm{y}}^{(k-1)}_{\textrm{sl}};\delta \bm{y}^{(k)}_{\textrm{sl}},\bm{z}_{\textrm{sl}}) \\
& \qquad = L_{\textrm{sl}}({\bm{y}}_{\textrm{sl}}^{(k-1)};\bm{z}_{\textrm{sl}})
- (\rho_{\textrm{sl}} {\bm{a}}^{(k-1)} , \bm{e}) \qquad \forall \bm{z}_{\textrm{sl}} \in Y_{{\rm sl},0}.
\end{align*}
\item[$\bullet$] Prescribe the transmission conditions:
\begin{align}
\bm{u}^{(k)} -  \delta_t\bm{d}^{(k)}  & = {\bf 0} && \quad \hbox{on}~\Gamma_I,\label{eq:fsi-dir}\\
\bm{n}_{\rm fl}\cdot (-p_{\rm fl}^{(k)} \bm{I}_{n_{\rm d}} + \bm{s}^{(k)}_{\rm fl}) + 
\bm{n}_{\rm sl}\cdot (p_{\rm sl}^{(k)} \bm{I}_{n_{\rm d}} + \bm{s}^{(k)}_{\rm sl}) & = {\bf 0} && \quad \hbox{on}~\Gamma_I.
\label{eq:fsi-neu}
\end{align}
\item[$\bullet$] Check convergence and update unknowns: ${\bm{y}}^{(k-1)}_{\textrm{fl}}\gets \bm{y}^{(k)}_{\textrm{fl}}$ and ${\bm{y}}^{(k-1)}_{\textrm{sl}} \gets {\bm{y}}^{(k-1)}_{\textrm{sl}} + \delta \bm{y}^{(k)}_{\textrm{sl}}$.
\end{itemize}
\item[] End iterative loop.
\end{itemize}
\item[] En loop over the number of time steps.
\end{itemize}
This is the monolithic version of the problem, in which all unknowns are solved at once, in a fully coupled way. It is understood that the initial conditions are prescribed at the first time step.

\subsection{Block-iterative coupling and comments on the fully discrete problem}

Rather than solving the monolithic version of the problem, the most popular approach is probably a block-iterative coupling, in which the solid and the fluid mechanics problems are solved sequentially. Using the also classical approach of prescribing Dirichlet conditions on $\Gamma_I$ coming from the solid when solving for the fluid, and Neumann conditions coming from the fluid when solving for the solid, the algorithm to solve the problem is:
\begin{itemize}
\item[] Loop over the number of time steps:
\begin{itemize}
\item[] At each time step, iterate until convergence, $(k)$ being the iteration counter:
\begin{itemize}
\item[$\bullet$] Solve the equations for the fluid, i.e., find $\bm{y}^{(k)}_{\textrm{fl}}\in Y_{\rm fl}$ such that:
\begin{align*} 
(\rho_{\textrm{fl}}\delta_{t}\bm{u}^{(k)},\bm{v}) +  B_{\textrm{fl}}({\bm{u}}^{(k-1)} - \bm{u}_{\rm dom} ; \bm{y}^{(k)}_{\textrm{fl}},\bm{z}_{\textrm{fl}}) = L_{\textrm{fl}}(\bm{z}_{\textrm{fl}})\qquad \forall \bm{z}_{\textrm{fl}} \in Y_{{\rm fl},0}.
\end{align*}
using the boundary condition:
\begin{align}
\bm{u}^{(k)} -  \delta_t\bm{d}^{(k-1)}  & = {\bf 0} && \quad \hbox{on}~\Gamma_I. \label{eq:dir-ale}
\end{align}
\item[$\bullet$] Solve the equations for the solid (written in terms of the iterative increments of the unknown), i.e., find $\delta\bm{y}^{(k)}_{\textrm{sl}} \in Y_{{\rm sl},0}$ such that:
\begin{align*} 
& (\rho_{\textrm{sl}}\delta \bm{a}^{(k)},\bm{e}) +  B_{\textrm{sl}}({\bm{y}}^{(k-1)}_{\textrm{sl}};\delta \bm{y}^{(k)}_{\textrm{sl}},\bm{z}_{\textrm{sl}}) \\
& \qquad = L_{\textrm{sl}}({\bm{y}}^{(k-1)}_{\textrm{sl}};\bm{z}_{\textrm{sl}})
- (\rho_{\textrm{sl}} {\bm{a}}^{(k-1)} , \bm{e}) \qquad \forall \bm{z}_{\textrm{sl}} \in Y_{{\rm sl},0}.
\end{align*}
using the boundary condition:
\begin{align*}
\bm{n}_{\rm fl}\cdot (-p^{(k)}_{\rm fl} \bm{I}_{n_{\rm d}} + \bm{s}^{(k)}_{\rm fl}) + 
\bm{n}_{\rm sl}\cdot (p^{(k)}_{\rm sl} \bm{I}_{n_{\rm d}} + \bm{s}^{(k)}_{\rm sl}) & = {\bf 0} && \quad \hbox{on}~\Gamma_I.
\end{align*}
\item[$\bullet$] Check convergence and update unknowns: ${\bm{y}}^{(k-1)}_{\textrm{fl}}\gets \bm{y}^{(k)}_{\textrm{fl}}$ and ${\bm{y}}^{(k-1)}_{\textrm{sl}} \gets {\bm{y}}^{(k-1)}_{\textrm{sl}} + \delta \bm{y}^{(k)}_{\textrm{sl}}$.
\end{itemize}
\item[] End iterative loop
\end{itemize}
\item[] En loop over the number of time steps
\end{itemize}
Note that we have used the current pressure and stress values in the fluid when solving for the solid, but also those of the previous iteration could have been employed.

The FE approximation in space of the equations to be solved follows directly from what has been explained in the paper. However, there are a few aspects of the space discretization that are particular of the FSI problem:

\begin{itemize}
\item[$\bullet$] The particular way the velocity $\bm{u}_{\rm dom}$ is computed depends strongly on the FE approximation, since what is needed in fact is the value of this velocity at the nodes of the FE mesh. The values at the nodes of $\Gamma_I$ are determined by imposing \cref{eq:dir-ale} for $\bm{u}_{\rm dom}$, implying that $\bm{c} = {\bf 0}$ on $\Gamma_I$, i.e., this boundary is a material one. To obtain $\bm{u}_{\rm dom}$ at the interior nodes of $\Omega_{\rm fl}$, the mesh movement algorithm
has been taken from \cite{Chiandussi2000}, which has proven simple, robust and reliable.
\item[$\bullet$] We have described the most classical Dirichlet-Neumann iteration-by-subdomain coupling, which suffices for the purposes of this paper to present a three-field approach for both the solid and the fluid. However, this coupling may suffer convergence difficulties when the solid is very soft or when the densities of the fluid and the solid are similar. In these situations, one might resort to more sophisticated coupling strategies (see for example \cite{codina-baiges-11}), including a Nitsche's-type method to prescribe conditions (\ref{eq:fsi-dir})-(\ref{eq:fsi-neu}). There are then numerous ways to segregate iteratively the calculation of $\bm{y}^{(k)}_{\textrm{fl}}$ and $\bm{y}^{(k)}_{\textrm{sl}}$ and to design an iteration-by-subdomain algorithm. However, we shall not pursue this analysis in this paper. 
\item[$\bullet$] A particularly relevant aspect of the three-field approach in FSI is the implementation of transmission conditions. Suppose first that the meshes for the solid and the fluid match at $\Gamma_I$, and that stresses from the fluid need to be transferred to the solid. Suppose also that a nodal Lagrangian interpolation is used. In a classical velocity-pressure approach, one should transmit normal stresses at the numerical integration points, which should be computed from pressures interpolated from the nodes to the integration points and velocity gradients at these same points. In a three-field approach, stresses are directly available at the nodes. If the meshes for the solid and the fluid do not match at $\Gamma_I$, there is an additional interpolation step. In a velocity-pressure approach, normal stresses at the integration points of the fluid have to transferred to the nodes of the solid, and from these to the integration points of the solid. In a three-field approach, normal stresses at the nodes of the fluid have to be transferred to the nodes of the solid; the number of operations involved in this case is significantly smaller than in the former. In both cases, the transfer of information can be used using the standard Lagrangian interpolation, as done in our implementation, or by imposing restrictions, as explained in \cite{Houzeaux2001}.
\end{itemize}

In the iterative process described, relaxation of the transmitted quantities is very often required if not mandatory. This allows one to minimize the number of block (fluid and solid) iterations. In this respect, we have used a relaxation of the position and velocity of the interface boundary that the solid solver transmits to the fluid solver. We denote this position as $\bm{d}_{\Gamma_{\rm I}}$; from it, one may compute the velocity of the fluid boundary and $\bm{u}_{\rm dom}$, as explained above. We have implemented an Aitken relaxation scheme, in particular Aitken $\Delta^{2}$, detailed in \cite{Kuttler2008}, which we describe now in our context. Suppose that from values at the $k$th iteration, the solid is solved, obtaining the boundary displacements $\bm{d}_{{\Gamma_{\rm I},{\rm s}}}^{(k+1)}$. Then, the fluid is solved from the boundary displacements $\bm{d}_{\Gamma_{\rm I}}^{(k+1)}$ computed as
\begin{align*}
\bm{d}_{\Gamma_{\rm I}}^{(k+1)} = \bm{d}_{\Gamma_{\rm I}}^{(k)} + \omega^{(k+1)} \bm{r}_{\Gamma_{\textrm{I}}}^{(k+1)},
\end{align*}
where
\begin{align*}
\bm{r}_{\Gamma_{\textrm{I}}}^{(k+1)}  := \bm{d}_{{\Gamma_{\rm I},{\rm s}}}^{(k+1)}- \bm{d}_{\Gamma_{\rm I}}^{(k)}, \qquad
\omega^{(k+1)}  = -\omega^{(k)} \frac{(\bm{r}_{\Gamma_{\textrm{I}}}^{(k)})^{T}(\bm{r}_{\Gamma_{\textrm{I}}}^{(k+1)}-\bm{r}_{\Gamma_{\textrm{I}}}^{(k)})}{|\bm{r}_{\Gamma_{\textrm{I}}}^{(k+1)}-\bm{r}_{\Gamma_{\textrm{I}}}^{(k)}|^{2}}.
\end{align*}

\section{Numerical Results}\label{numresection}

In this section, numerical results are shown for stationary and dynamic cases, first for the solid in order to benchmark the new formulation, and followed by some well known FSI benchmarks in order to compare with the traditional standard coupling. 

\subsection{Three-field elasto-dynamic benchmarking}

We will compare in what follows the behavior of the irreducible formulation, in which the only unknown is the displacement field, with the stabilized three-field formulation proposed in this paper, both in static and in dynamic cases. In all cases, equal continuous interpolations will be used for displacements, pressures and stresses in the case of the three-field formulation.

\subsubsection{Cook's membrane}\label{test1}

\begin{figure}[h!]
    \centering
    \includegraphics[width = 0.45\textwidth]{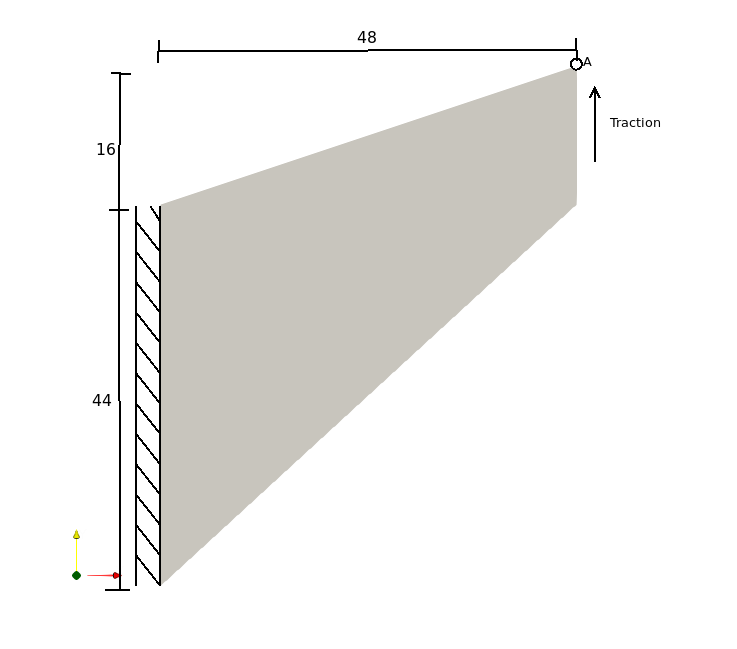}
    \caption{Geometry}
    \label{fig:1_cook_diag}
\end{figure}

The following example is a typical benchmark for solid mechanics. A tapered beam is subjected to a shearing load on one of its sides. In our case the shear traction is taken as $30$ GN. Fig.~\ref{fig:1_cook_diag} shows the geometry of the beam. 
 
The properties of the material for all tests in these benchmarks are shown in Table~\ref{1_physparam}.  For the irreducible formulation the Poisson coefficient was taken as $\nu=0.49$, so as to have a material with low compressibility but without locking. For the three-field formulation we deal with an incompressible material ($\nu = 0.5$).
 
\begin{table}[h!]  \caption{Physical parameters}
\centering
    \begin{tabular}{cc} 
        \hline                              &          \\ 
  \hline         $\rho_{\textrm{sl}}$     &  7850.0 [Kg/m$^{3}$] \\ 
                    $\mu_{\textrm{sl}} $    & 80$\times 10^9$  [Pa] \\
                    $\lambda_{\textrm{sl}}$  & $\infty$  \\
                                  Model     & Neo-Hookean \\
         \hline 
     \end{tabular}
     \label{1_physparam}
 \end{table}

\begin{figure}[h!]
    \centering
    \begin{subfigure}[h!]{0.45\textwidth}
        \includegraphics[width=\linewidth]{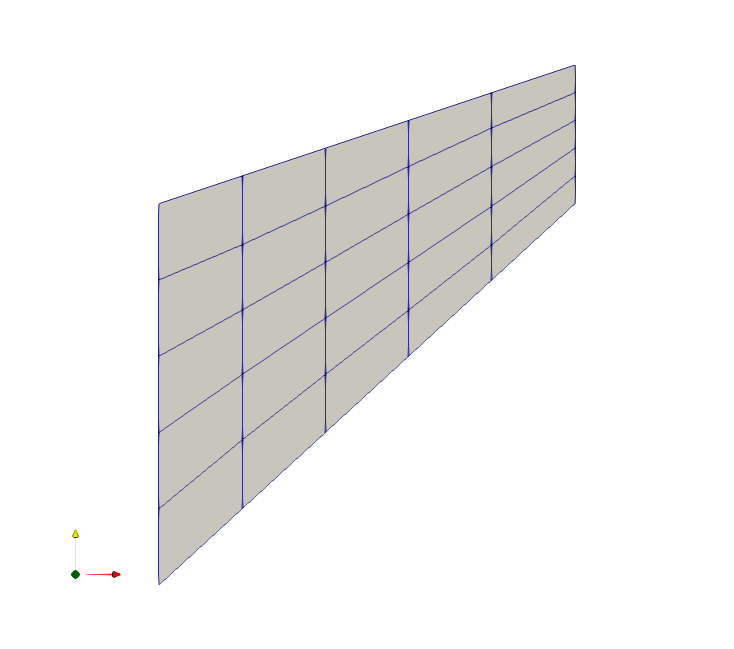}
        \caption{Quad mesh}
        \label{fig:1_quad_mesh} 
    \end{subfigure}
    \begin{subfigure}[h!]{0.45\textwidth}
        \includegraphics[width = \linewidth]{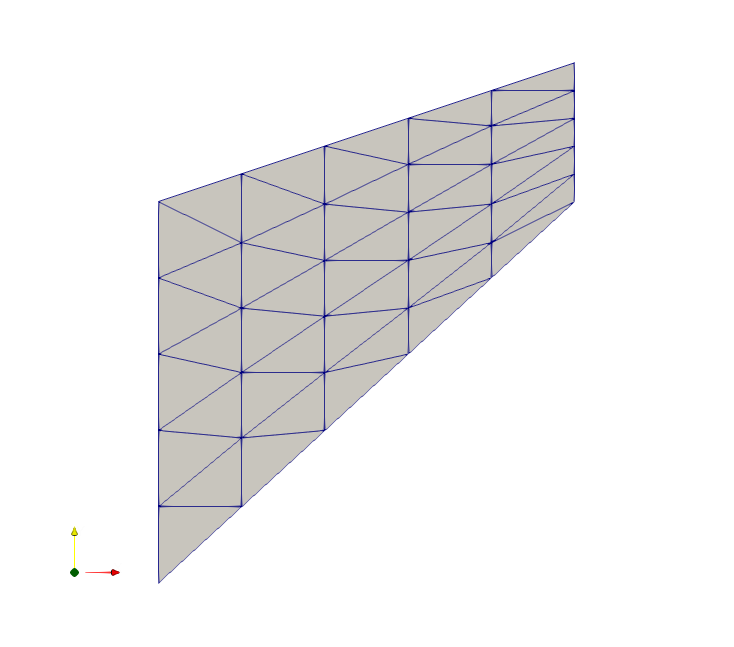}
        \caption{Triangle mesh}
        \label{fig:1_tri_mesh} 
    \end{subfigure}
    \caption{Mesh examples for benchmark}
    \label{fig:1_mesh} 
\end{figure}

Convergence tests were run for different mesh sizes. Figs.~\ref{fig:1_quad_mesh}  and \ref{fig:1_tri_mesh} show examples for quadrilateral (4 node bilinear) and triangular (3 node linear) elements.  The notation for the results is detailed in Table~\ref{1_caseparams}.
 
 \begin{table}[h!] \caption{Case parameters}
\centering
    \begin{tabular}{ccc} 
    \hline      Name        &  Formulation & Type of elem        \\ 
    \hline      irr\_tri     &  irreducible & linear triangle     \\ 
                irr\_sq      &  irreducible & bi-linear square    \\
                sup\_tri     &  three-field & linear triangle     \\
                sup\_sq      &  three-field & bi-linear square    \\ 
    \hline 
    \end{tabular} 
    \label{1_caseparams}
 \end{table}

\begin{figure}[h!]
    \centering
    \begin{subfigure}[h!]{0.45\textwidth}
        \includegraphics[width=\linewidth]{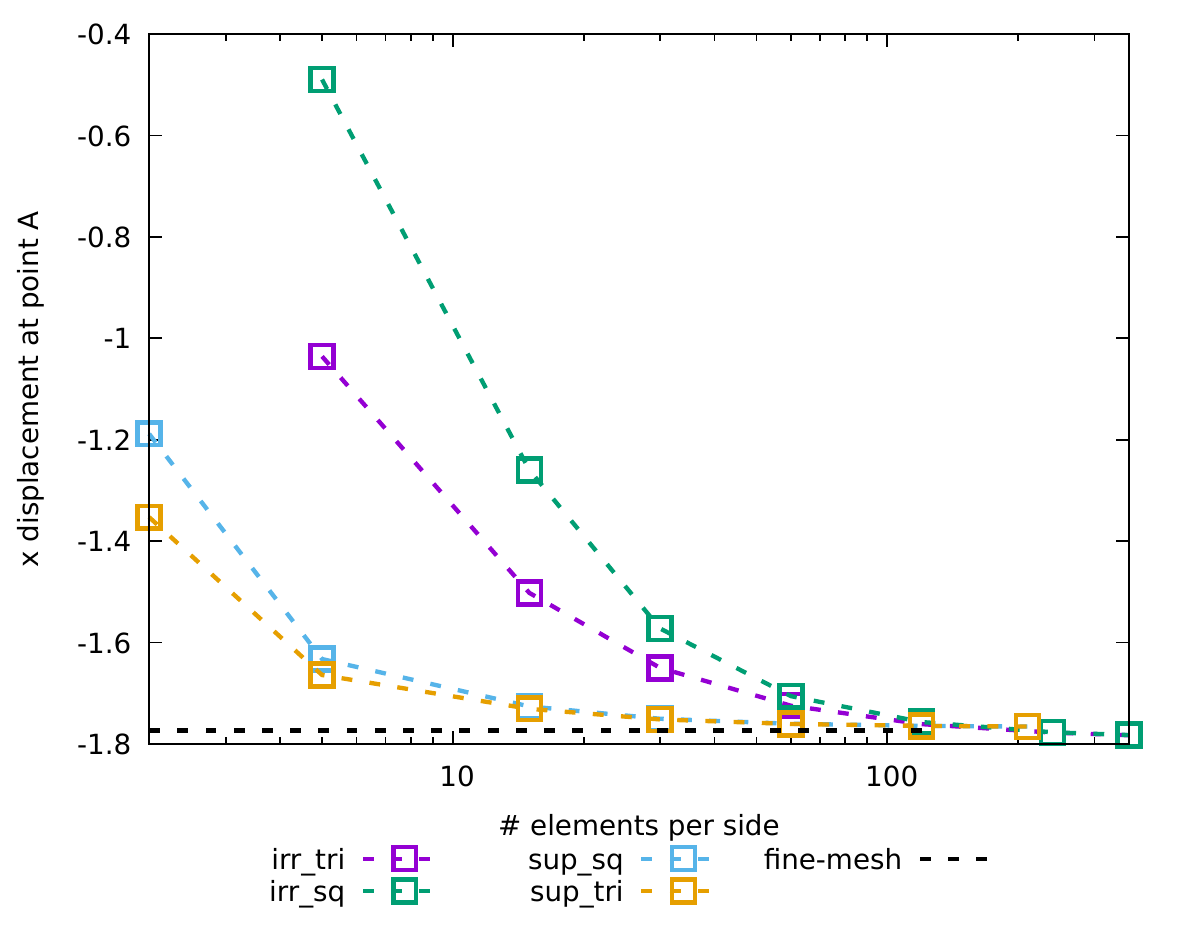}
        \caption{Displacement in $x$ axis}
        \label{fig:1_disp_x} 
    \end{subfigure}
    \begin{subfigure}[h!]{0.45\textwidth}
        \includegraphics[width = \linewidth]{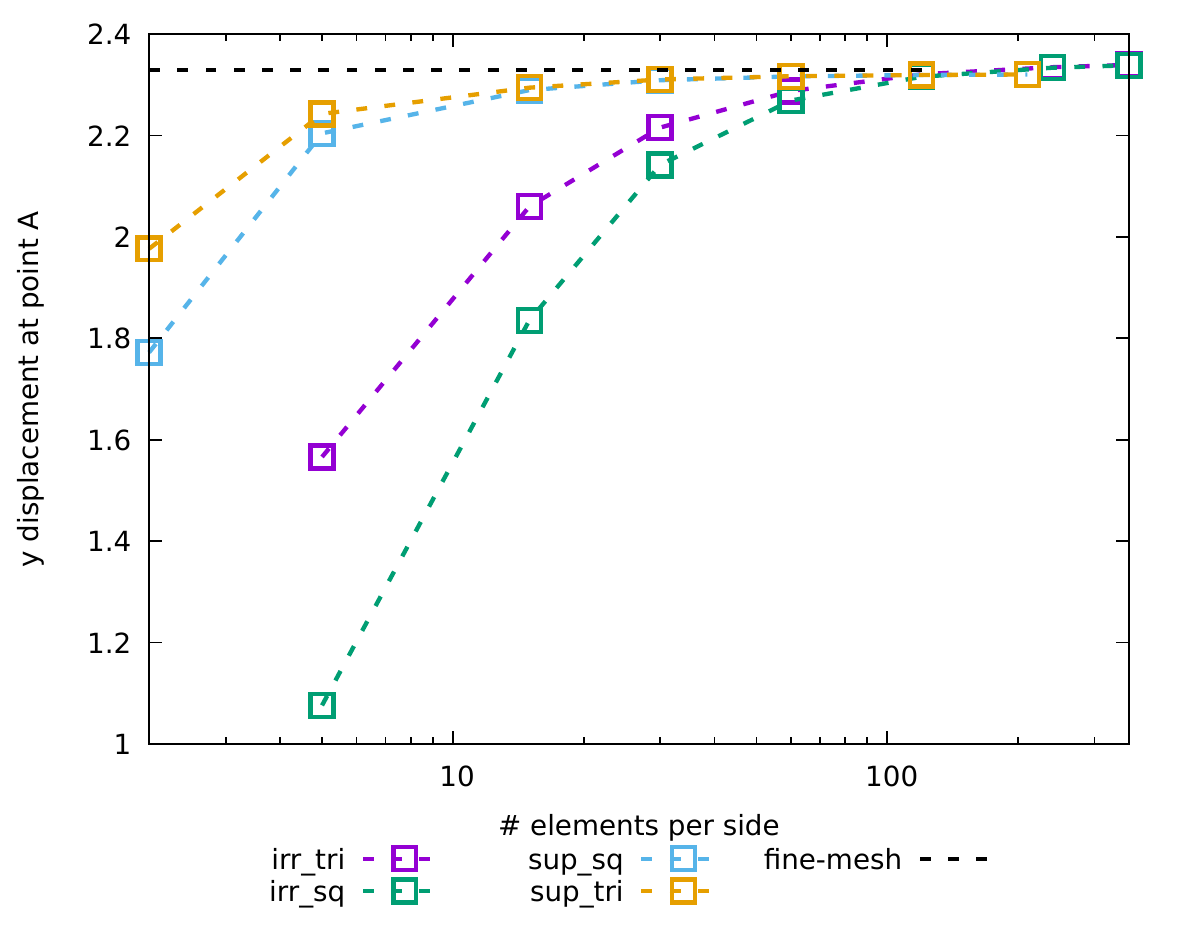}
        \caption{Displacement in $y$ axis}
        \label{fig:1_disp_y} 
    \end{subfigure}
    \caption{Displacement of point A}
    \label{fig:1_disp} 
\end{figure}

As it is a bending dominated test, it is of interest to see the displacement at the tip of the beam (point A) both in $x$ and $y$ directions. Fig.~\ref{fig:1_disp} shows the convergence for different mesh sizes in comparison with a solution obtained with a very fine mesh, shown in black.

\begin{figure}[h!]
    \centering
    \begin{subfigure}[h!]{0.45\textwidth}
        \includegraphics[width=\linewidth]{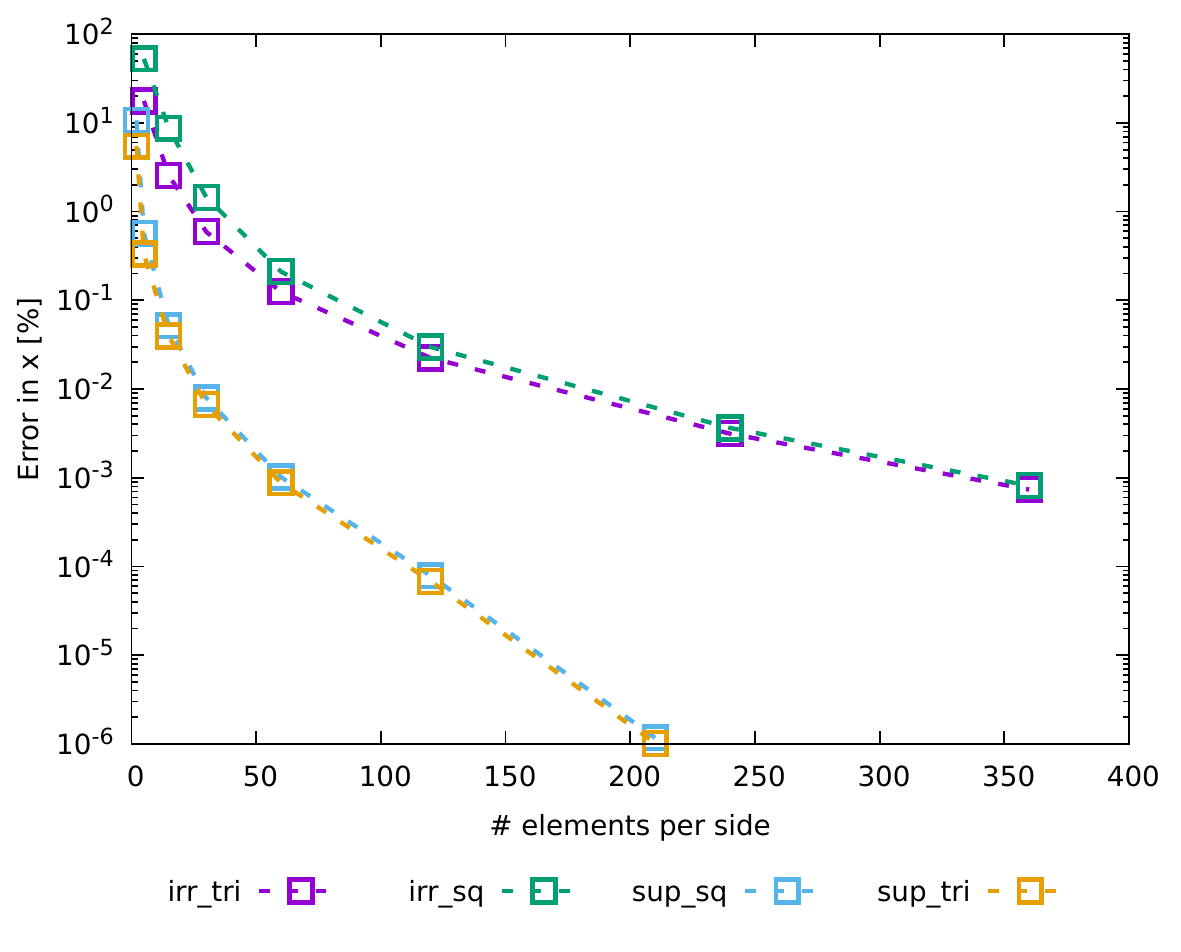}
        \caption{Error in $x$ axis}
        \label{fig:1_err_x} 
    \end{subfigure}
    \begin{subfigure}[h!]{0.45\textwidth}
        \includegraphics[width = \linewidth]{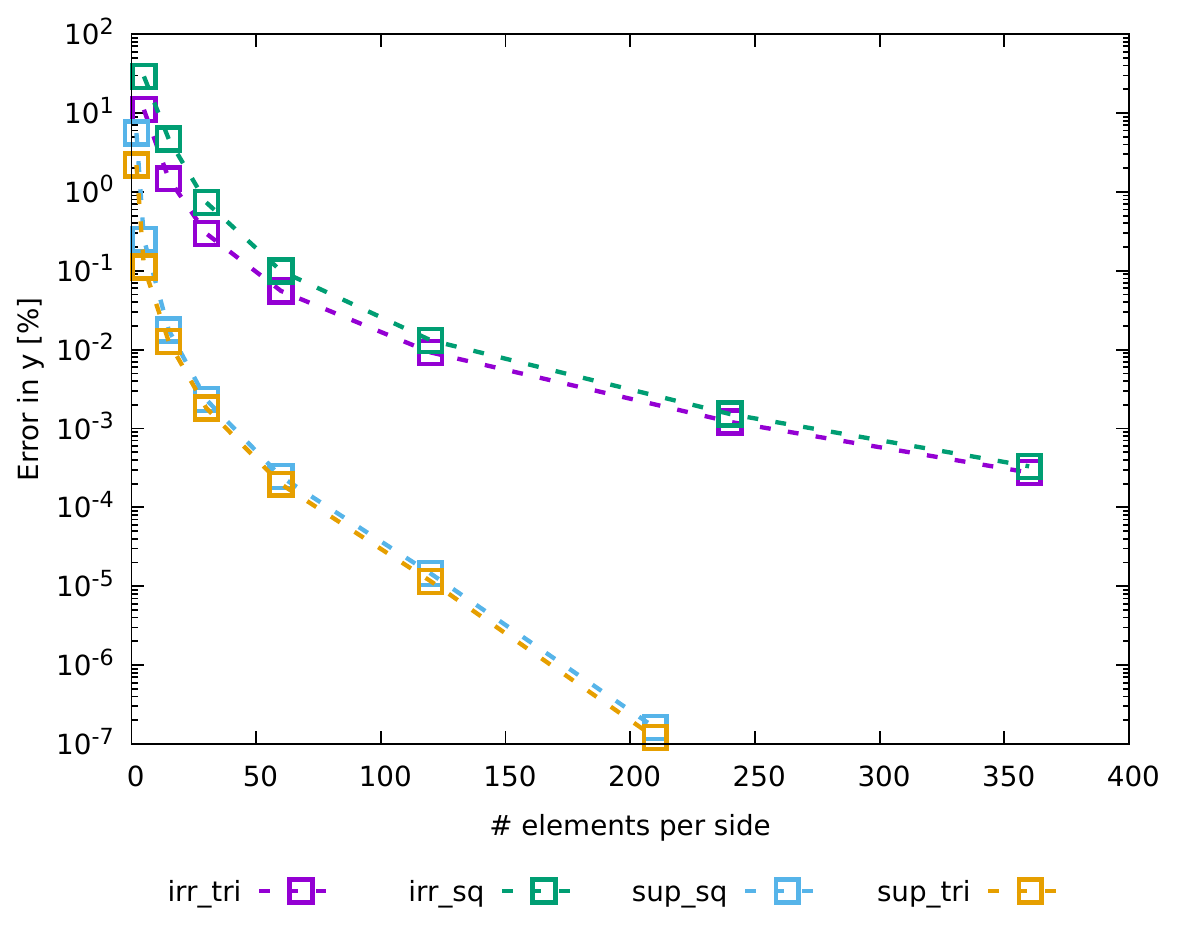}
        \caption{Error in $y$ axis}
        \label{fig:1_err_y} 
    \end{subfigure}
    \caption{Error for displacement at point A according to \# of elements}
    \label{fig:1_err} 
\end{figure}

\begin{figure}[h!]
    \centering
    \begin{subfigure}[h!]{0.45\textwidth}
        \includegraphics[width=\linewidth]{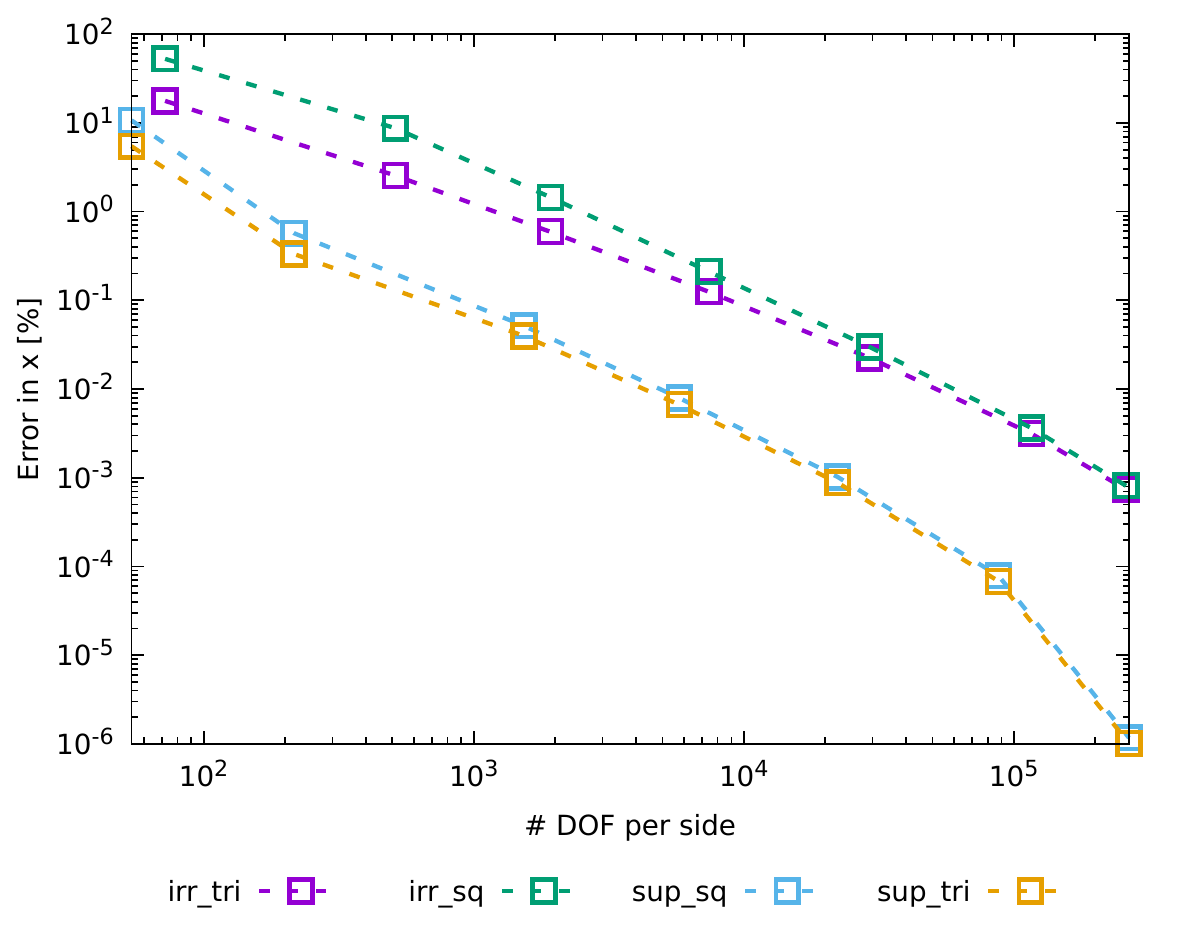}
        \caption{Error in $x$ axis}
        \label{fig:1_err_x_dof} 
    \end{subfigure}
    \begin{subfigure}[h!]{0.45\textwidth}
        \includegraphics[width = \linewidth]{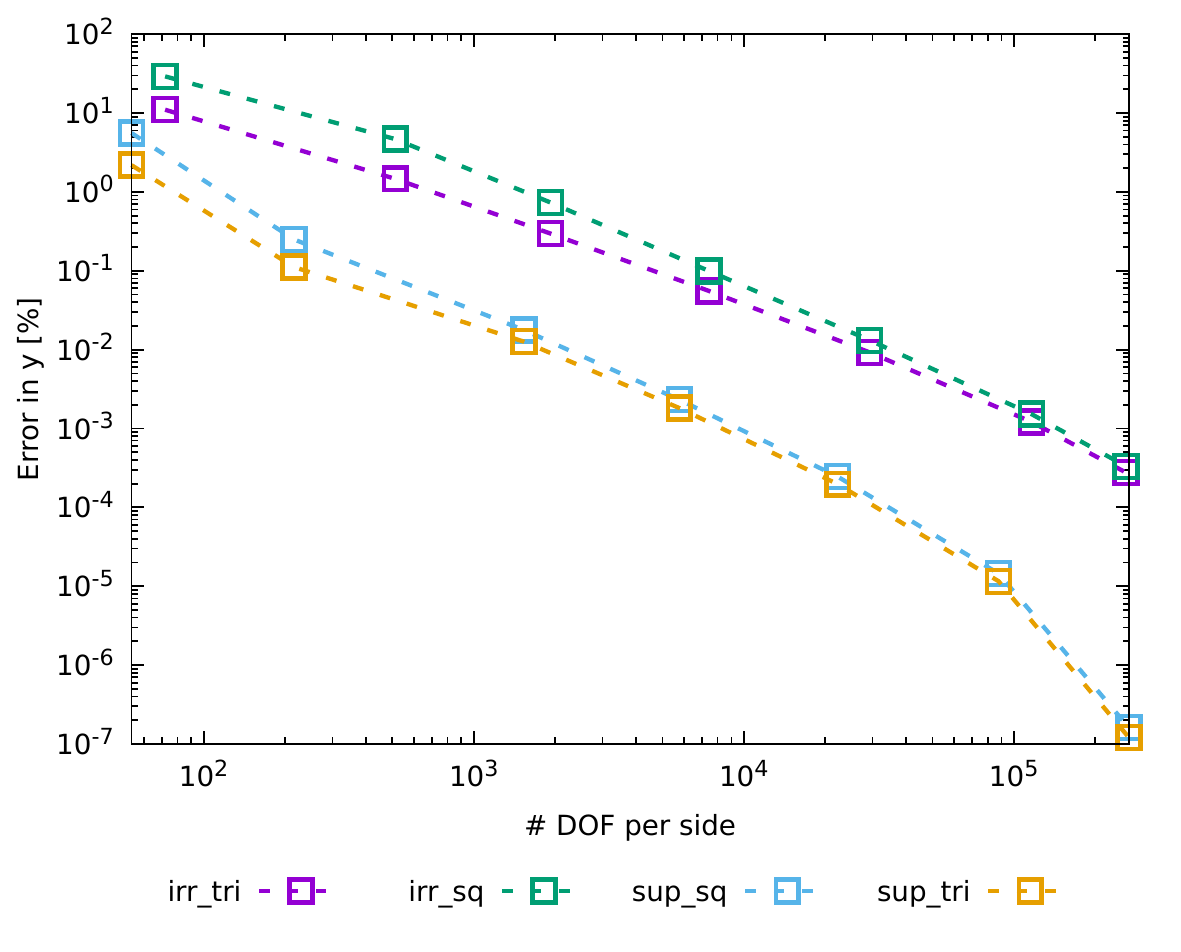}
        \caption{Error in $y$ axis}
        \label{fig:1_err_y_dof} 
    \end{subfigure}
    \caption{Error for displacement at point A according to \# of DOF}
    \label{fig:1_err_dof} 
\end{figure}

Fig.~\ref{fig:1_err} shows the evolution of the error for the results shown previously in terms of number of elements; it can be seen that the convergence for the three-field is better, although a fairer comparison can be seen from Fig.~\ref{fig:1_err_dof}, which shows error in terms of number of DOFs. Results are in agreement with \cite{Chiumenti2015}, keeping in mind that their results are shown for the linear case. Both the irreducible and the three-field formulation show good convergence upon mesh refinement, with the three-field model being more precise and having faster convergence overall. For the three-field formulation, bi-linear square elements and triangular elements show very similar convergence properties; however, in the irreducible case triangular elements appear to be more precise than squares.

\begin{figure}[h!]
    \centering
    \begin{subfigure}[h!]{0.45\textwidth}
        \includegraphics[width=\linewidth]{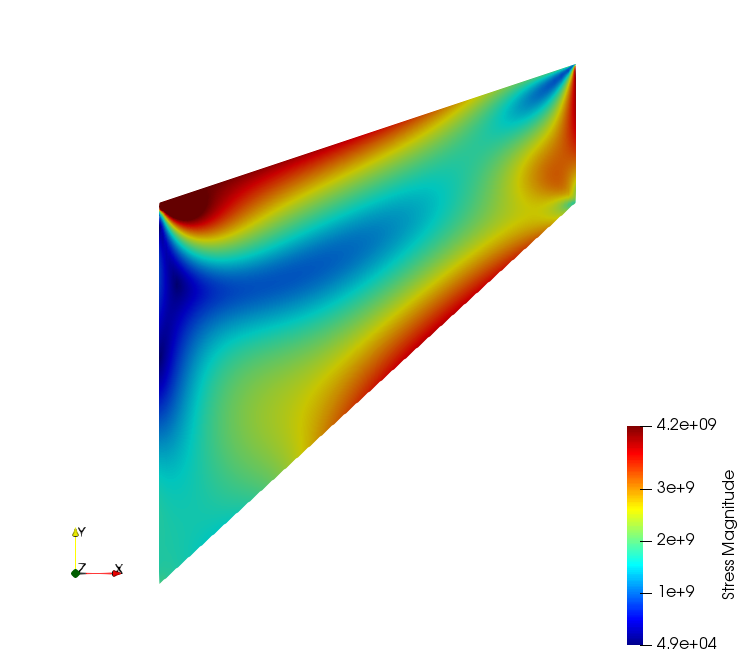}
        \caption{Deviatoric stress magnitude contours}
        \label{fig:1_stress} 
    \end{subfigure}
    \begin{subfigure}[h!]{0.45\textwidth}
        \includegraphics[width = \linewidth]{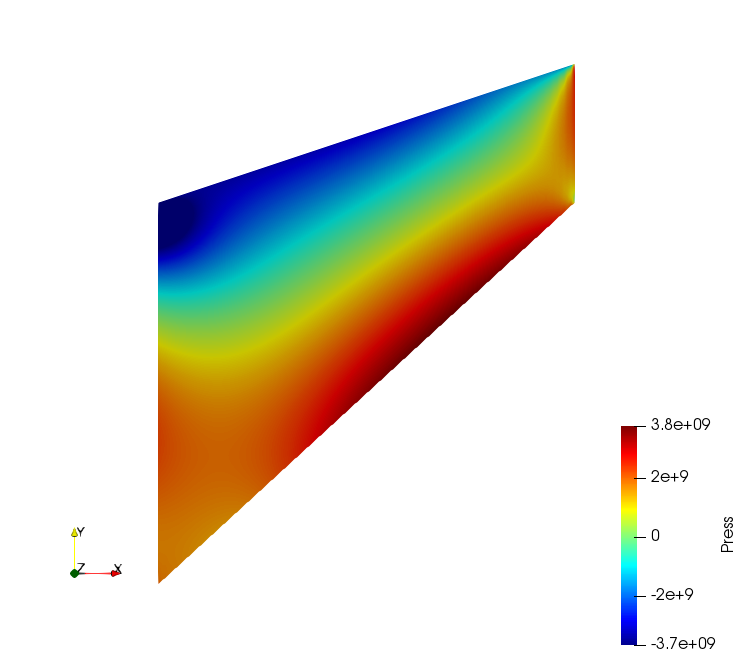}
        \caption{Pressure contours}
        \label{fig:1_press} 
    \end{subfigure}
    \caption{Deviatoric stress and pressure}
    \label{fig:1_stress_press} 
\end{figure}

Fig.~\ref{fig:1_stress_press} shows the stress and pressure distribution for the beam using the three-field formulation. It can be seen that smooth and continuous fields have been obtained, without any oscillation in spite of using equal interpolation for all fields.

\subsubsection{Dynamic oscillation of a cantilever bar}\label{test2}

\begin{figure}[h!]
    \centering
    \includegraphics[width = 0.55\textwidth]{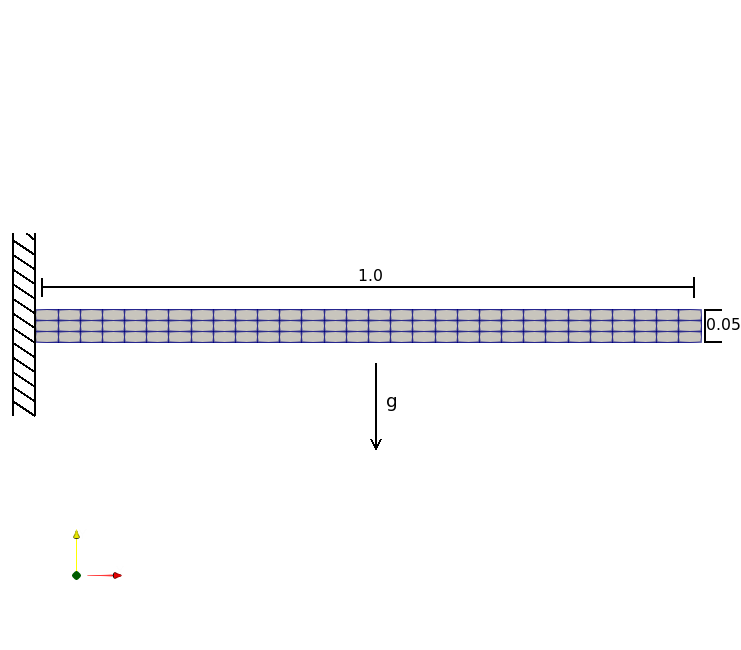}
    \caption{Geometry}
    \label{fig:2_beam_diag}
\end{figure}

In this section we analyze the time evolution of a clamped beam under the effect of gravity. As in the previous case, we compare our results against the same case for the irreducible formulation. Fig.~\ref{fig:2_beam_diag} shows the geometry of the beam and an example of a mesh used. For the initial conditions, the bar starts at rest and then suddenly gravity is applied, so the bar falls in the direction shown.

The properties of the material and parameters for all tests in this benchmark are shown in Table~\ref{2_physparam}. Note that for the irreducible formulation the Poisson coefficient was taken as $\nu=0.499$,  whereas for the three-field formulation we deal with an incompressible material. The time interval of analysis is $[0,1]$, with a time step $\delta t = 10^{-3}$.
 
\begin{table}[h!] \caption{Physical parameters}
\centering
    \begin{tabular}{cc} 
        \hline                              &          \\ 
  \hline            $\rho_{\textrm{sl}}$     & $100.0 ~[{{\rm Kg}}/{{\rm m}^{3}}]$ \\ 
                    $\mu_{\textrm{sl}} $     & $2.135 \times 10^{7} ~[{\rm Pa}]$  \\
                    $\lambda_{\textrm{sl}}$  & $\infty$  \\
                                  Model     & Neo-Hookean \\
                                  Gravity   & $2.0 [{{\rm m}}/{{\rm s}^2}]$ \\
         \hline 
     \end{tabular} 
     \label{2_physparam}
 \end{table}
 
\begin{figure}[h!]
    \centering
    \begin{subfigure}[h!]{0.45\textwidth}
        \includegraphics[width=\linewidth]{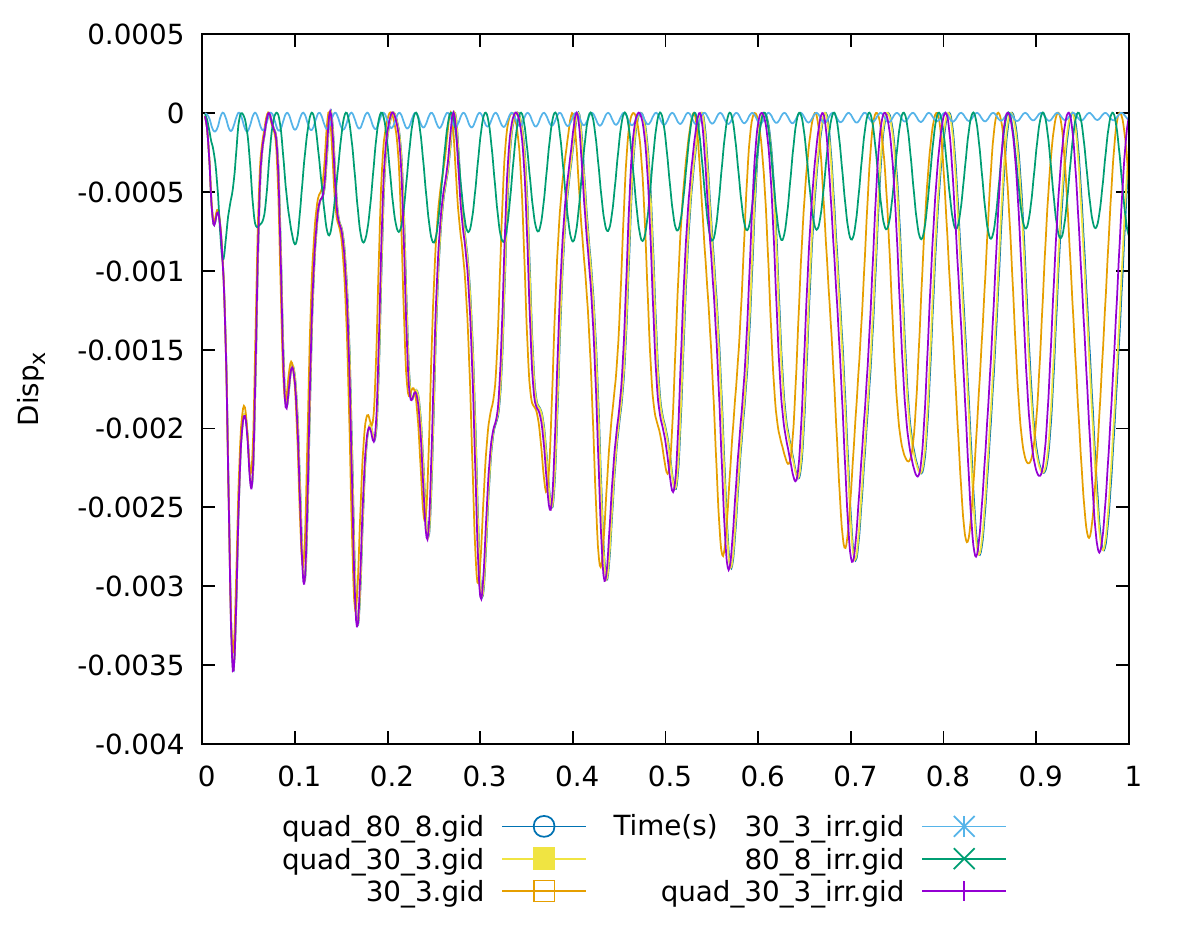}
        \caption{Displacement in $x$}
        \label{fig:2_disp_x_full} 
    \end{subfigure}
    \begin{subfigure}[h!]{0.45\textwidth}
        \includegraphics[width = \linewidth]{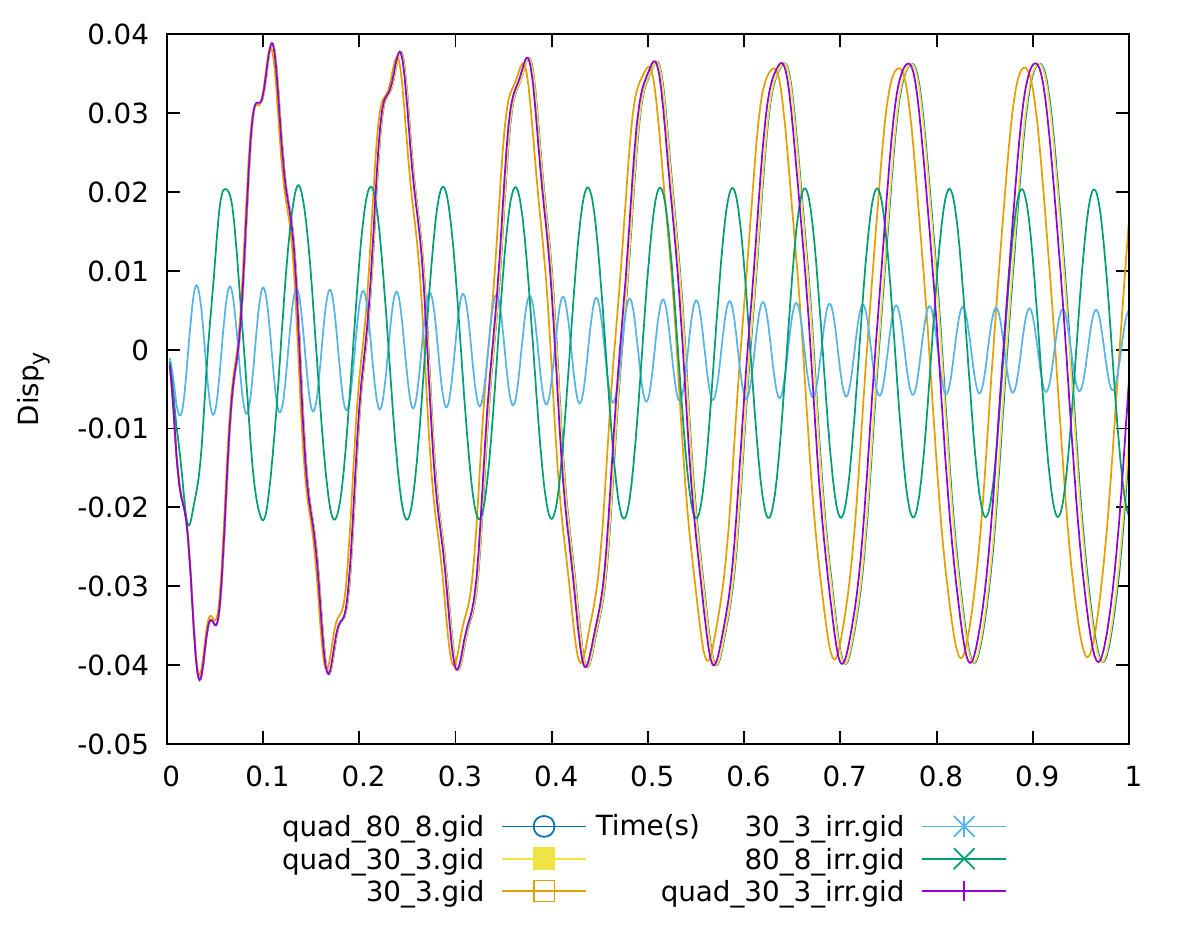}
        \caption{Displacement in $y$}
        \label{fig:2_disp_y_full} 
    \end{subfigure}
    \caption{Displacement at the tip of the beam}
    \label{fig:2_disp_full} 
\end{figure}

Fig.~\ref{fig:2_disp_full} shows the time evolution of the displacements for different cases run. Solutions were compared for both the three-field and the irreducible formulation for different mesh sizes and elements (linear or quadratic); only squares were used in this example. The notation for the results is detailed in Table~\ref{2_caseparams}.
 
 \begin{table}[h!] \caption{Case parameters}
\centering
    \begin{tabular}{ccccc} 
        \hline      Name             &  formulation & \# elems length wise & \# elems height & type of elem        \\ 
  \hline            30\_3            &  three-field &  30                 &        3       & linear square     \\ 
                    30\_3\_irr       &  irreducible &  30                 &        3       & linear square     \\
                    80\_8\_irr       &  irreducible &  80                 &        8       & linear square     \\
                    quad\_30\_3      &  three-field &  30                 &        3       & quadratic square     \\ 
                    quad\_30\_3\_irr &  irreducible &  30                 &        3       & quadratic square     \\
                    quad\_80\_8      &  three-field &  80                 &        8       & quadratic square     \\ 
         \hline 
     \end{tabular} 
     \label{2_caseparams}
 \end{table}
 
\begin{figure}[h!]
    \centering
    \begin{subfigure}[h!]{0.45\textwidth}
        \includegraphics[width=\linewidth]{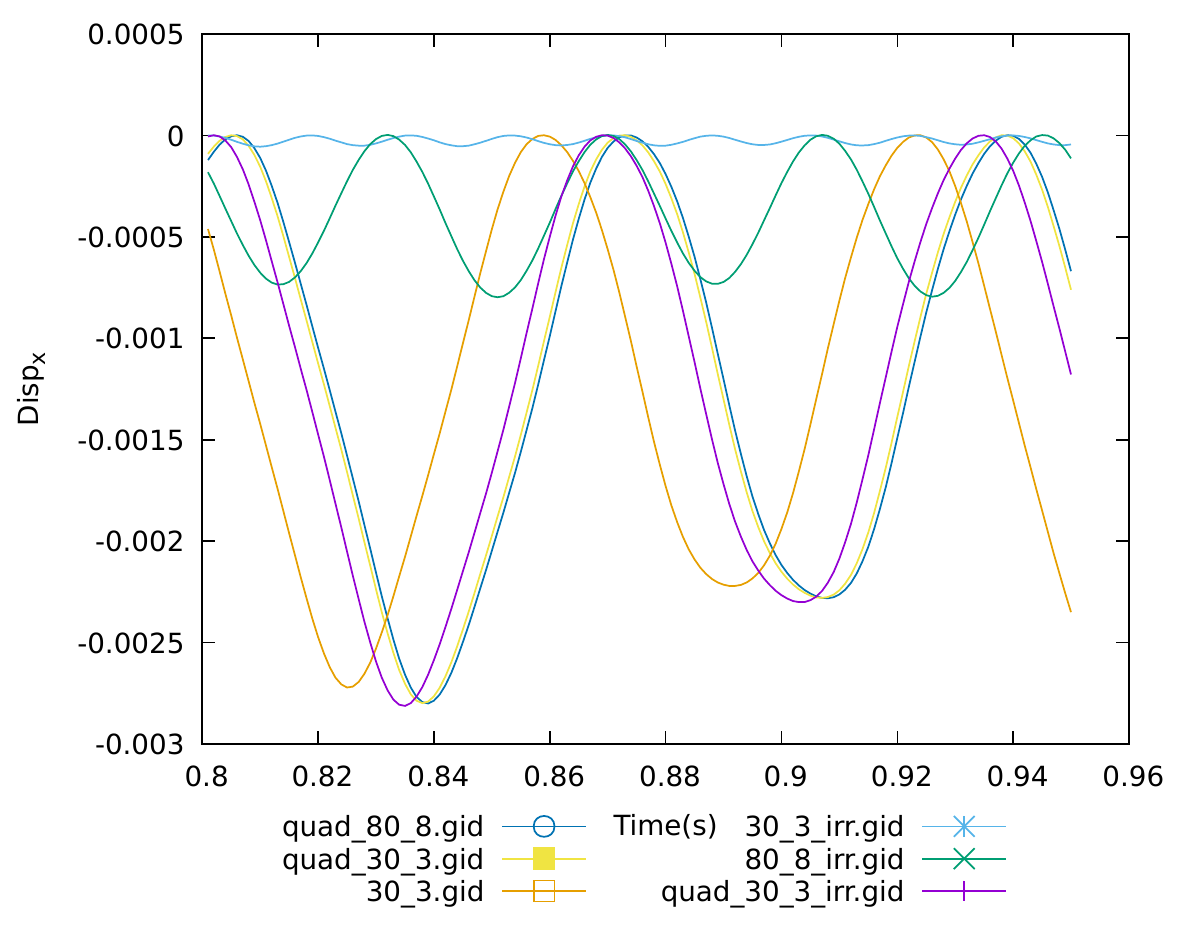}
        \caption{Displacement in $x$}
        \label{fig:2_disp_x} 
    \end{subfigure}
    \begin{subfigure}[h!]{0.45\textwidth}
        \includegraphics[width = \linewidth]{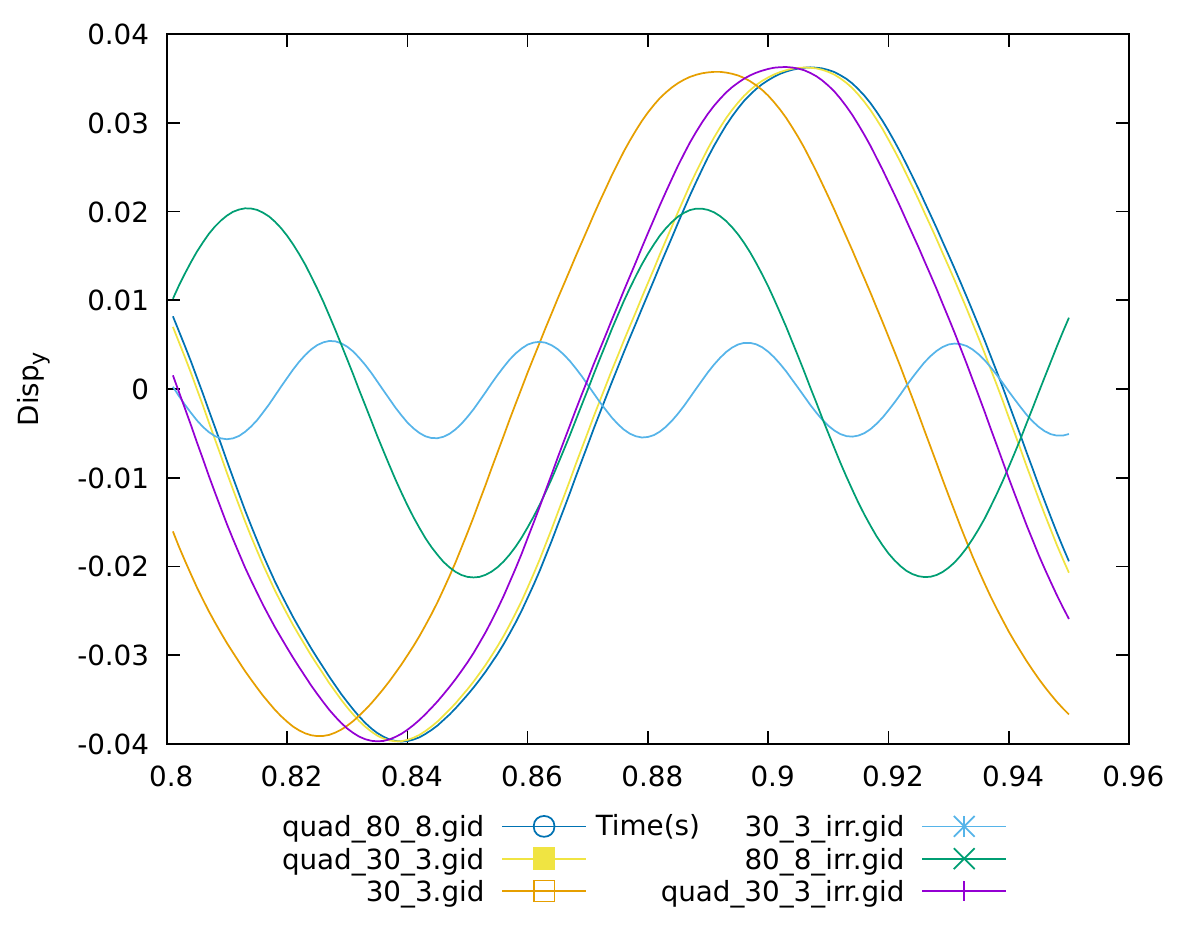}
        \caption{Displacement in $y$}
        \label{fig:2_disp_y} 
    \end{subfigure}
    \caption{Zoom of displacement at the tip of the beam}
    \label{fig:2_disp} 
\end{figure}

From Fig.~\ref{fig:2_disp_full} it can be seen that there are marked differences between the two formulations for a dynamic case. This is better observed in Fig.~\ref{fig:2_disp}, which shows a zoom of a portion of the time interval. The three-field formulation approximates a much finer reference solution obtained with linear elements, while the irreducible one is over-diffusive both in time and space. When quadratic elements are used, the irreducible formulation performs more accurately but the three-field is always more precise, in conserving both phase and amplitude.

\subsection{Three-field fluid-structure interaction tests}

In this section we show the benchmarking of the FSI problem by means of two well known test cases. The first one is a dynamic problem that converges to a stationary solution, and the second one is a fully transient case.

\subsubsection{Semi-stationary bending of a beam}\label{test3}

This semi-stationary problem, taken from \cite{Baiges2009}, consists of a supported
beam perpendicular to a fluid stream. Once the flow starts from the left wall it will
bend the beam. For the particular conditions of the test, a force balance between the
tractions imposed by the fluid and the stress on the beam will be achieved, in which the
beam will then remain bent without significant oscillation. The test conditions are shown in Table~\ref{3_physparam}.

\begin{table}[h!] \caption{Physical parameters (SI units)}\label{3_physparam}
\centering
    \begin{tabular}{cccc} 
        \hline      & Fluid     &            & Solid      \\ 
  \hline $\rho_{\textrm{fl}}$ & 2.0       & $\rho_{\textrm{sl}}$ &  10.0      \\ 
         $\nu_{\textrm{fl}} $ & 0.2       & $\nu_{\textrm{sl}} $ & 0.142857   \\
                    &           & Young      & 55\,428.0   \\
          Model     & Newtonian &            & Neo-Hookean \\
         \hline 
     \end{tabular} 
 \end{table}

\begin{figure}[h!]
    \centering
    \includegraphics[width=0.95\textwidth]{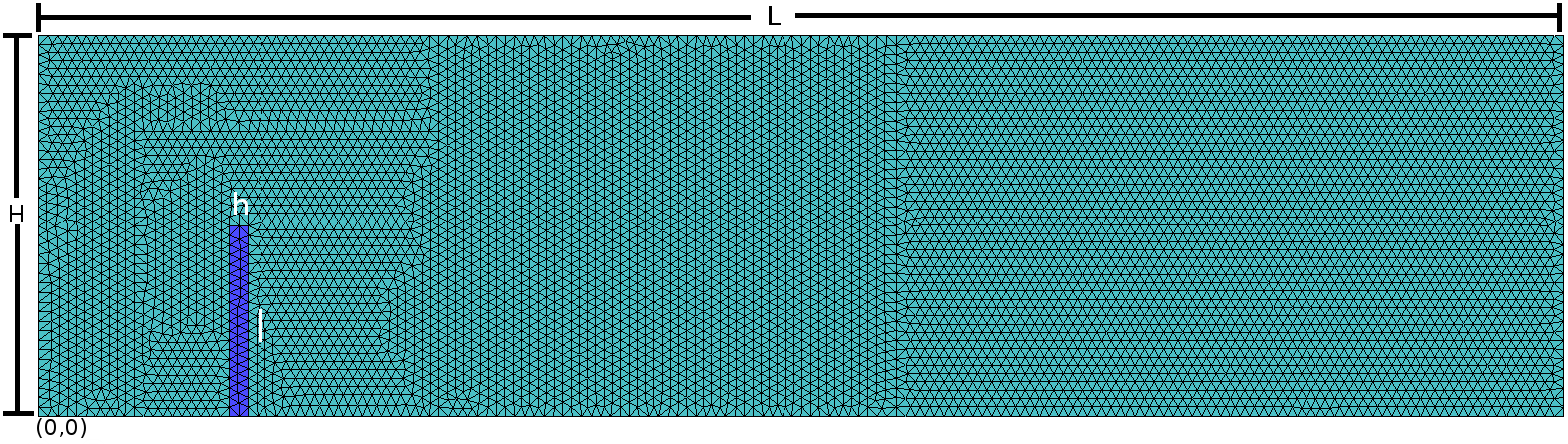}
    \caption{Geometry and mesh used for semi-stationary bending of a beam}
    \label{fig:3_stat_fsi_geomesh} 
\end{figure}

Fig.~\ref{fig:3_stat_fsi_geomesh} shows the geometry and mesh for the test, where $H=20$, $L=80$, $h=1$, $l=10$. Table~\ref{3_meshparams} shows important mesh parameters and Table~\ref{bcex3} the boundary conditions. To improve readability, results are assigned a suffix that corresponds to Table~\ref{3_caseparams}; all element used are quadratic triangles. ``Three-field" formulation refers to fully coupled three-field interaction, which means that both fluid and solid use the three-field approach. On the other hand ``standard" refers to the usual coupling methodology, in which a displacement formulation is used for the solid and a velocity-pressure formulation for the fluid. The last column in Table~\ref{3_caseparams} shows the number of DOFs used for each case. Note that even though cases A and B have different amount of DOFs, these cases have been run with the same mesh. Case C uses a much coarser mesh.

\begin{table}[h!] \caption{Mesh parameters}\label{3_meshparams}
\centering
    \begin{tabular}{cccc} 
        \hline          & Fluid             & Solid           \\ 
  \hline  Element type  & Quadratic triangle   & Quadratic triangle \\ 
          Nodes per element &  6                & 6               \\
          \# of elements&  14\,308           & 78              \\
          \# of nodes   &  29\,057           & 201             \\
         \hline 
     \end{tabular} 
 \end{table}

\begin{table}[h!] \caption{Boundary conditions}\label{bcex3}
  \centering
    \begin{tabular}{cccc} 
        \hline               Fluid                     & Solid               \\ 
  \hline              
  $x = 0$: $u_x$ = 1, $u_y$ = 0  & $y = 0$: $d_x$ = $d_y$ = 0 \\ 
  $y = 0, H$: Free slip                 &         Other boundaries: fluid tractions            \\ 
  $x = L$:  Free                      &                     \\ 
  Other boundaries: solid velocities & \\
          \hline 
     \end{tabular} 
 \end{table}
 
 \begin{table}[h!] \caption{Case parameters}
\centering
    \begin{tabular}{cccc} 
        \hline      Name          &  FSI formulation & time step & Degrees of freedom \\ 
  \hline            dt\_0075\_A    &  three-field     &  0.075    &        175\,548  \\ 
                    dt\_01\_A      &  three-field     &  0.01     &        175\,548   \\
                    dt\_02\_A      &  three-field     &  0.02     &        175\,548   \\
                    dt\_03\_A      &  three-field     &  0.03     &        175\,548    \\ 
  \hline            dt\_0075\_B    &  Standard        &  0.075    &        87\,573  \\ 
                    dt\_01\_B      &  Standard        &  0.01     &        87\,573   \\
                    dt\_02\_B      &  Standard        &  0.02     &        87\,573   \\
                    dt\_03\_B      &  Standard        &  0.03     &        87\,573    \\ 
  \hline            dt\_0075\_C    &  three-field     &  0.075    &        58\,320  \\ 
                    dt\_01\_C      &  three-field     &  0.01     &        58\,320   \\
                    dt\_02\_C      &  three-field     &  0.02     &        58\,320   \\
                    dt\_03\_C      &  three-field     &  0.03     &        58\,320    \\ 
         \hline 
     \end{tabular} 
     \label{3_caseparams}
 \end{table}
  
In the next figures we compare cases A and B for velocity and pressure around the top of the beam for the fluid, and displacement and acceleration at the top of the beam for the solid.
 
\begin{figure}[h!]
    \centering
    \begin{subfigure}[h!]{0.45\textwidth}
        \includegraphics[width=\linewidth]{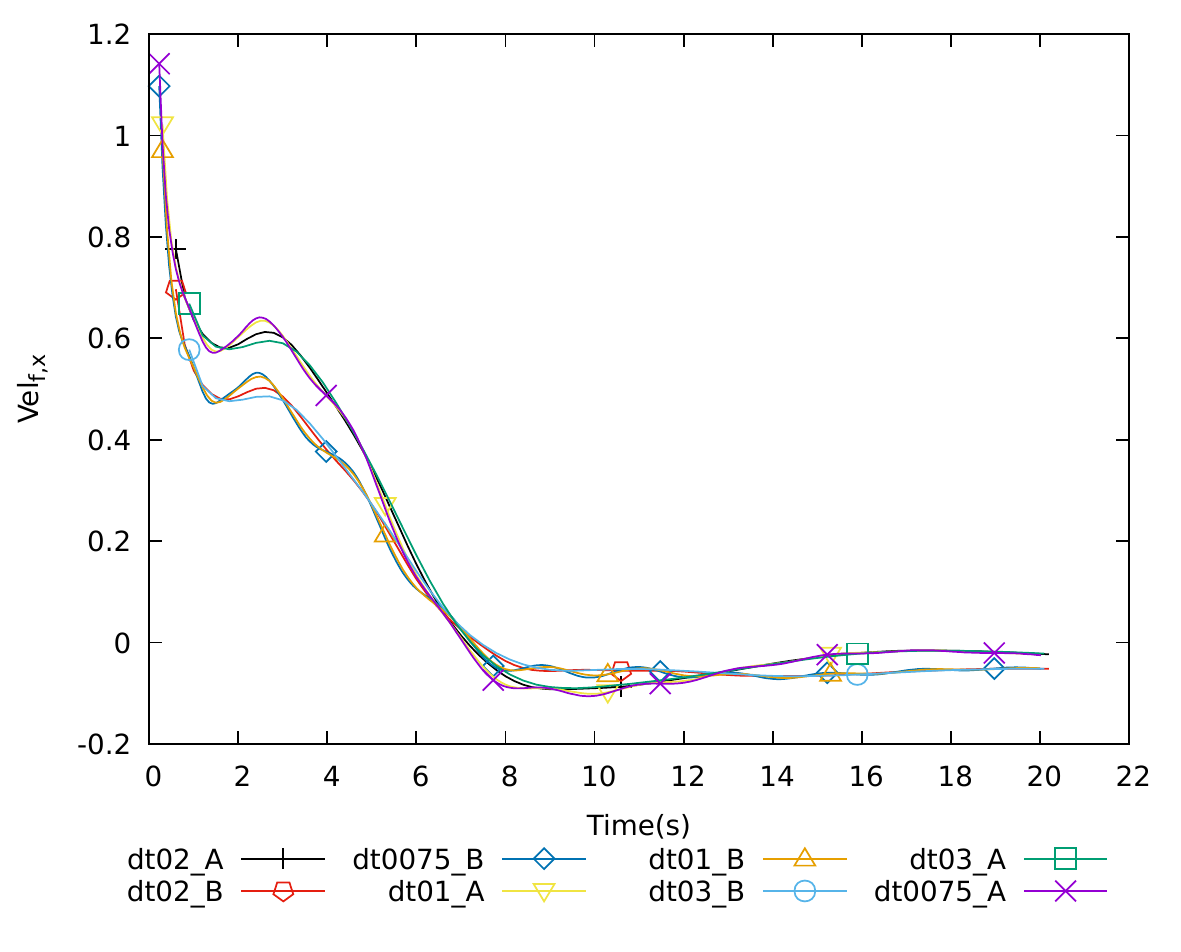}
        \caption{Velocity in $x$}
        \label{fig:3_vel_x} 
    \end{subfigure}
    \begin{subfigure}[h!]{0.45\textwidth}
        \includegraphics[width = \linewidth]{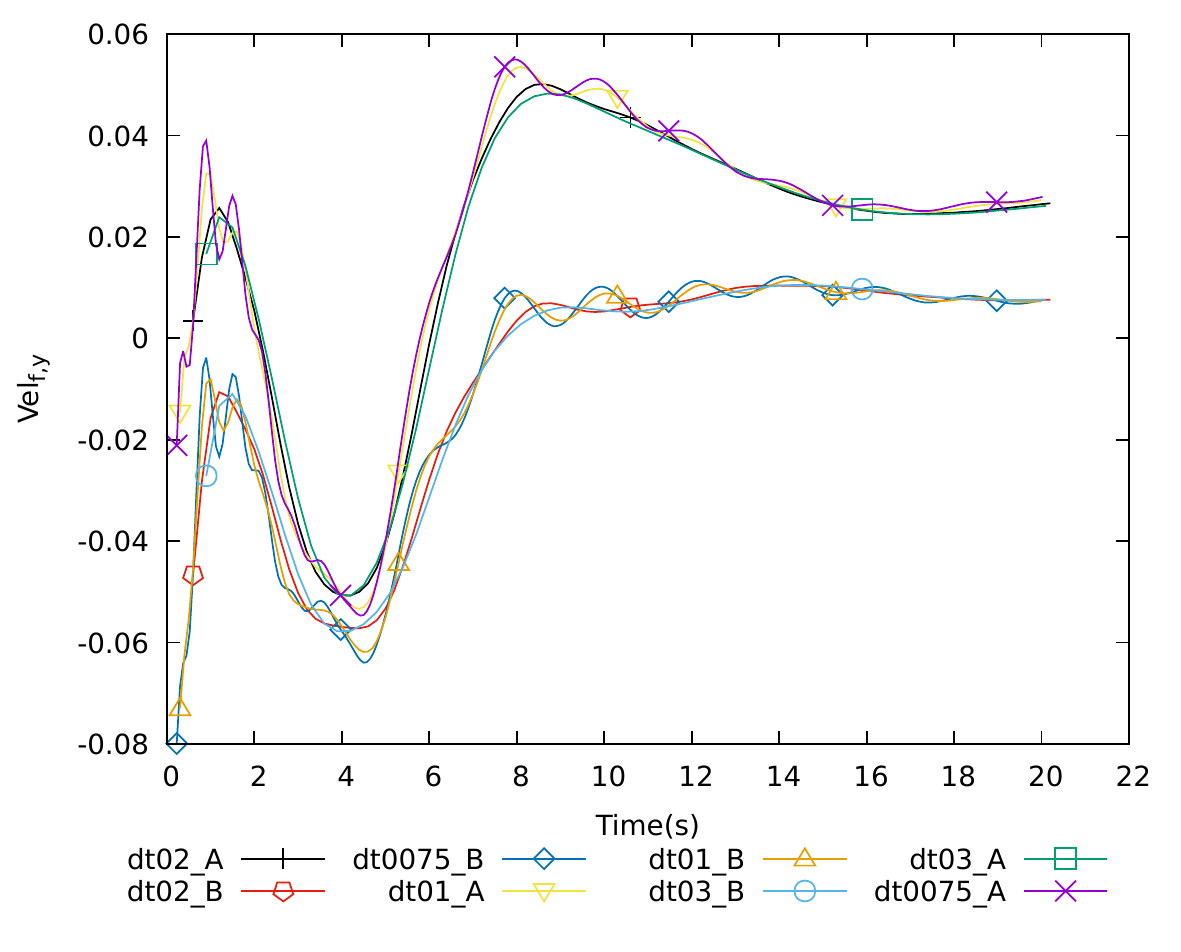}
        \caption{Velocity in $y$}
        \label{fig:3_vel_y} 
    \end{subfigure}
    \caption{Velocity around the tip of the beam}
    \label{fig:3_vel} 
\end{figure}

\begin{figure}[h!]
    \centering
    \includegraphics[width=0.45\textwidth]{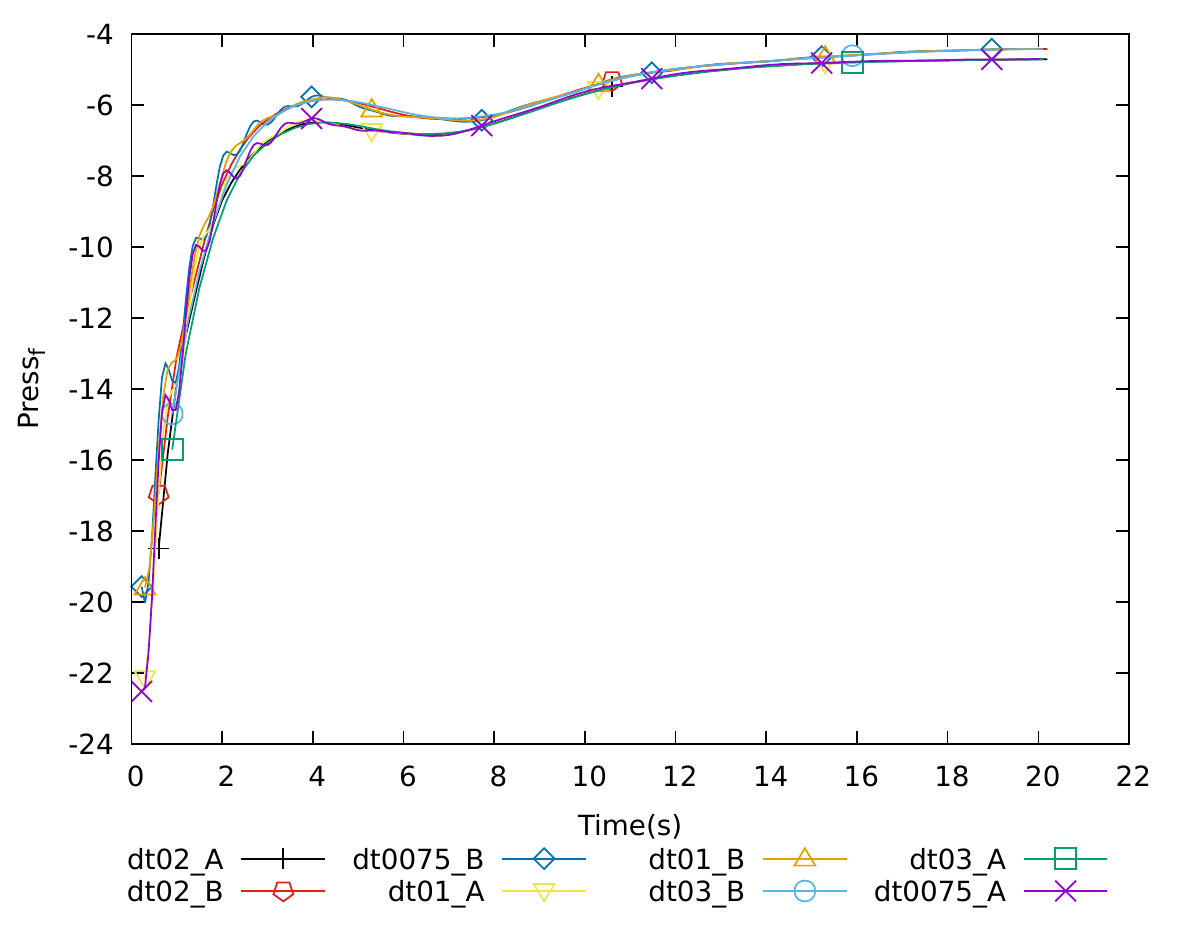}
    \caption{Pressure around the tip of the beam}
    \label{fig:3_press} 
\end{figure}

Fig.~\ref{fig:3_vel} shows the velocity around the top of the beam for the fluid. It can be seen that both formulation are able to capture higher frequency modes at lower time steps and show very similar behavior, specially for the pressure, shown in Fig.~\ref{fig:3_press}. Even though the initial transient has some differences, the solution tends to converge to a stationary solution. Overall similar responses can be seen from the fluid for both formulations, except for the initial transient. 
 
\begin{figure}[h!]
    \centering
    \begin{subfigure}[h!]{0.45\textwidth}
        \includegraphics[width=\linewidth]{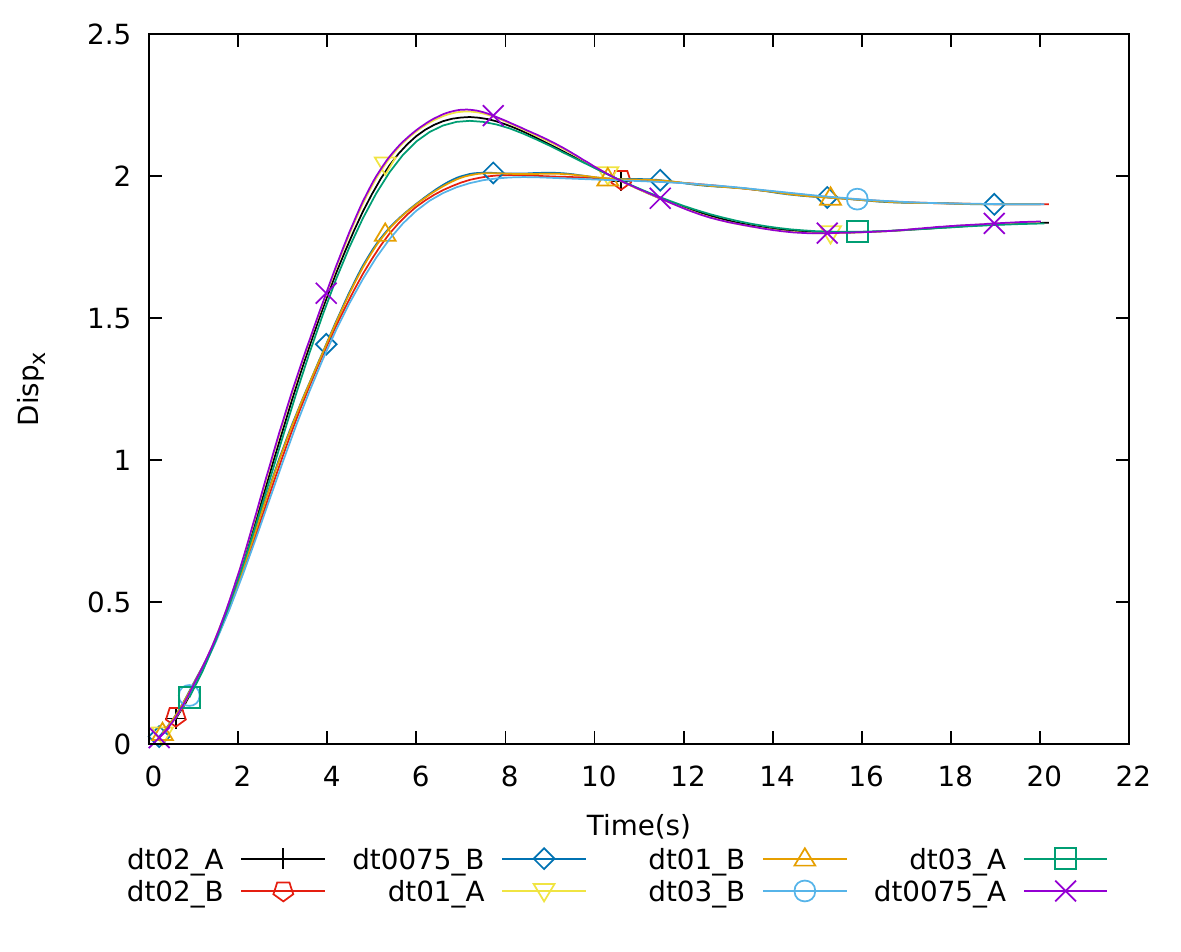}
        \caption{Displacement in $x$}
        \label{fig:3_disp_x} 
    \end{subfigure}
    \begin{subfigure}[h!]{0.45\textwidth}
        \includegraphics[width = \linewidth]{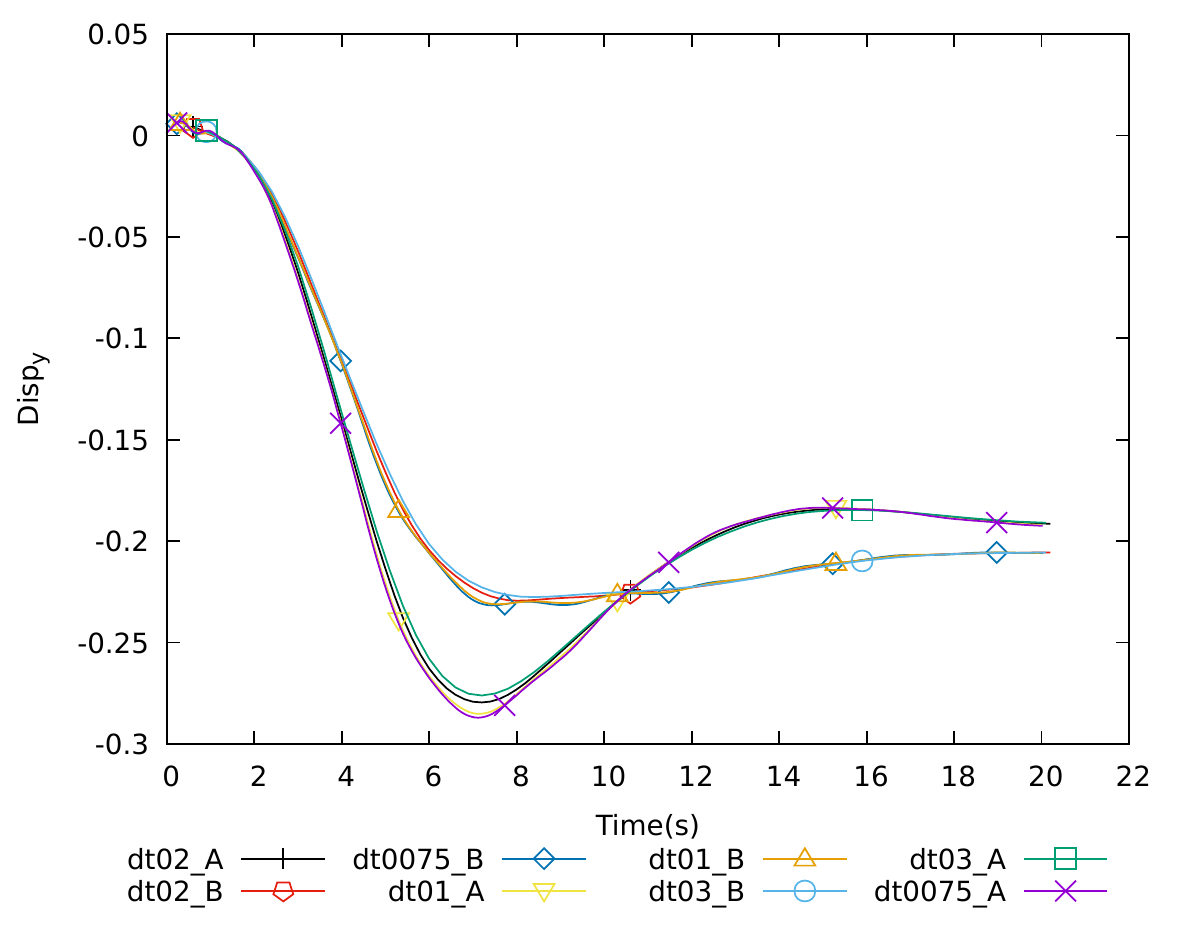}
        \caption{Displacement in $y$}
        \label{fig:3_disp_y} 
    \end{subfigure}
    \caption{Displacement at the tip of the beam}
    \label{fig:3_disp} 
\end{figure}

For the solid the standard formulation tends to be over-dissipative in the initial transient with regard to the displacements, as seen in Fig.~\ref{fig:3_disp}, while the three-field is much less dissipative thus capturing more frequencies of the response. While the stationary solution appears to be achieved faster for the standard method, the three-field approach would need longer time sampling to achieve a possible steady state. Both formulations appear to converge to the same stationary value.
 
In order to make a correct comparison between methods, a similar number of DOFs should be used. In the following results we explore the effect of a reduction of DOFs in the solution by comparing cases A and C, this is, for the same three-field formulation we explore the effect of reducing more than three times the amount of DOFs. 

\begin{figure}[h!]
    \centering
    \begin{subfigure}[h!]{0.45\textwidth}
        \includegraphics[width=\linewidth]{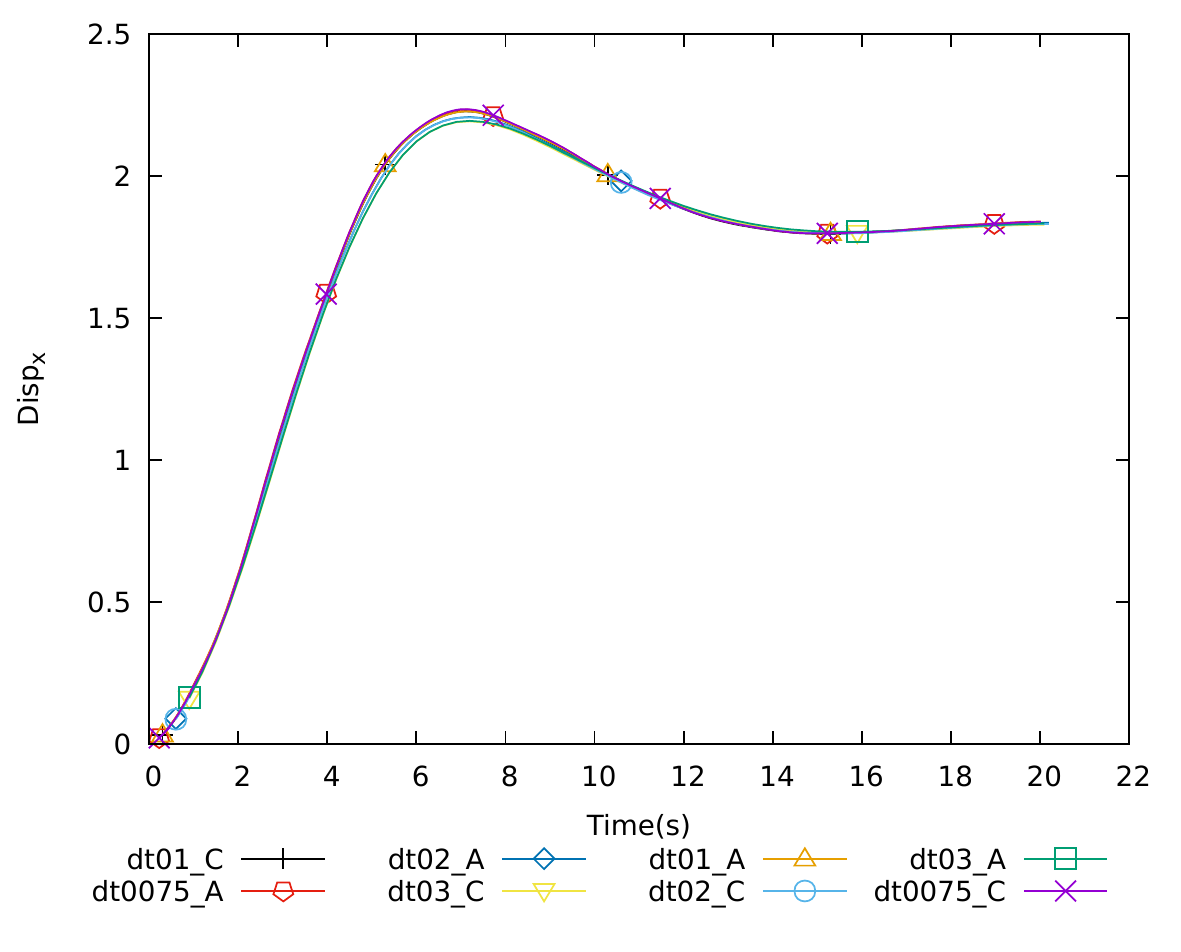}
        \caption{Displacement in $x$}
        \label{fig:3_disp_x_dof} 
    \end{subfigure}
    \begin{subfigure}[h!]{0.45\textwidth}
        \includegraphics[width = \linewidth]{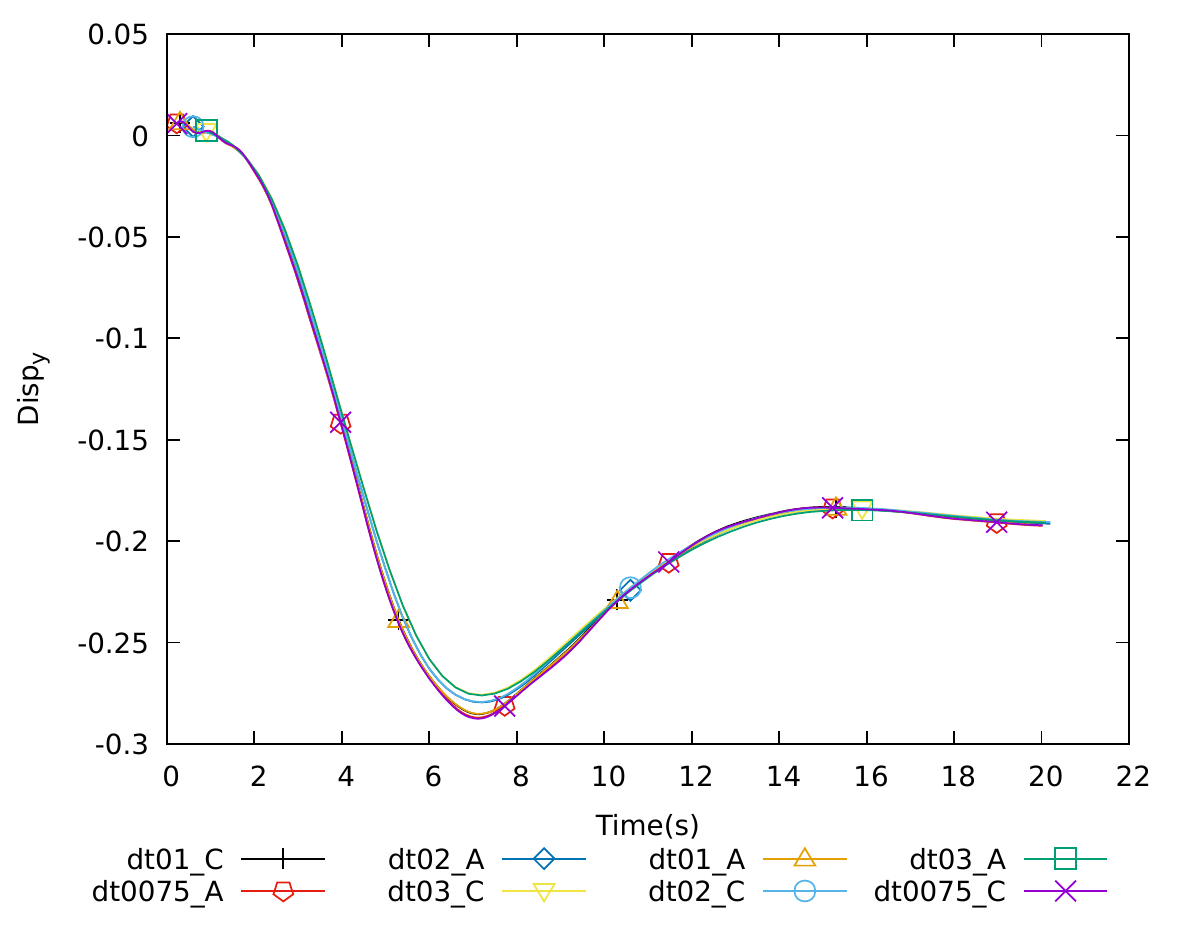}
        \caption{Displacement in $y$}
        \label{fig:3_disp_y_dof} 
    \end{subfigure}
    \caption{Displacement for the tip of the beam}
    \label{fig:3-disp-dof} 
\end{figure}

\begin{figure}[h!]
    \centering
    \begin{subfigure}[h!]{0.45\textwidth}
        \includegraphics[width=\linewidth]{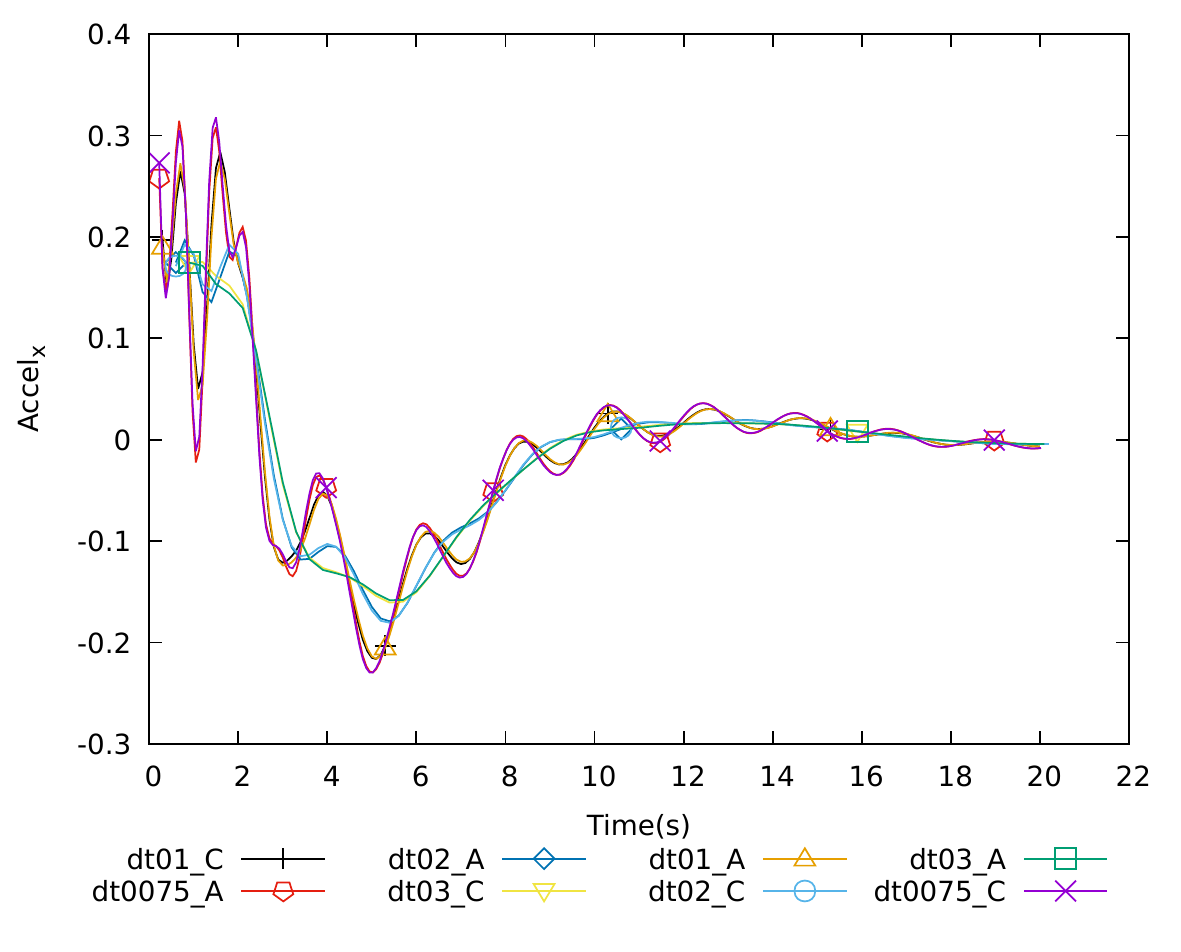}
        \caption{Acceleration in $x$}
        \label{fig:3_accel_x_dof} 
    \end{subfigure}
    \begin{subfigure}[h!]{0.45\textwidth}
        \includegraphics[width = \linewidth]{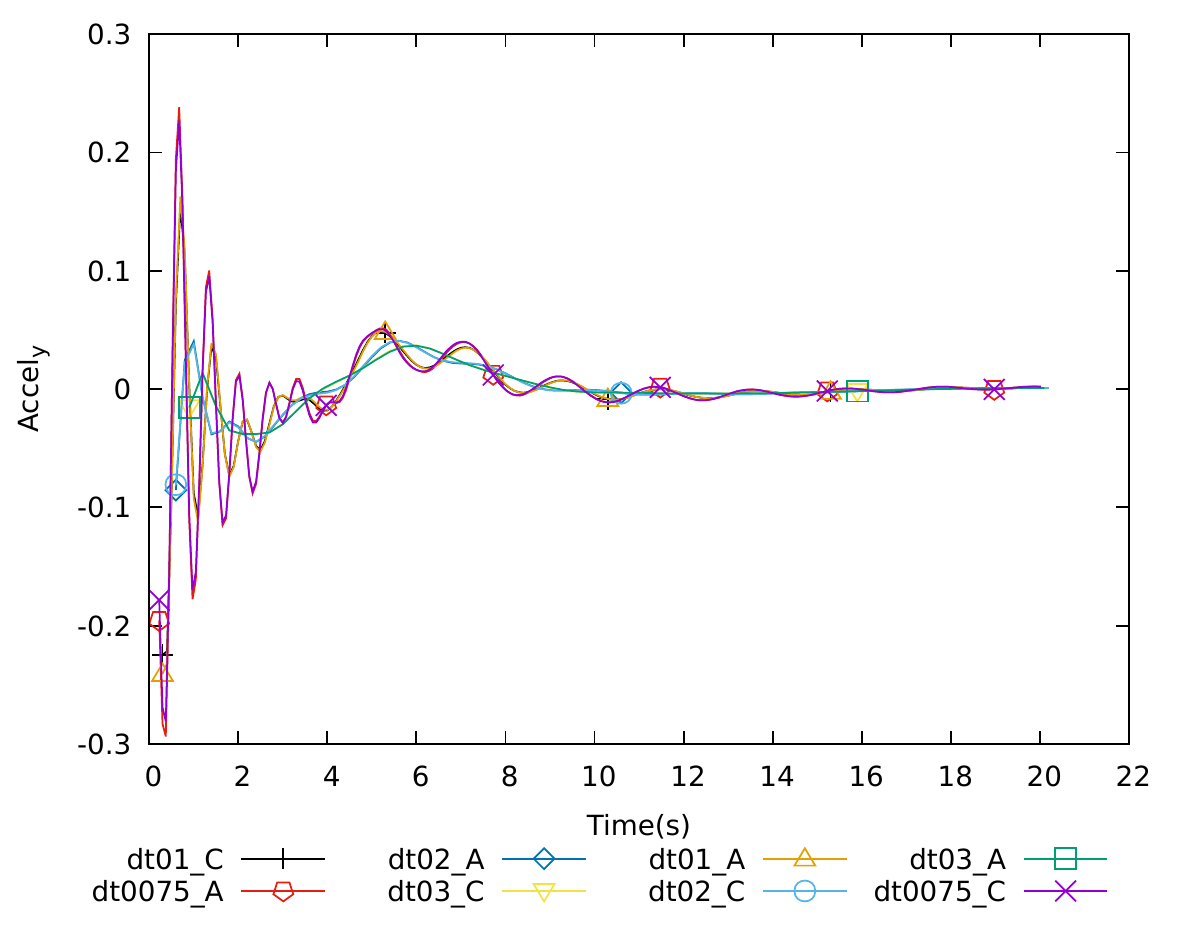}
        \caption{Acceleration in $y$}
        \label{fig:3_accel_y_dof} 
    \end{subfigure}
    \caption{Acceleration at the tip of the beam}
    \label{fig:3-accel-dof} 
\end{figure}

Figs.~\ref{fig:3-disp-dof} and \ref{fig:3-accel-dof} show the displacement and acceleration evolution, respectively, over the selected time frame. It can be seen that the solution does not appear to change importantly for a significant change in time step or mesh size. For a coarser mesh and a higher time step there is a slight loss of amplitude. Regardless of the mesh size, a finer time step captures successfully higher frequency modes of the solution.

\begin{figure}[h!]
    \centering
    \begin{subfigure}[h!]{0.45\textwidth}
        \includegraphics[width=\linewidth]{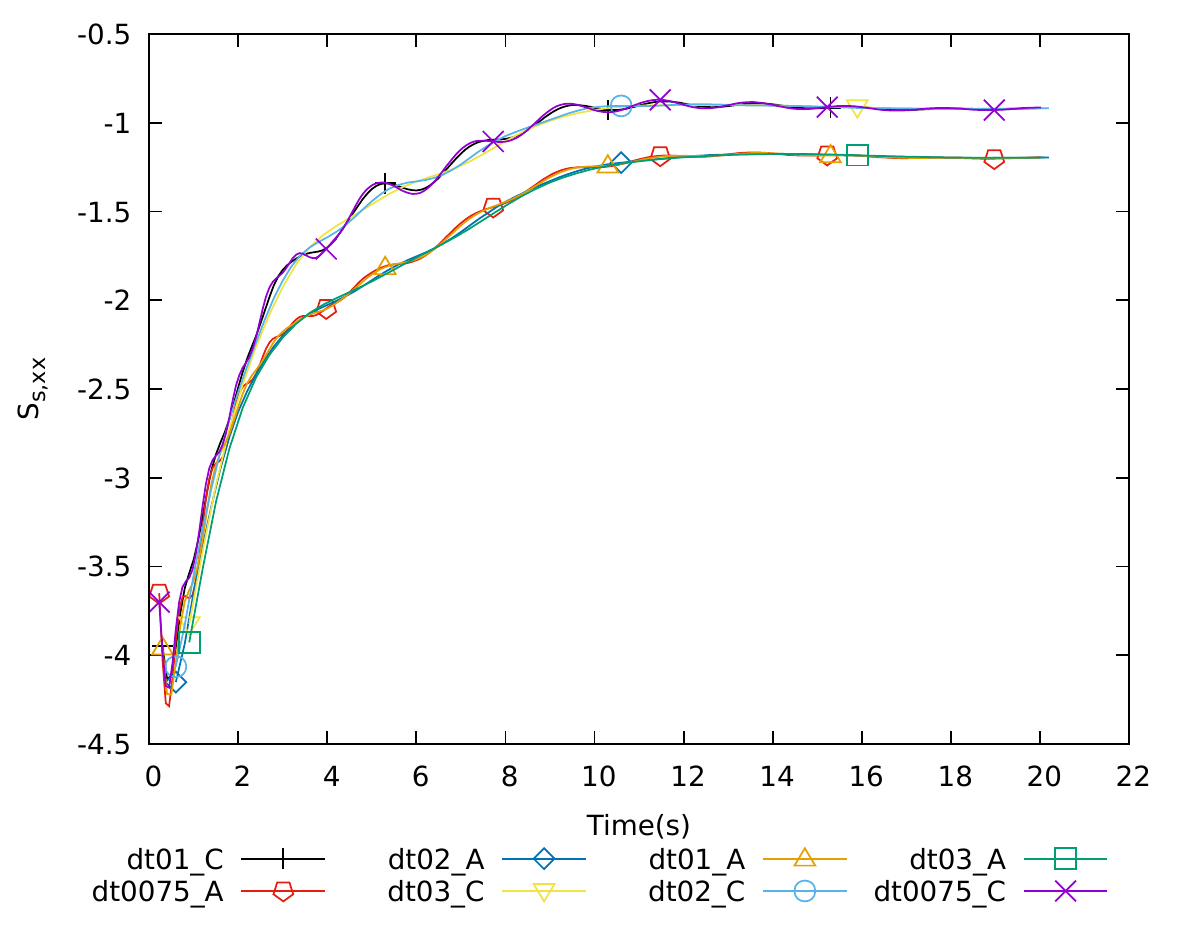}
        \caption{Stress $S_{xx}$}
        \label{fig:3_stress_xx_dof} 
    \end{subfigure}
    \begin{subfigure}[h!]{0.45\textwidth}
        \includegraphics[width = \linewidth]{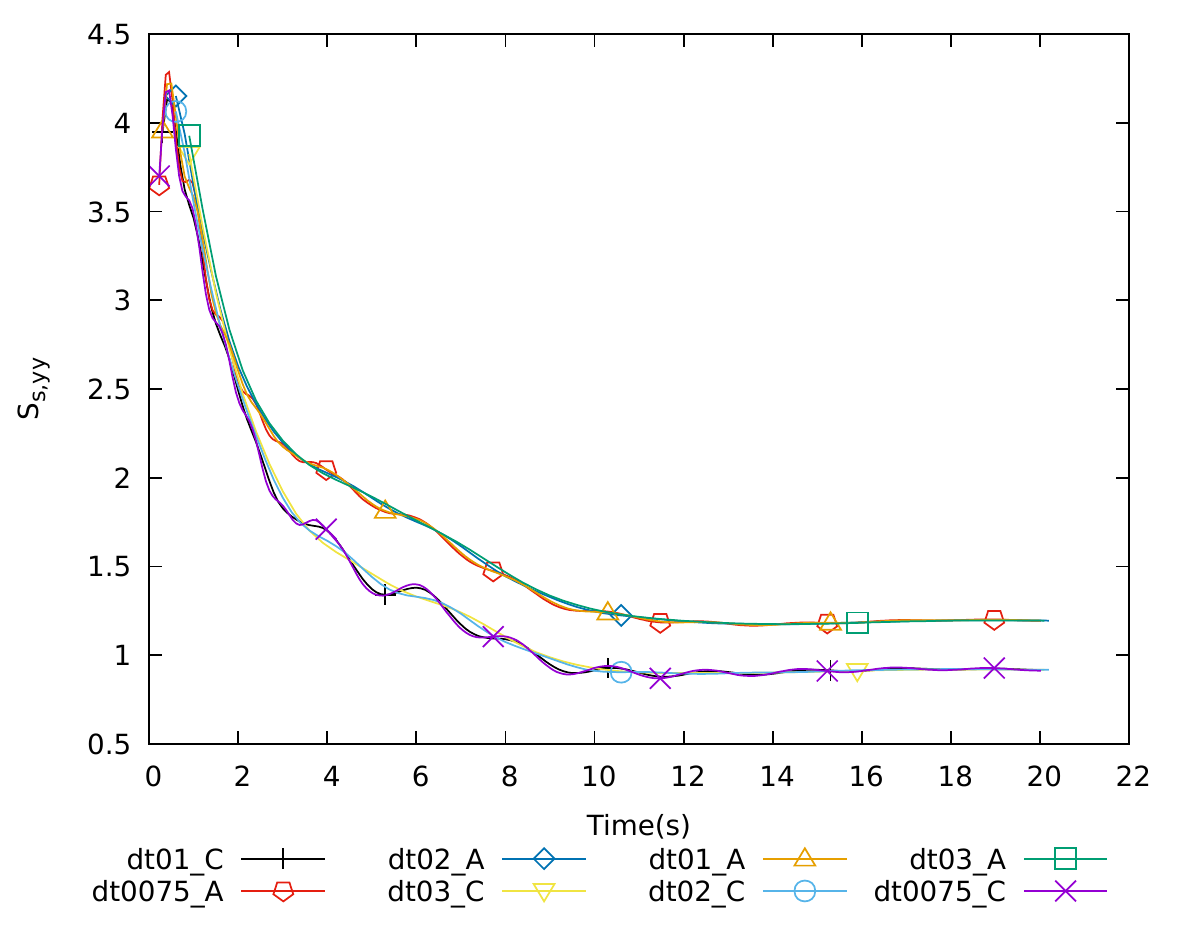}
        \caption{Stress $S_{yy}$}
        \label{fig:3_stress_yy_dof} 
    \end{subfigure}
    \caption{Deviatoric stress at the tip of the beam}
    \label{fig:3_stress_dof} 
\end{figure}

\begin{figure}[h!]
    \centering
    \begin{subfigure}[h!]{0.45\textwidth}
        \includegraphics[width=\linewidth]{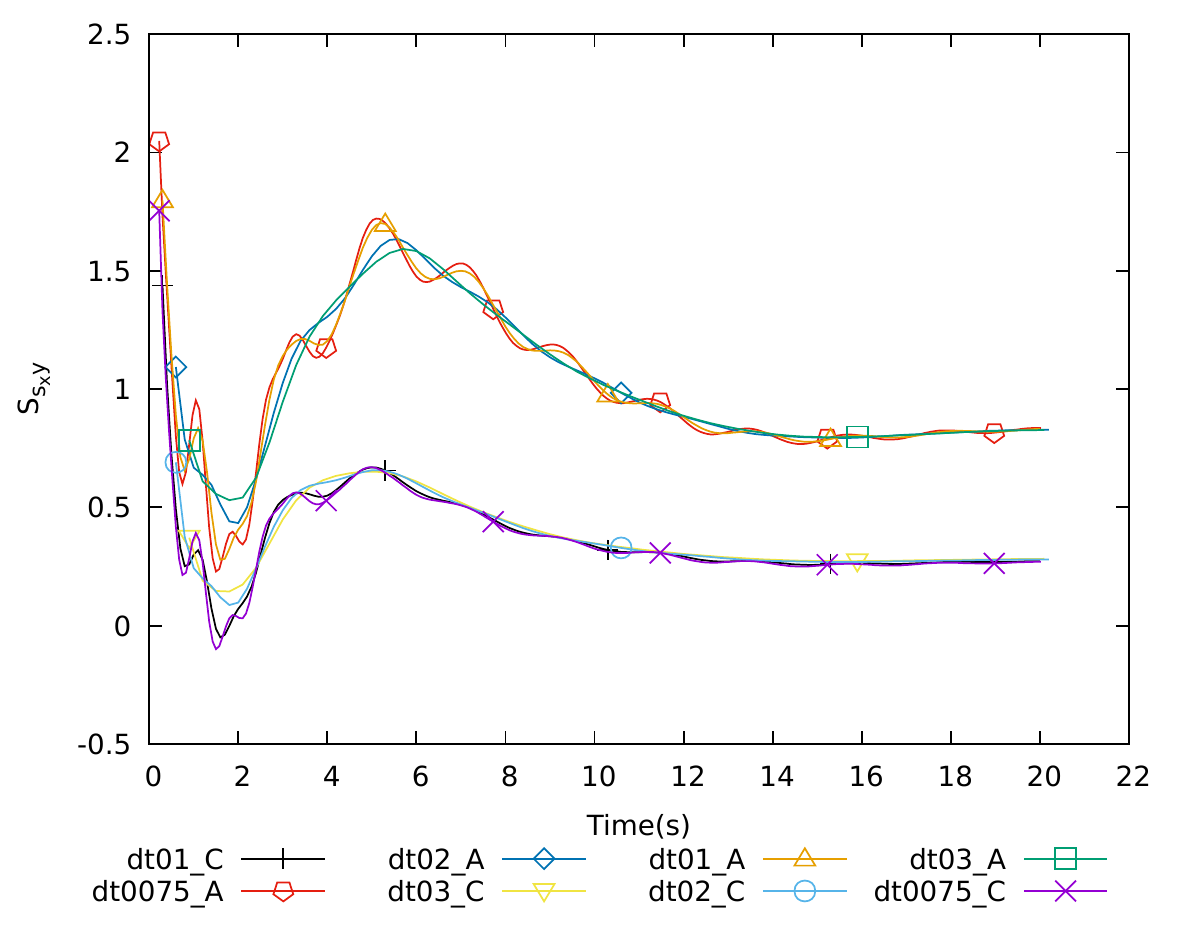}
        \caption{Stress $S_{xy}$}
        \label{fig:3_stress_xy_dof} 
    \end{subfigure}
    \begin{subfigure}[h!]{0.45\textwidth}
        \includegraphics[width = \linewidth]{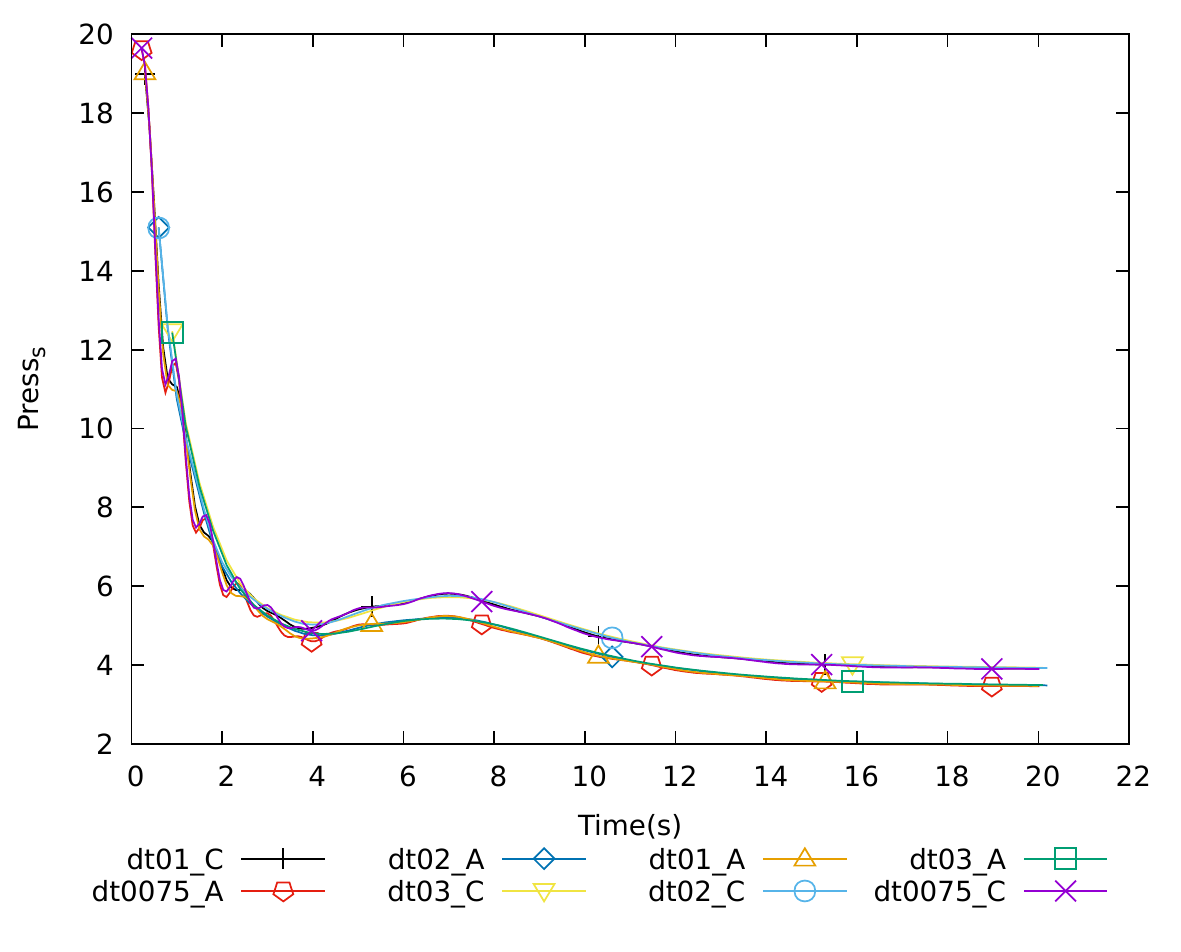}
        \caption{Pressure}
        \label{fig:3_pressure_dof} 
    \end{subfigure}
    \caption{Deviatoric stress and pressure at the tip of the beam}
    \label{fig:3_stress_press_dof} 
\end{figure}

Figs.~\ref{fig:3_stress_dof} and \ref{fig:3_stress_press_dof} show the evolution of deviatoric stress and pressure at the top of the beam. It can be seen that the solution is continuous and stable. Overall a coarser mesh produces a slightly higher magnitude for the stresses $S_{xx}$ and $S_{yy}$. The case is opposite for the shear stress $S_{xy}$ where the finer mesh produces more pronounced maxima and minima during the first ten seconds of simulation. Approximation of the solution seems to be accurate enough for both tests.

\subsubsection{3D Flow around a plate}\label{test5}

The following example is a 3D version of the one shown in Section~\ref{test3}. The test conditions are shown in Table~\ref{5_physparam}.

\begin{table}[h!] \caption{Physical parameters}\label{5_physparam}
\centering
    \begin{tabular}{cccc} 
        \hline            & Fluid         &                     & Solid      \\ 
  \hline $\rho_{\textrm{fl}}$ & 100.0     & $\rho_{\textrm{sl}}$ & 1\,000.0   \\ 
         $\nu_{\textrm{fl}} $ & 1.0       & $\nu_{\textrm{sl}} $ & 0.48        \\
                          &               & Young               & $300\times 10^3$      \\
          model           & Newtonian     &                     & NeoHookean \\
         \hline 
     \end{tabular} 
 \end{table}

\begin{figure}[h!]
    \centering
    \includegraphics[width=0.85\textwidth]{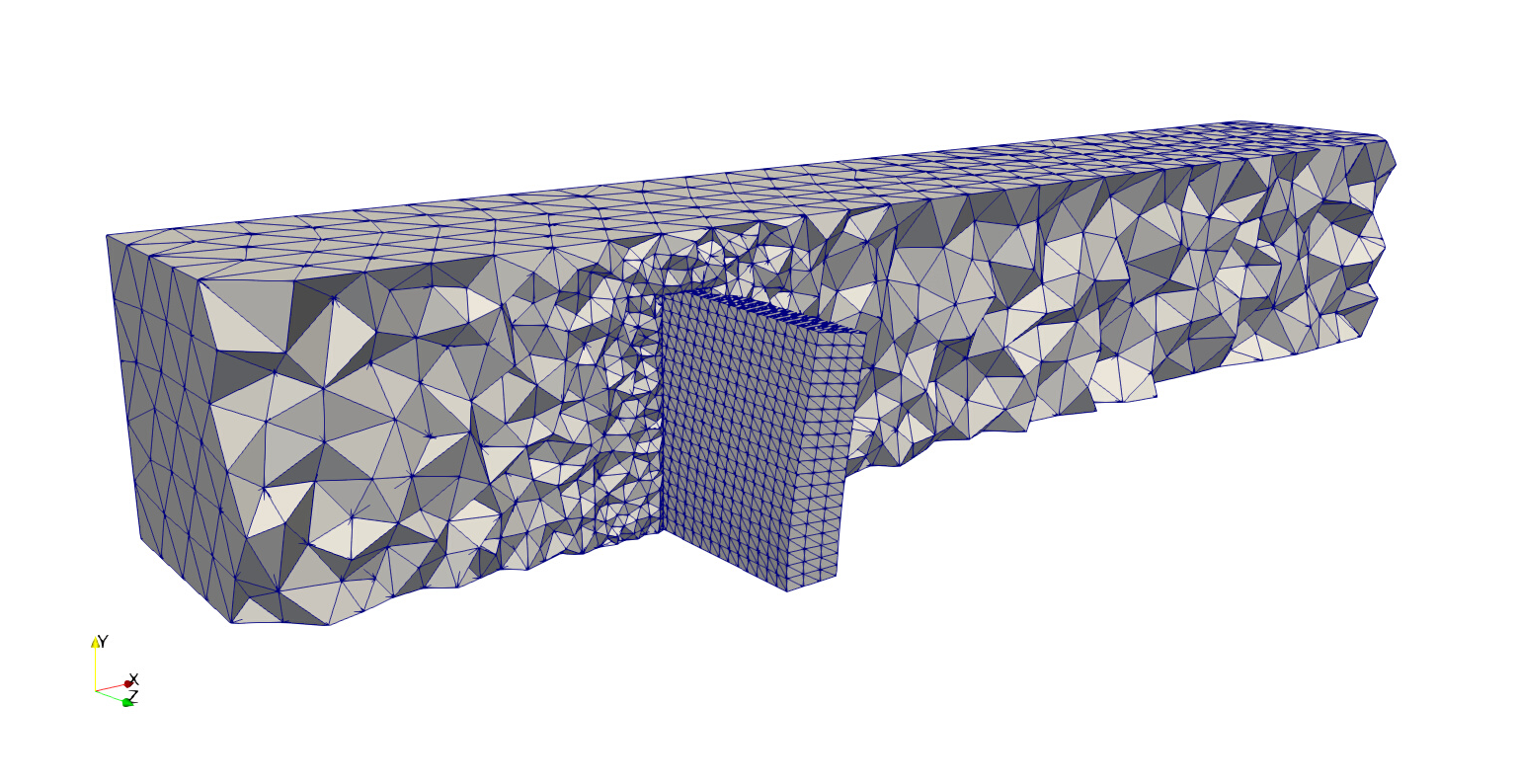}
    \caption{Case geometry}
    \label{fig:5_test_geo} 
\end{figure}

The geometry is shown in Fig.~\ref{fig:5_test_geo}. Table~\ref{5_meshparams} shows important mesh parameters, geometrical parameters are shown in Table~\ref{5_params} and Table~\ref{5_bcs} shows the boundary conditions.

\begin{table}[h!] \caption{Mesh parameters}\label{5_meshparams}
\centering
    \begin{tabular}{cccc} 
        \hline          & Fluid                 & Solid                \\ 
  \hline  Element type  & Linear Tetrahedra   & Linear Tetrahedra  \\ 
          Nodes per element &  4                & 4                    \\
          \# of elements&  87\,941               & 10\,833               \\
          \# of nodes   &  16\,785               & 2\,491                \\
         \hline 
     \end{tabular} 
 \end{table}
 
\begin{table}[h!] \caption{Channel dimensions and flow parameters}\label{5_params}
  \centering
    \begin{tabular}{cc} 
        \hline                 &                    \\          
        Height (H):        & 0.5                    \\    
        Width :            & 1.0                    \\                    
        Length:            & 3.0                    \\                    
        Inlet mean velocity ($\overline{U}$) : & 1.0  \\
        \hline 
    \end{tabular} 
 \end{table}
 
\begin{table}[h!] \caption{Boundary conditions}\label{5_bcs}
  \centering
    \begin{tabular}{cccc} 
        \hline                   & Fluid                 & Solid               \\ 
        \hline              
        Flow inlet:              & $1.5\overline{U}\cdot \frac{y\cdot(0.5-y)\cdot z\cdot(1.0-z)}{\left( \frac{H}{2} \right)^{2}}$  \\ 
        Flow outlet:             & Free                  &                     \\ 
        Channel walls:           & No slip               &                     \\
        Plate sides:             & Solid velocities      & Fluid tractions     \\
        Plate bottom:            &                       & Fixed               \\
        \hline 
    \end{tabular} 
 \end{table}

\begin{figure}[h!]
    \centering
    \begin{subfigure}[h!]{0.3\textwidth}
        \includegraphics[width=\linewidth]{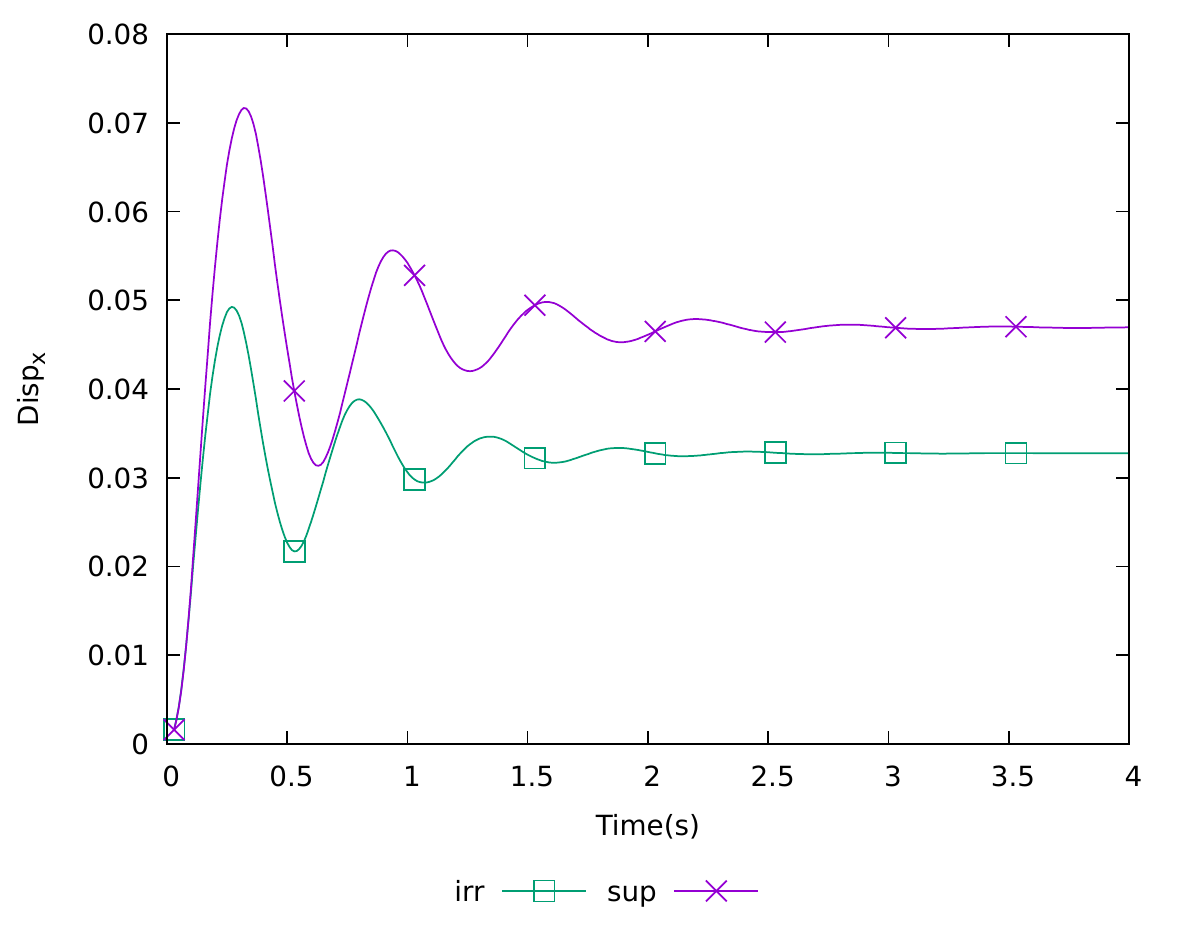}
        \caption{Displacement in $x$}
        \label{fig:5_disp_x} 
    \end{subfigure}
    \begin{subfigure}[h!]{0.3\textwidth}
        \includegraphics[width = \linewidth]{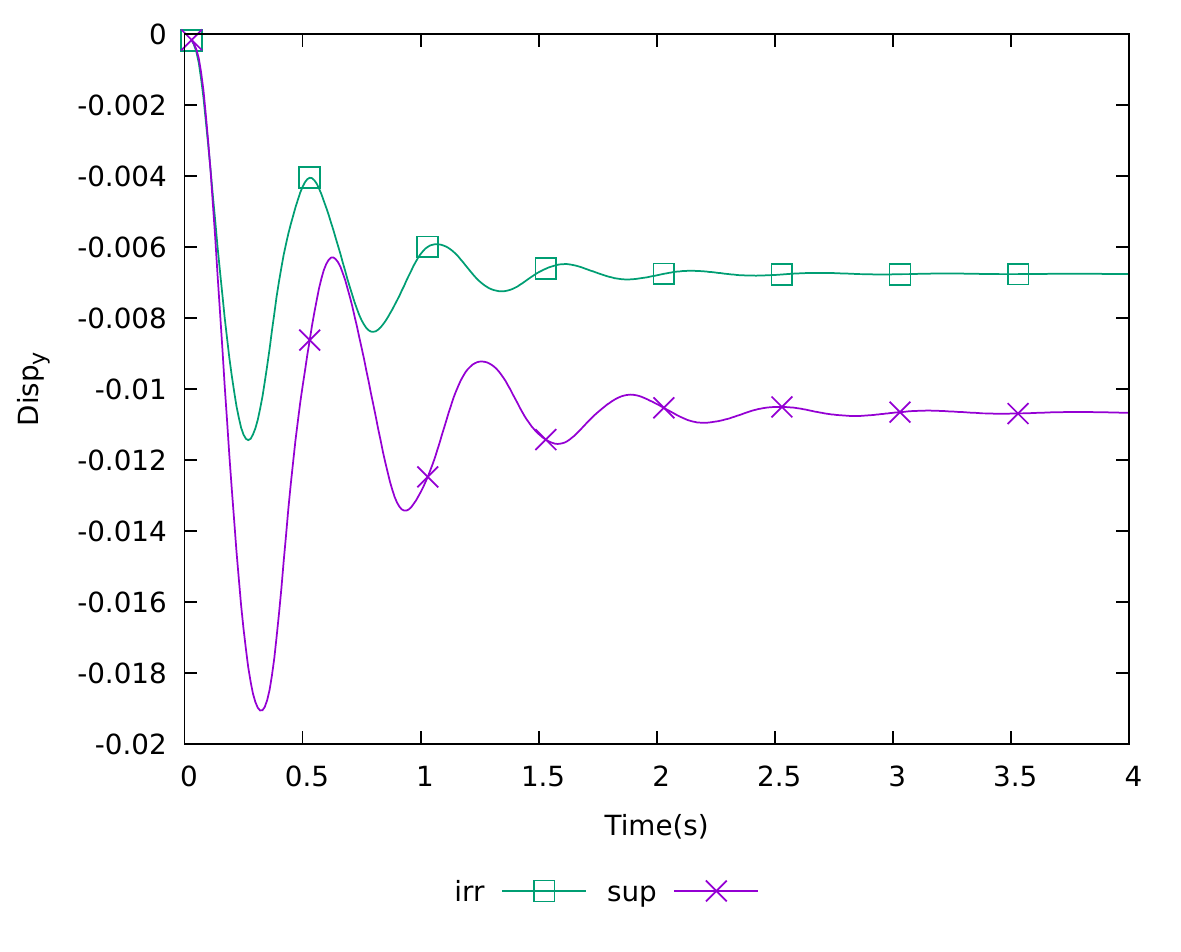}
        \caption{Displacement in $y$}
        \label{fig:5_disp_y} 
    \end{subfigure}
    \begin{subfigure}[h!]{0.3\textwidth}
        \includegraphics[width = \linewidth]{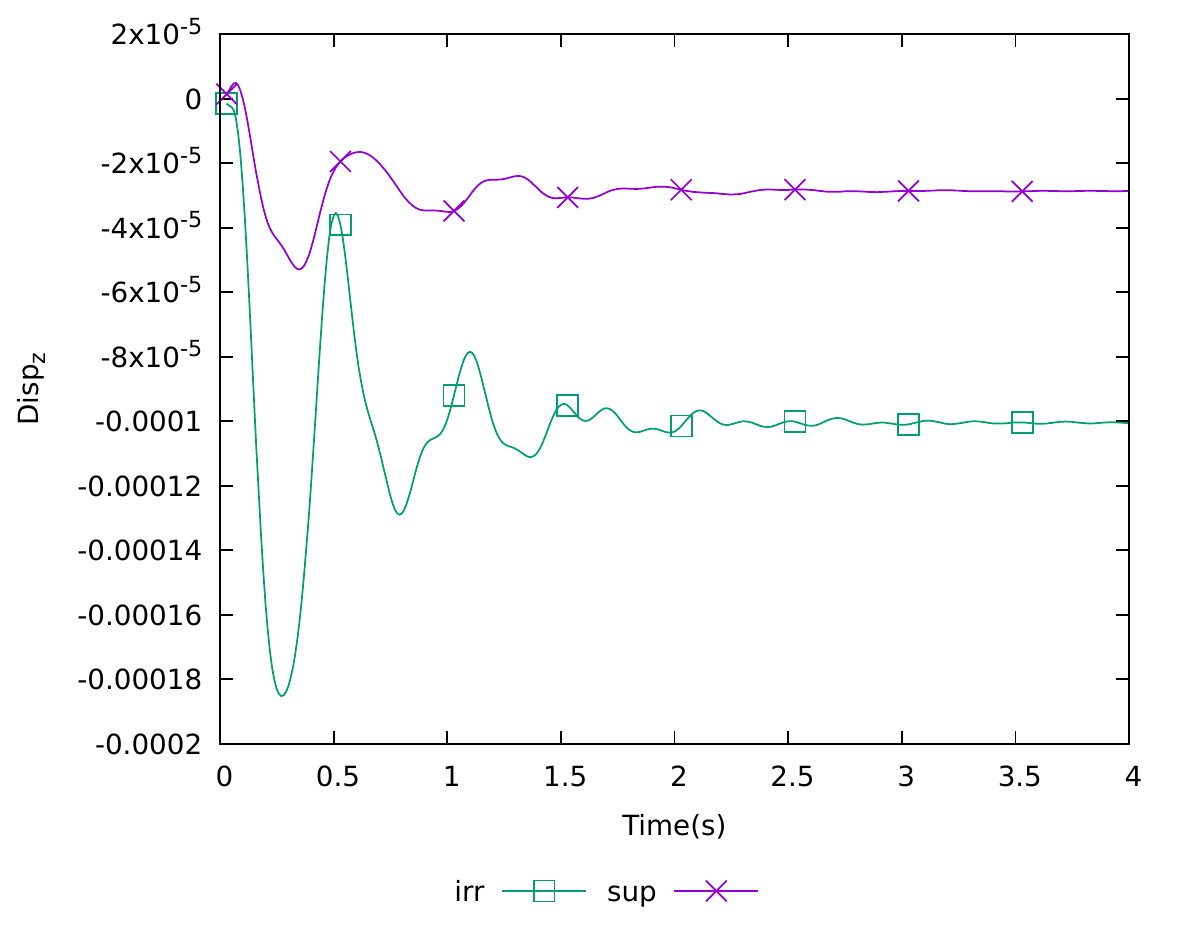}
        \caption{Displacement in $z$}
        \label{fig:5_disp_z} 
    \end{subfigure}
    \caption{Displacement at the tip of the plate}
    \label{fig:5_disp} 
\end{figure}

\begin{figure}[h!]
    \centering
    \begin{subfigure}[h!]{0.3\textwidth}
        \includegraphics[width=\linewidth]{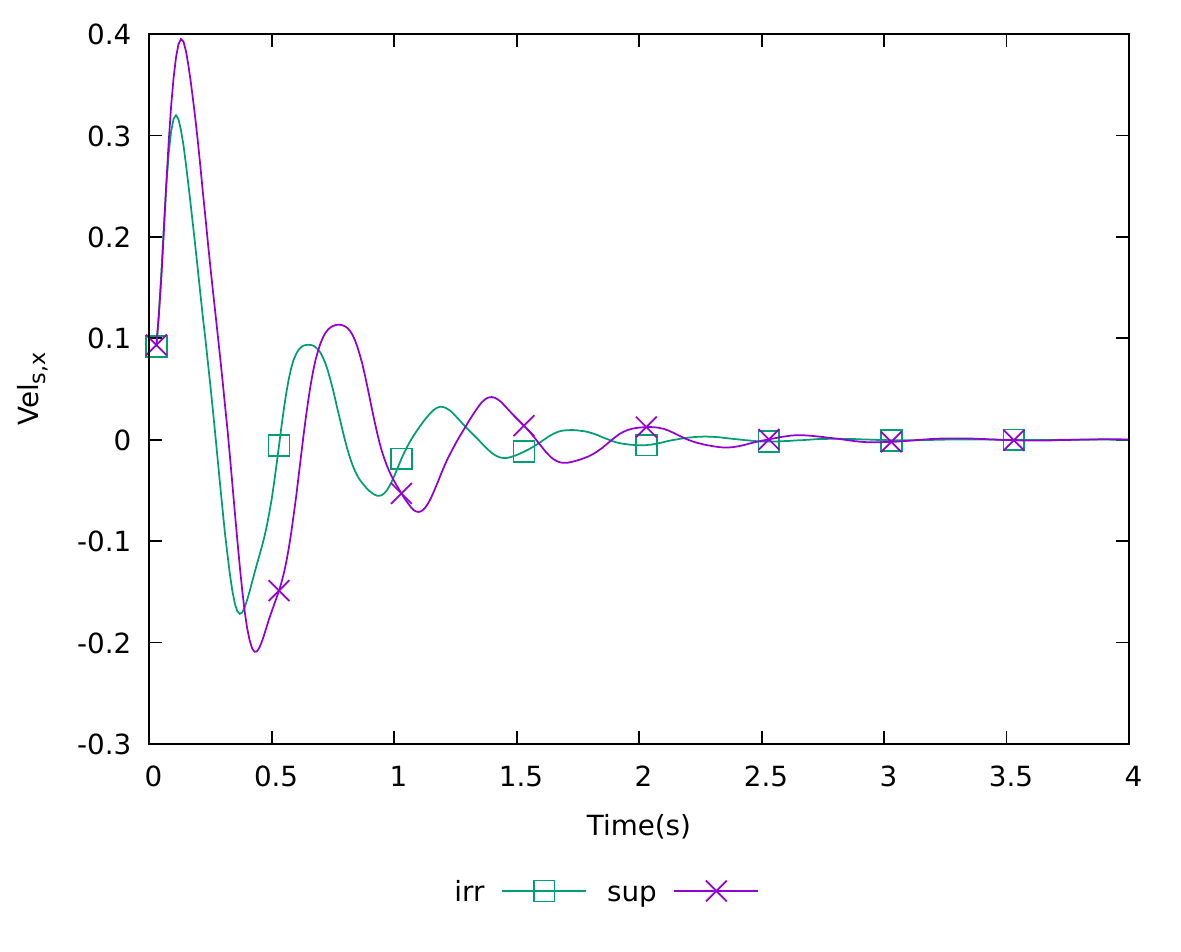}
        \caption{Velocity in $x$}
        \label{fig:5_vel_x} 
    \end{subfigure}
    \begin{subfigure}[h!]{0.3\textwidth}
        \includegraphics[width = \linewidth]{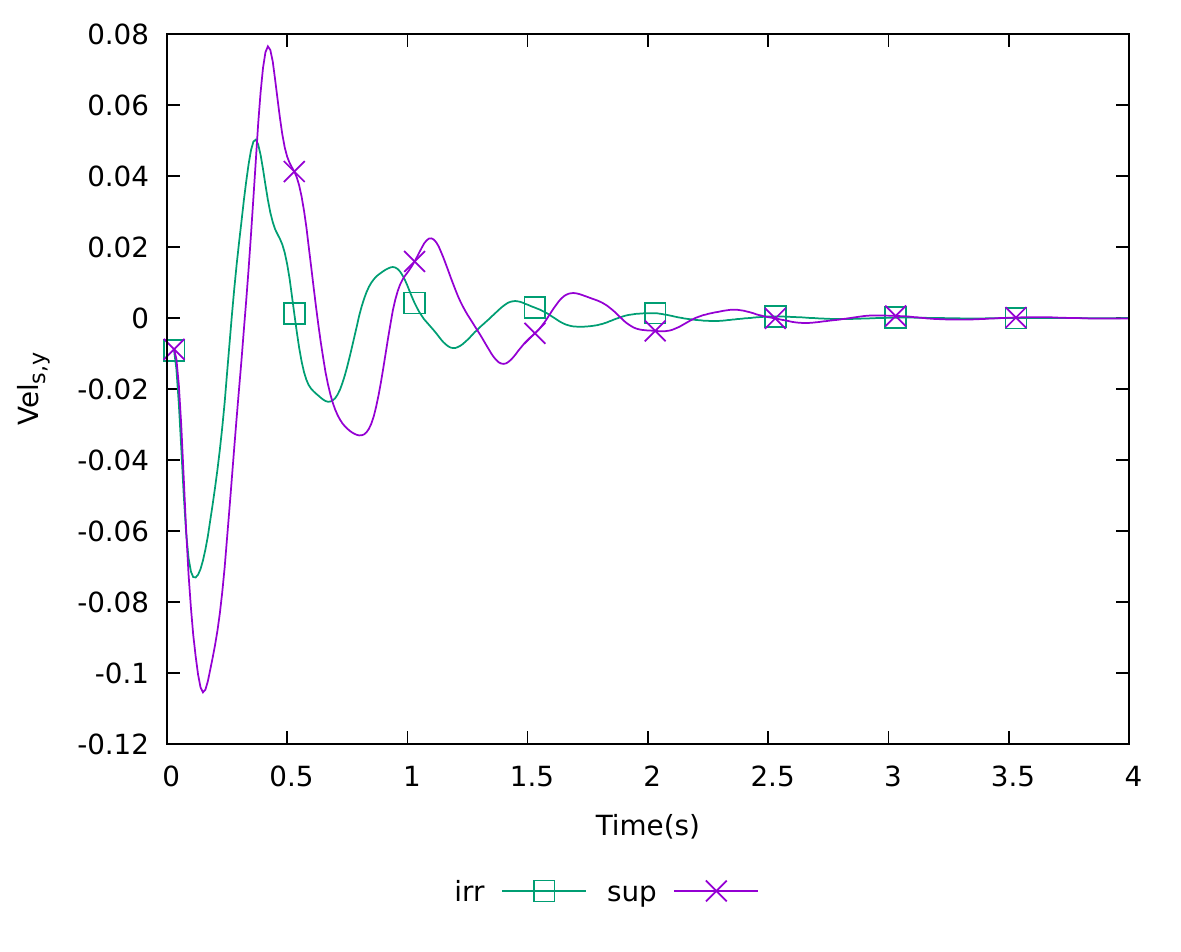}
        \caption{Velocity in $y$}
        \label{fig:5_vel_y} 
    \end{subfigure}
    \begin{subfigure}[h!]{0.3\textwidth}
        \includegraphics[width = \linewidth]{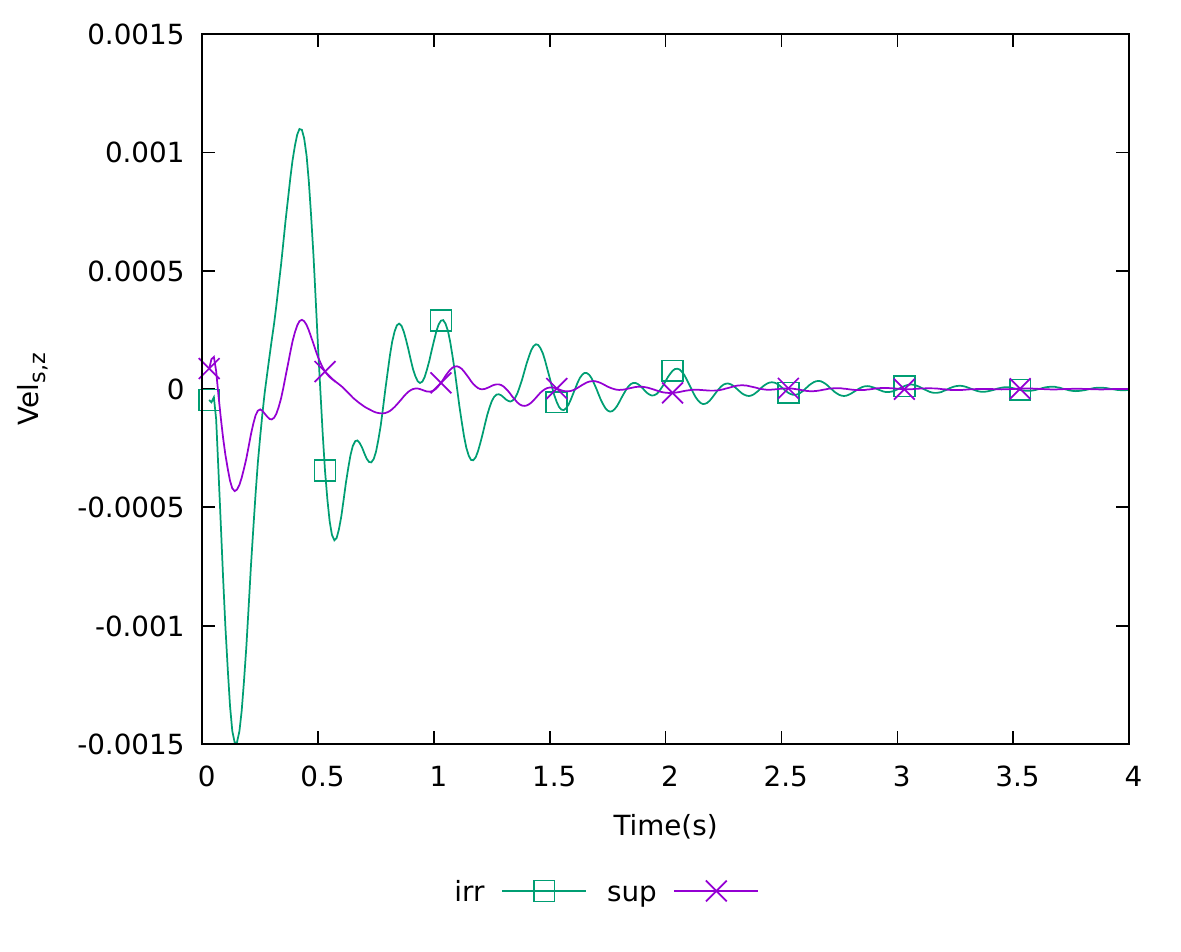}
        \caption{Velocity in $z$}
        \label{fig:5_vel_z} 
    \end{subfigure}
    \caption{Velocity at the tip of the plate}
    \label{fig:5_vel} 
\end{figure}

\begin{figure}[h!]
    \centering
    \begin{subfigure}[h!]{0.3\textwidth}
        \includegraphics[width=\linewidth]{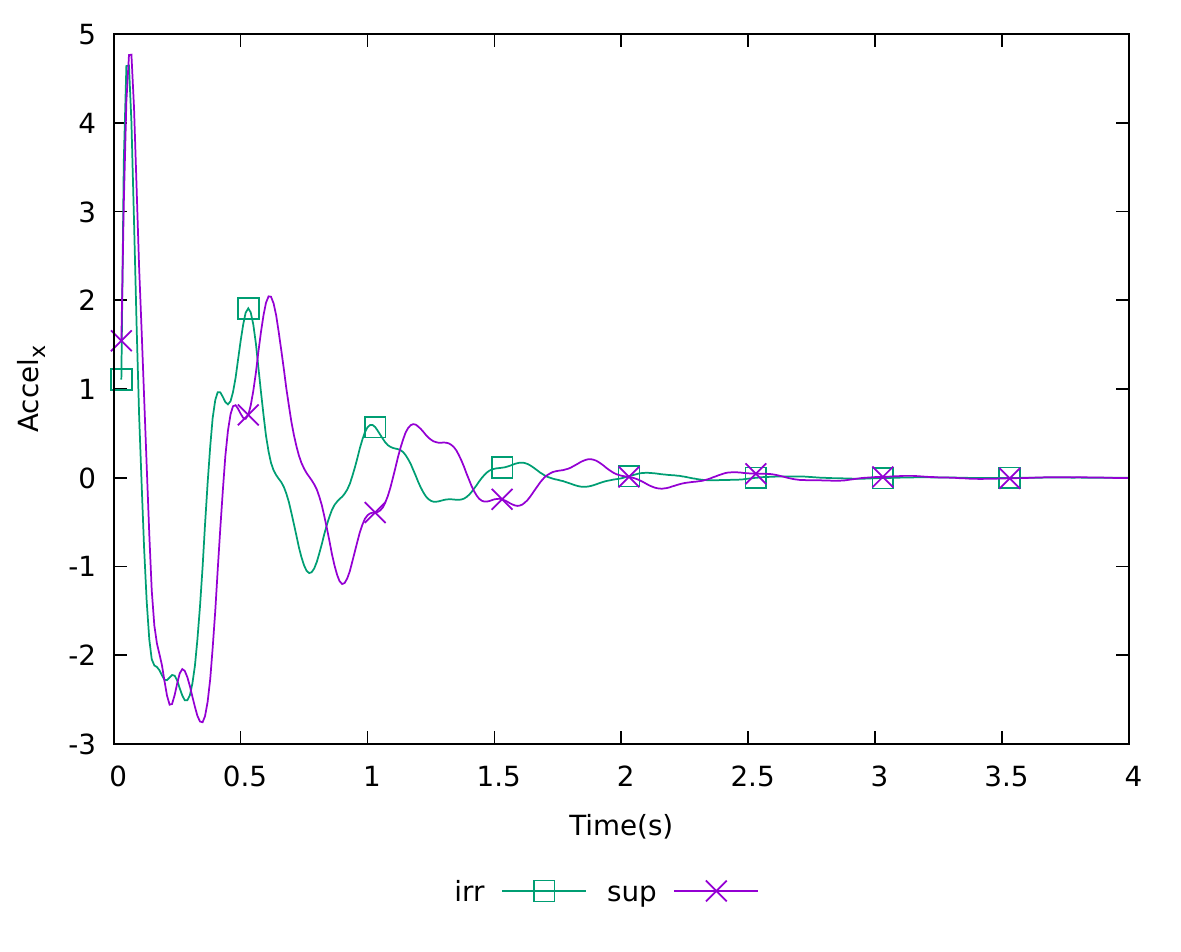}
        \caption{Acceleration in $x$}
        \label{fig:5_accel_x} 
    \end{subfigure}
    \begin{subfigure}[h!]{0.3\textwidth}
        \includegraphics[width = \linewidth]{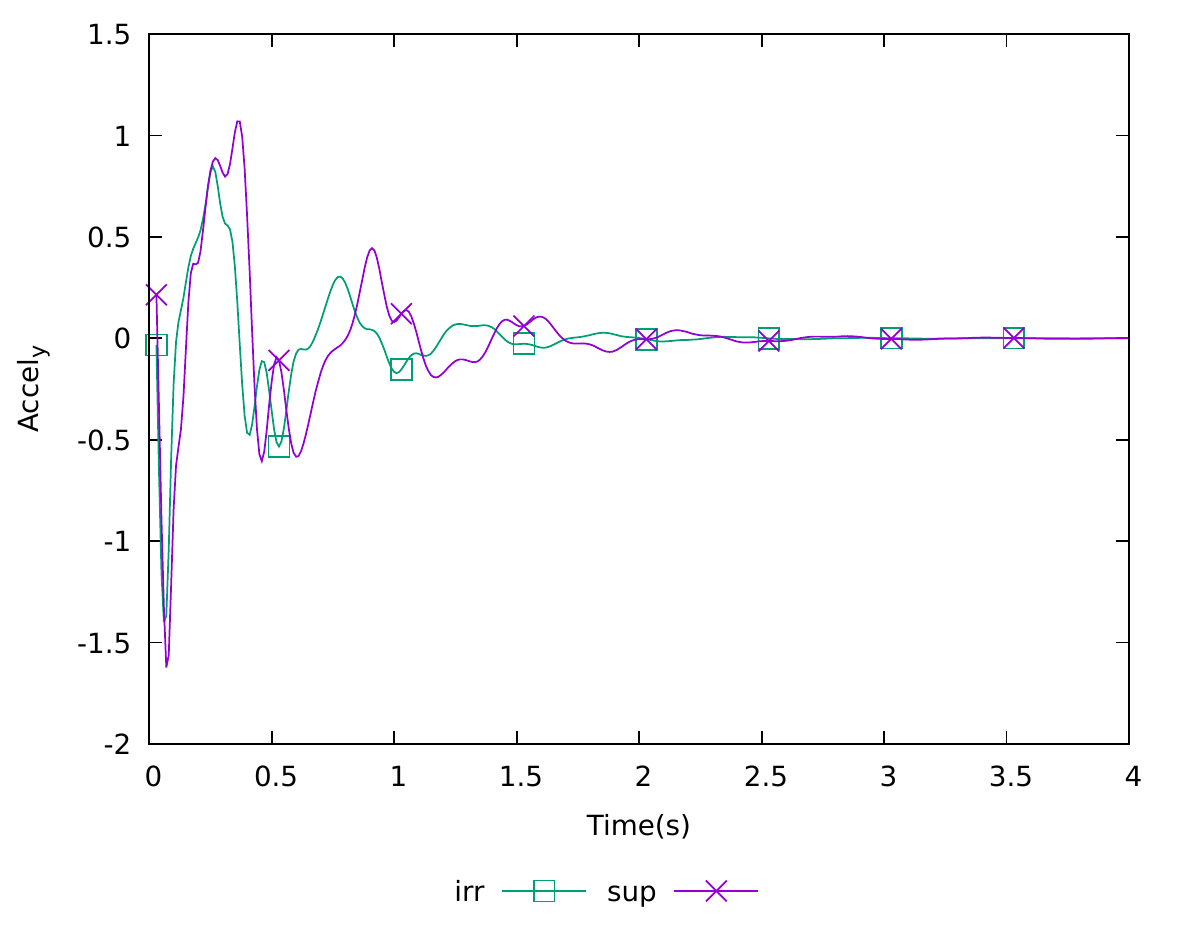}
        \caption{Acceleration in $y$}
        \label{fig:5_accel_y} 
    \end{subfigure}
    \begin{subfigure}[h!]{0.3\textwidth}
        \includegraphics[width = \linewidth]{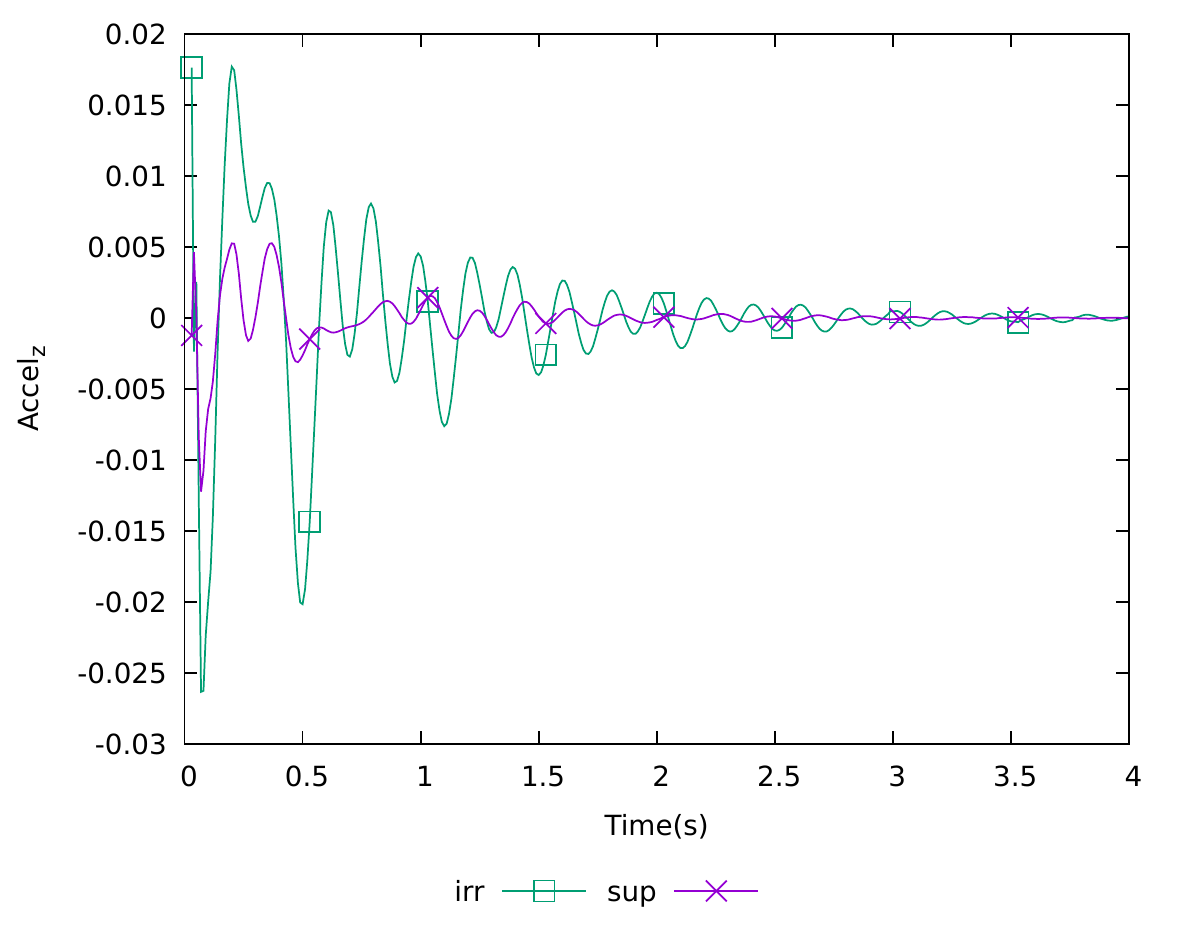}
        \caption{Acceleration in $z$}
        \label{fig:5_accel_z} 
    \end{subfigure}
    \caption{Acceleration at the tip of the plate}
    \label{fig:5_accel} 
\end{figure}

The next series of graphs in Figs.~\ref{fig:5_disp}, \ref{fig:5_vel} and \ref{fig:5_accel} show a comparison between the solution obtained by means of the three-field formulation (labeled `sup') and using the displacement based one (labeled `irr').

\begin{figure}[h!]
    \centering
    \begin{subfigure}[h!]{0.3\textwidth}
        \includegraphics[width=\linewidth]{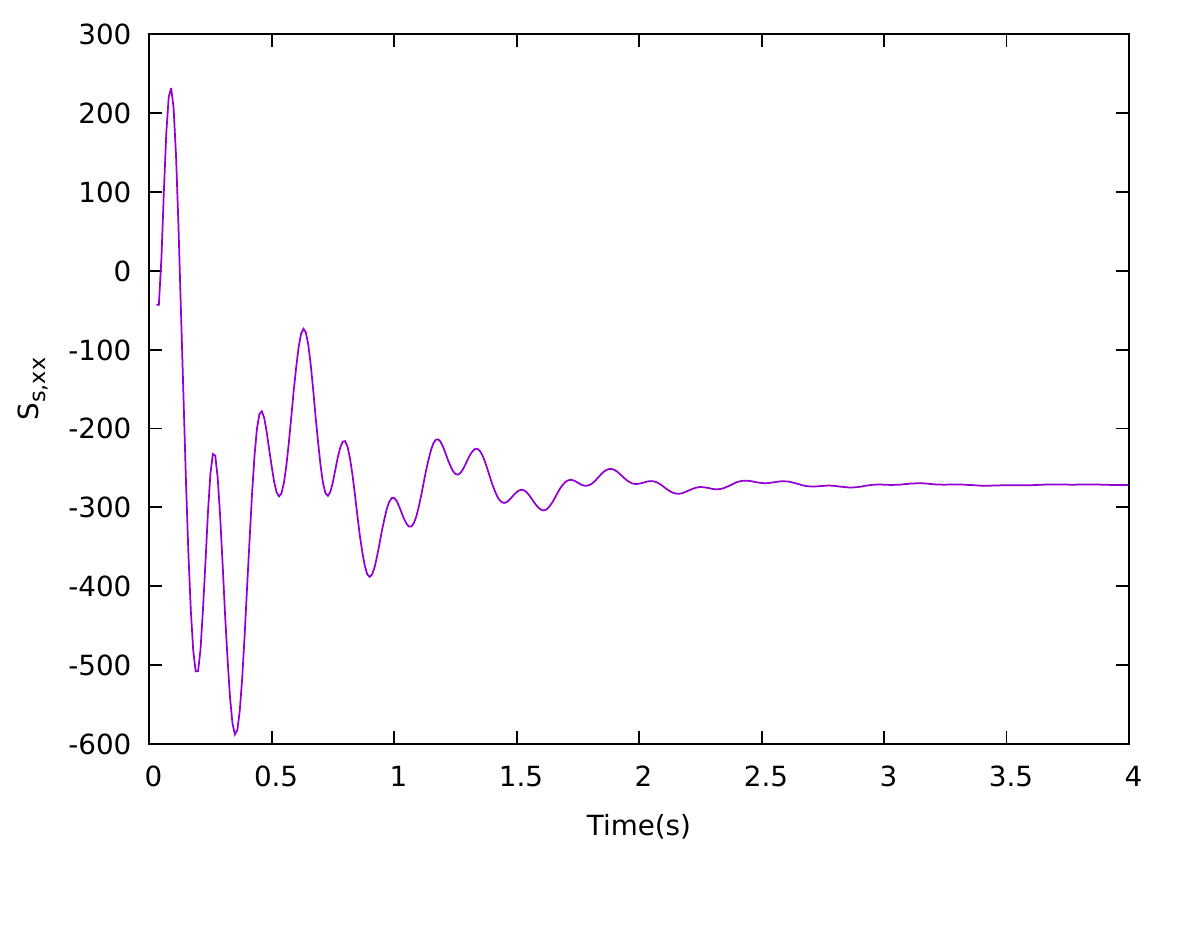}
        \caption{Stress in $xx$}
        \label{fig:5_stress_xx} 
    \end{subfigure}
    \begin{subfigure}[h!]{0.3\textwidth}
        \includegraphics[width = \linewidth]{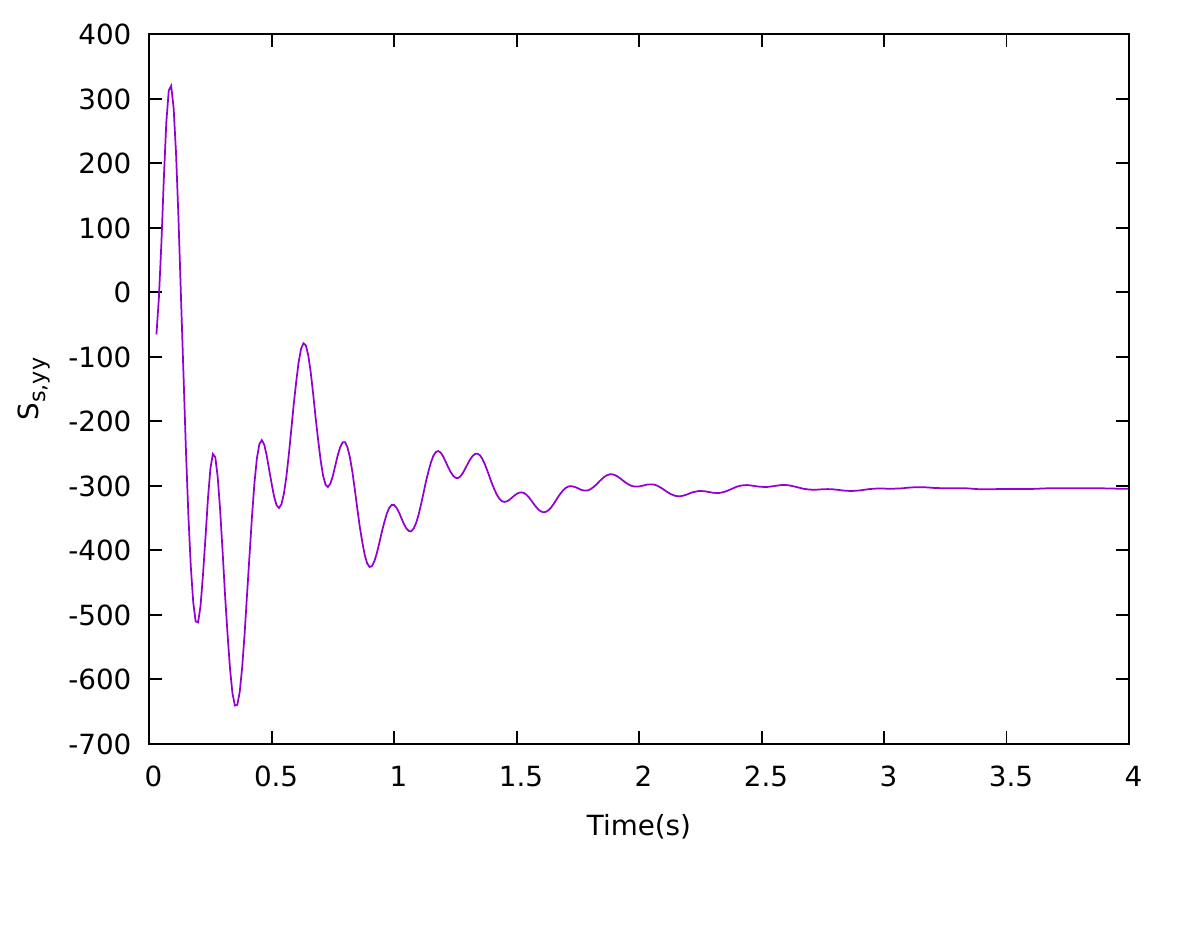}
        \caption{Stress in $yy$}
        \label{fig:5_stress_yy} 
    \end{subfigure}
    \begin{subfigure}[h!]{0.3\textwidth}
        \includegraphics[width = \linewidth]{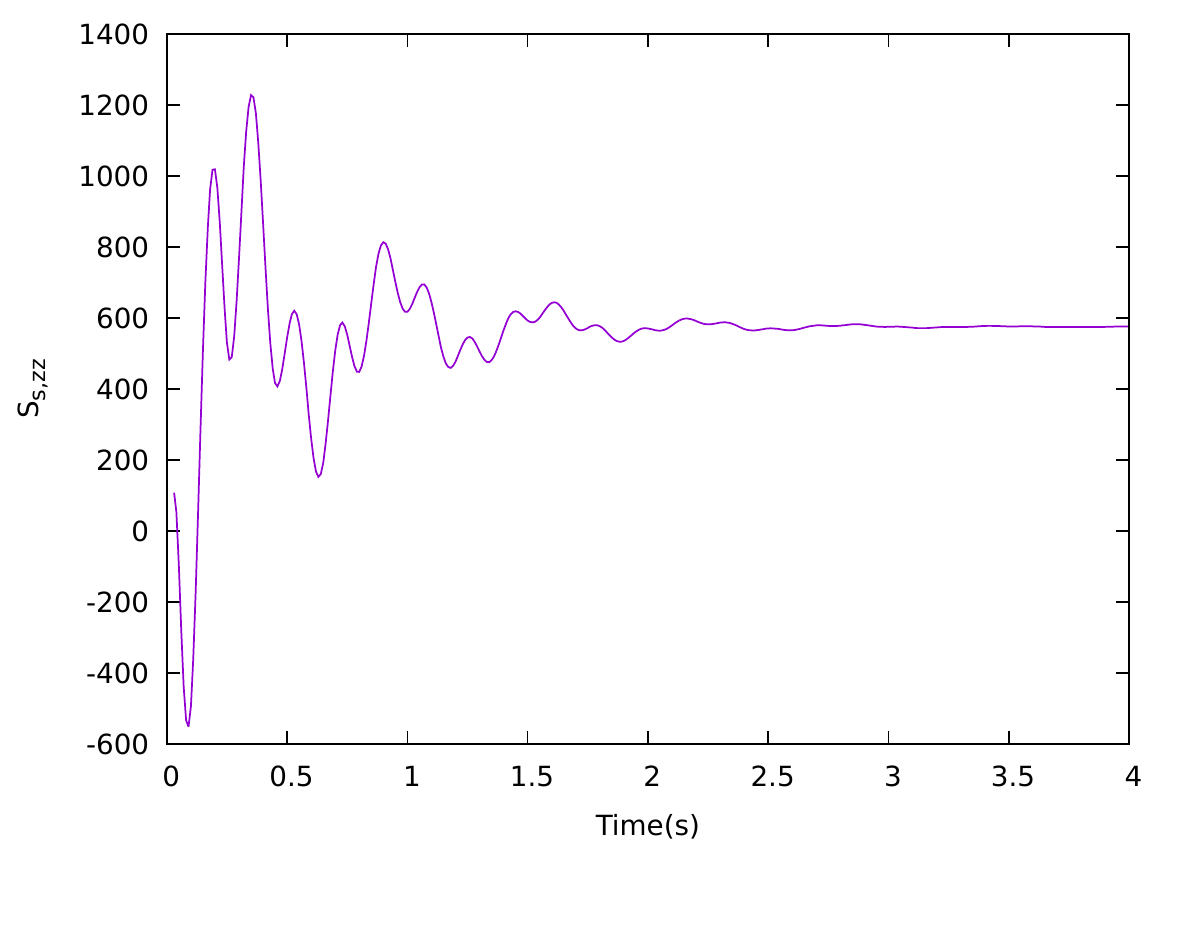}
        \caption{Stress in $zz$}
        \label{fig:5_stress_zz} 
    \end{subfigure}
    \caption{Normal stresses at the tip of the plate}
    \label{fig:5_stress_nor} 
\end{figure}

\begin{figure}[h!]
    \centering
    \begin{subfigure}[h!]{0.3\textwidth}
        \includegraphics[width=\linewidth]{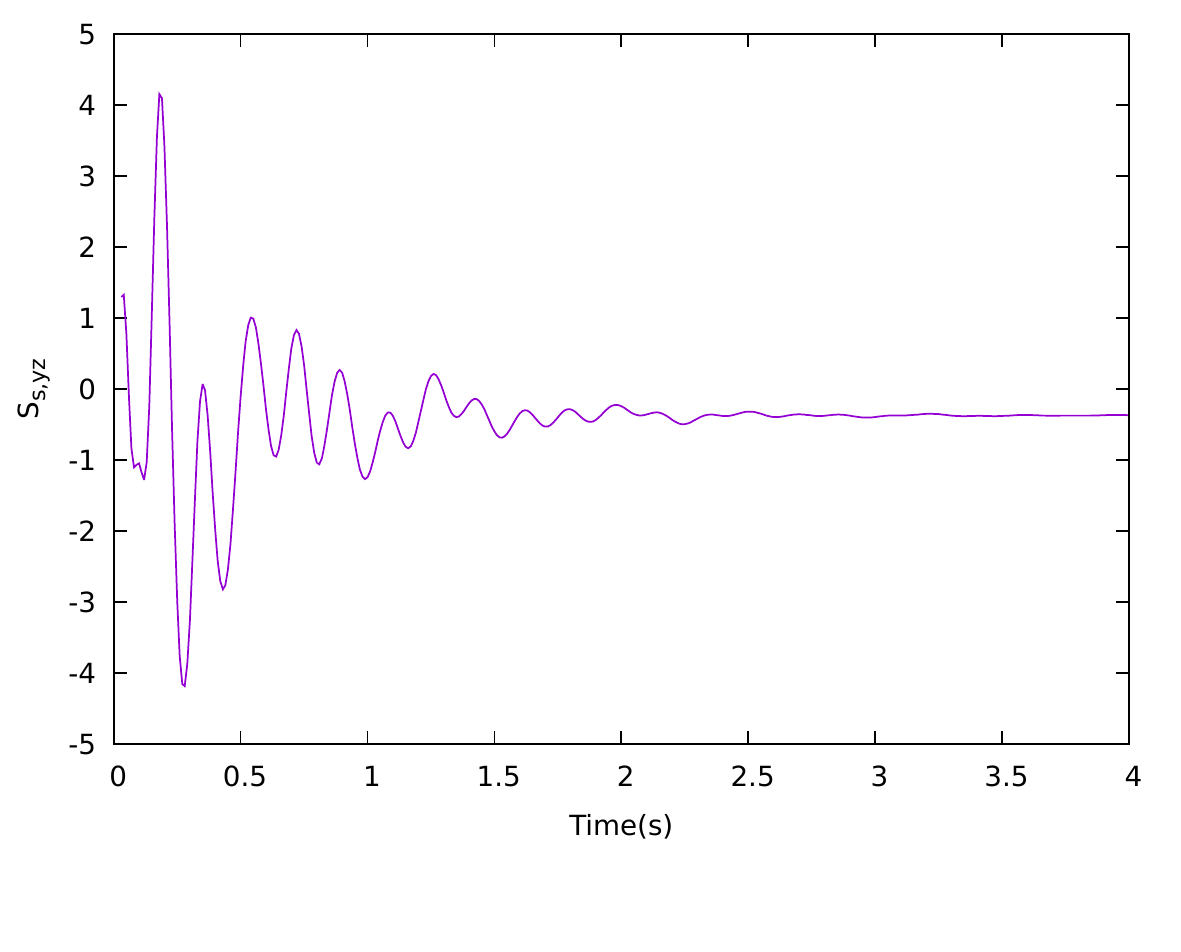}
        \caption{Stress in $yz$}
        \label{fig:5_stress_yz} 
    \end{subfigure}
    \begin{subfigure}[h!]{0.3\textwidth}
        \includegraphics[width = \linewidth]{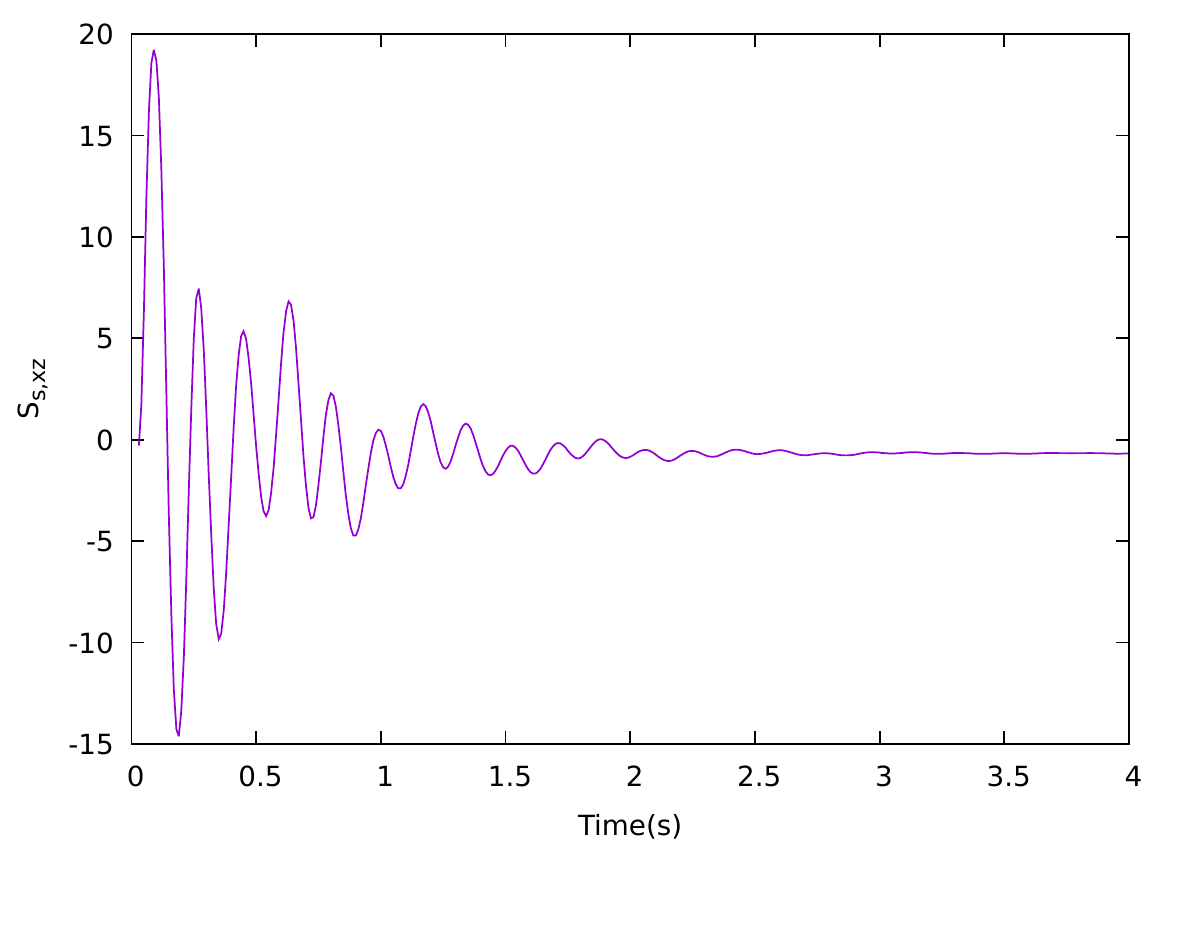}
        \caption{Stress in $xz$}
        \label{fig:5_stress_xz} 
    \end{subfigure}
    \begin{subfigure}[h!]{0.3\textwidth}
        \includegraphics[width = \linewidth]{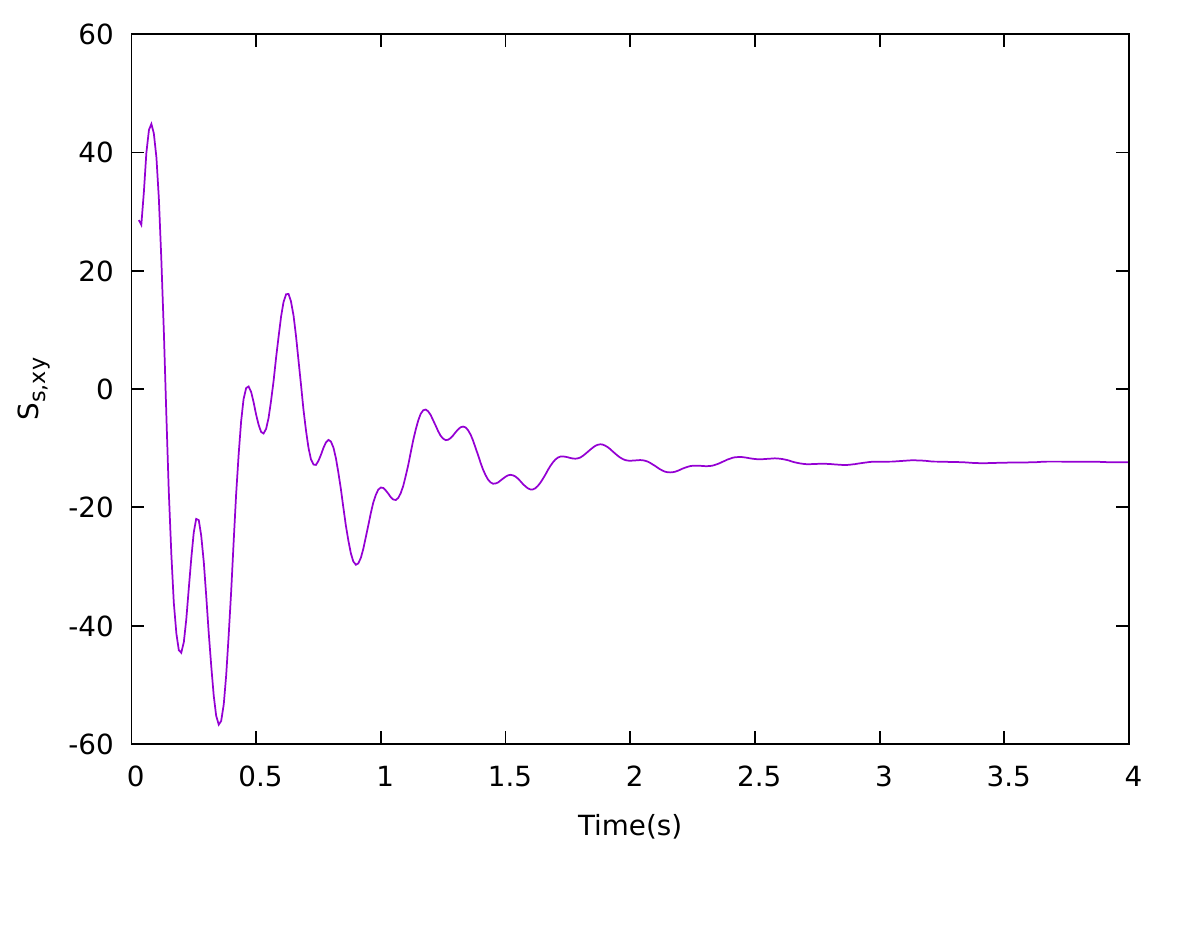}
        \caption{Stress in $xy$}
        \label{fig:5_stress_xy} 
    \end{subfigure}
    \caption{Normal stresses for the tip of the plate}
    \label{fig:5_stress_shear} 
\end{figure}

\begin{figure}[h!]
    \centering
    \includegraphics[width=0.5\textwidth]{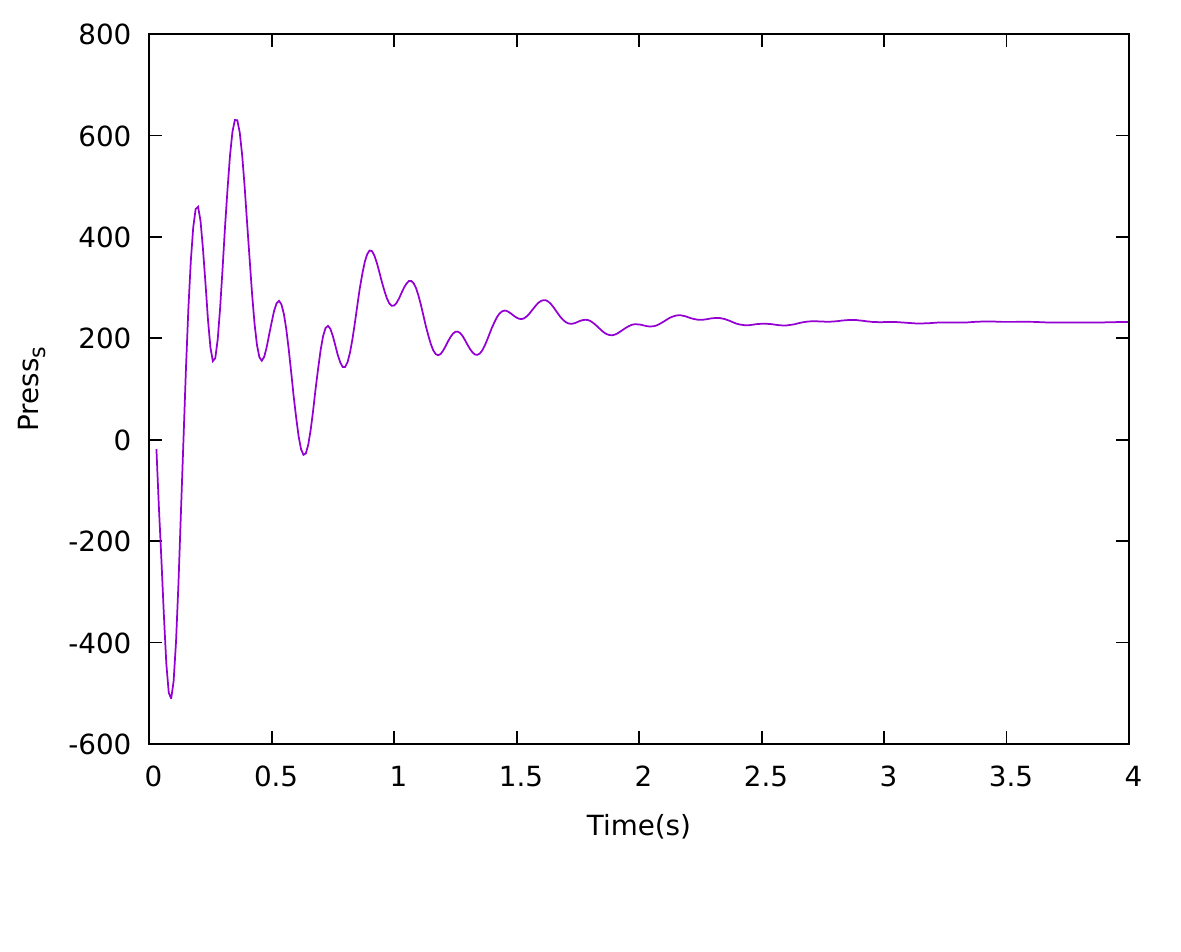}
    \caption{Pressure}
    \label{fig:5_press} 
\end{figure}

\begin{figure}[h!]
    \centering
    \includegraphics[width=0.85\textwidth]{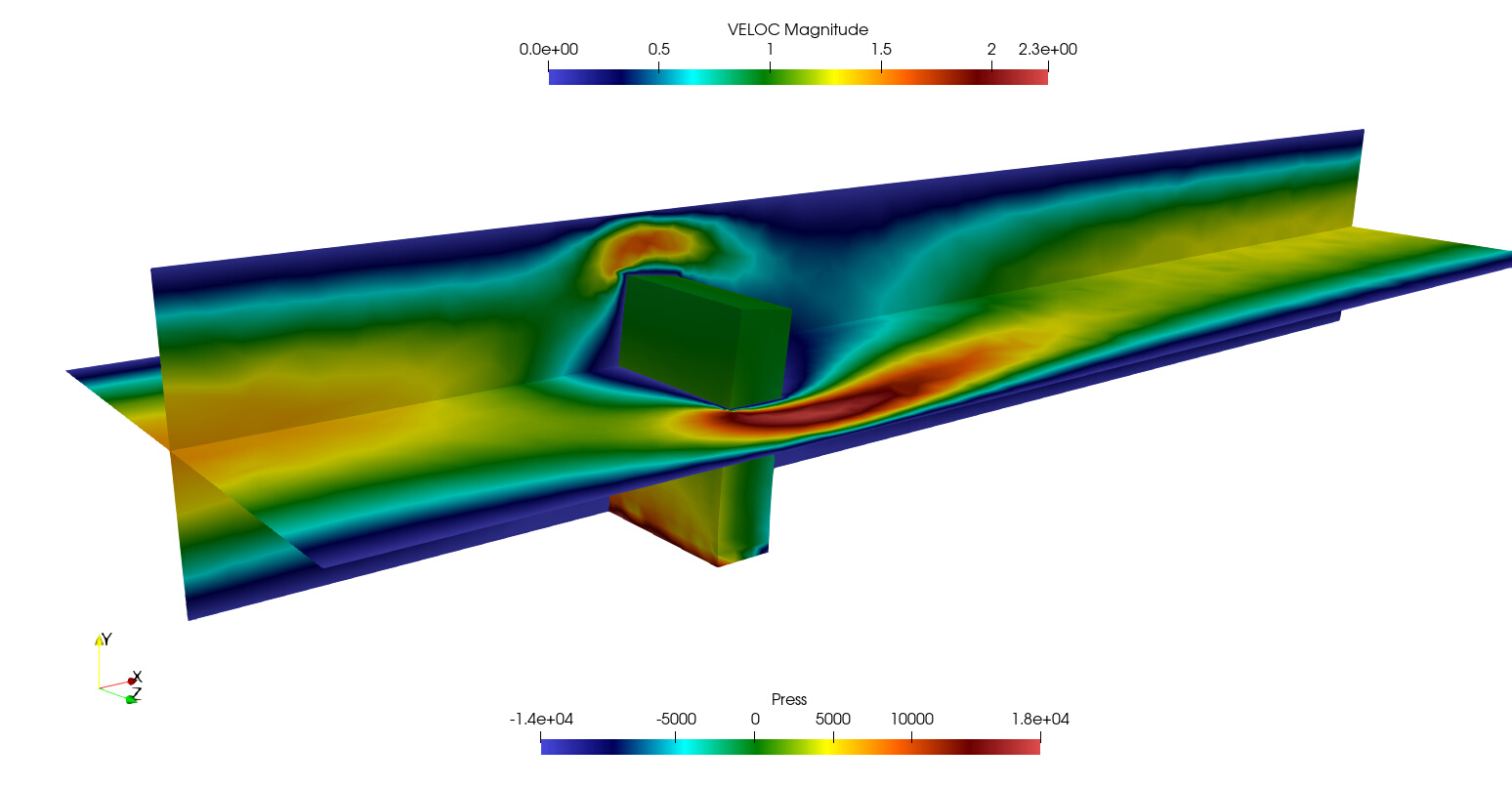}
    \caption{Fluid: Velocity contours; Solid: Pressure contours}
    \label{fig:5_cont_vel_press} 
\end{figure}

\begin{figure}[h!]
    \centering
    \includegraphics[width=0.85\textwidth]{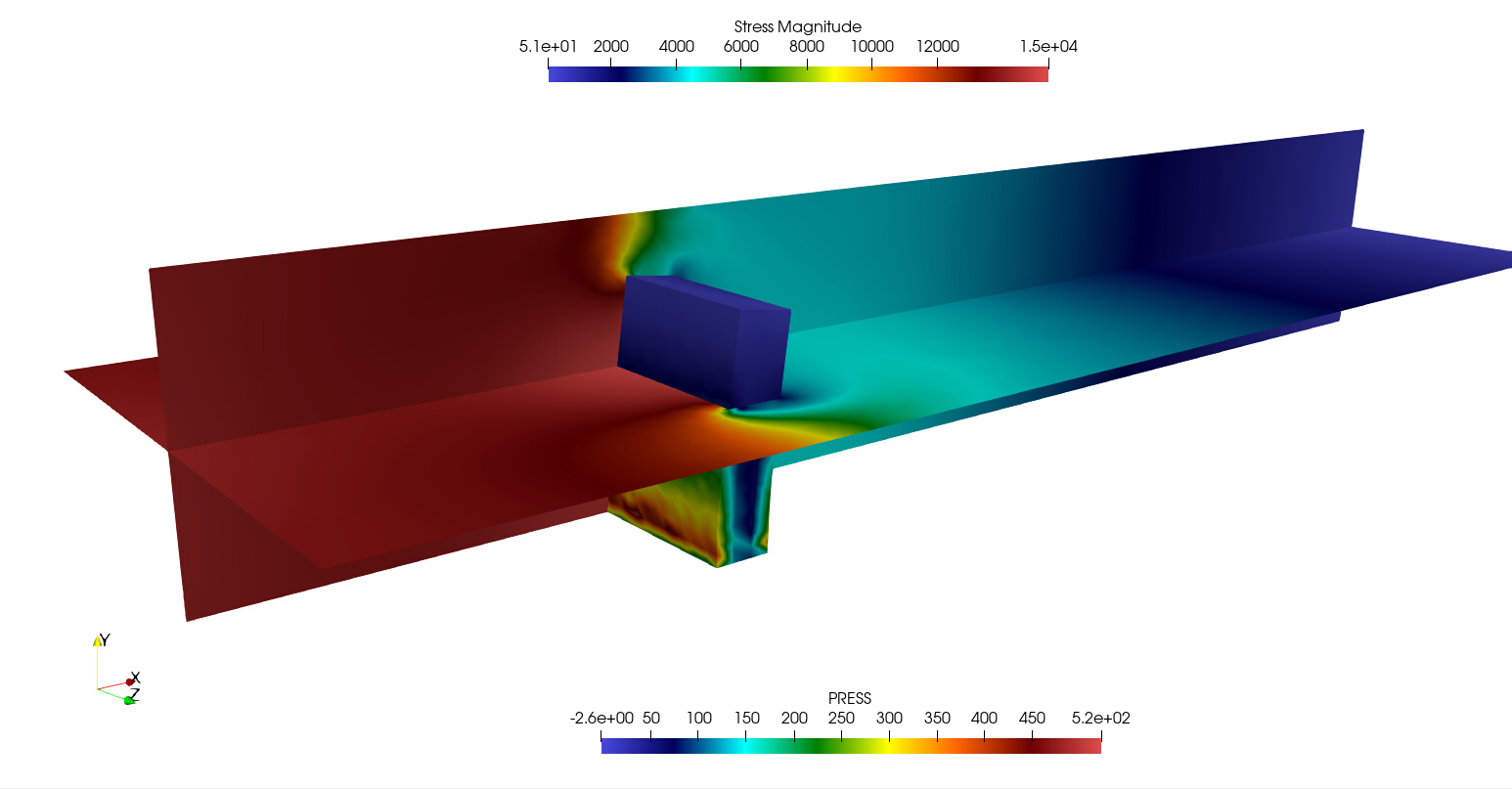}
    \caption{Fluid: Pressure contours; Solid: Stress magnitude contours}
    \label{fig:5_cont_press_stress} 
\end{figure}

\begin{figure}[h!]
    \centering
    \includegraphics[width=0.85\textwidth]{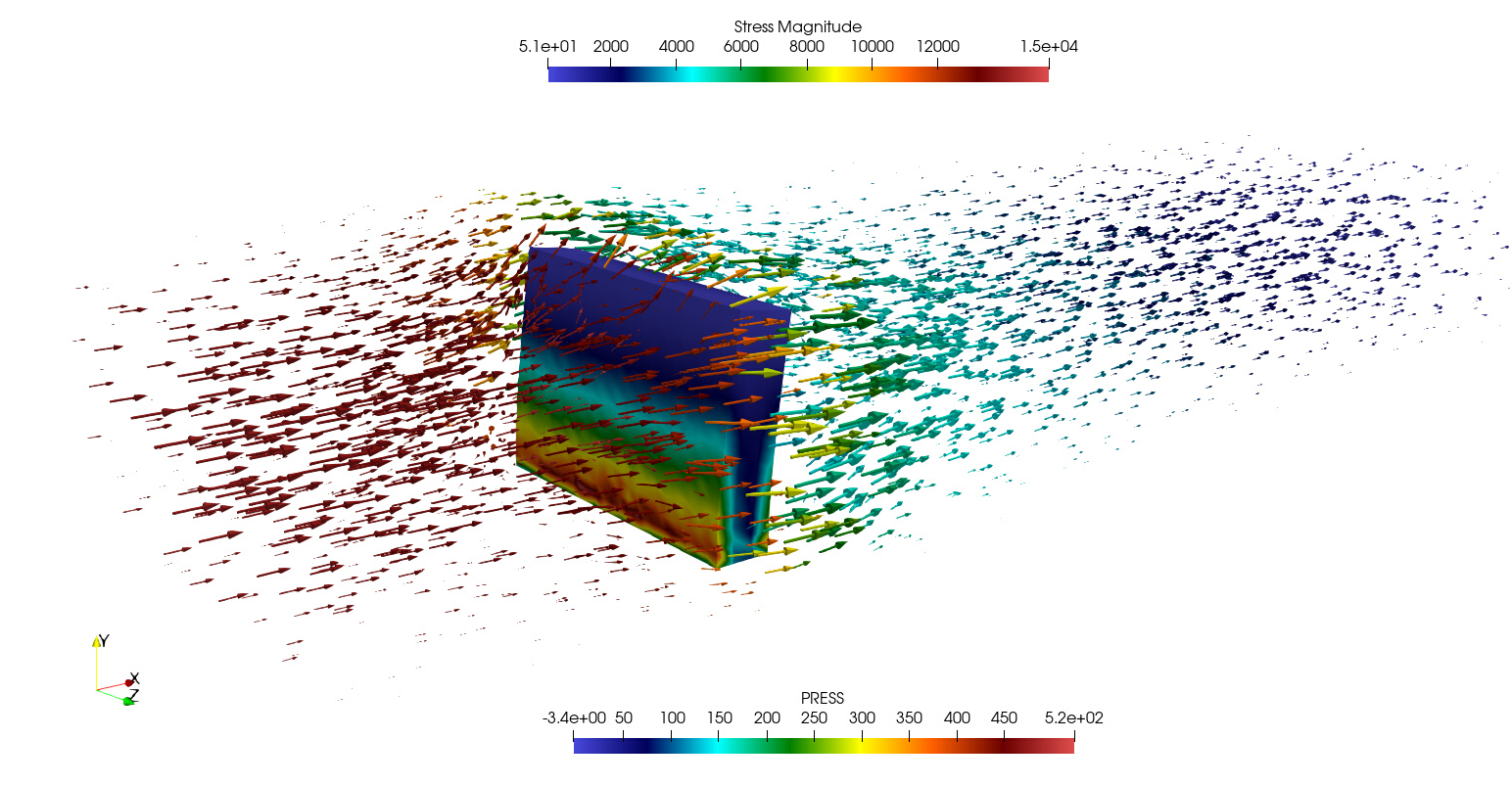}
    \caption{Fluid: Velocity vectors show pressure; Solid: Stress magnitude contours}
    \label{fig:5_glyph} 
\end{figure}

Notice how the overall behavior of the standard displacement based formulation is over-diffusive and tends to over-dampen the motion of the plate, see Fig.~\ref{fig:5_disp}. In turn, the velocities and accelerations of the plate in the three-field formulation have a higher amplitude, see Figs.~\ref{fig:5_vel} and \ref{fig:5_accel}. The displacement based formulation reaches a stationary state earlier than the three-field formulation, and it also produces a result with lower frequency, which leads to important phase differences.

Even if we do not have any reference for comparison, it is important to note that the stress and pressure fields shown in Figs.~\ref{fig:5_stress_nor}, \ref{fig:5_stress_shear} and \ref{fig:5_press} are smooth and continuous and decay to reach a stationary state, as the displacement field.

Finally, and in order to visualize the solution to this problem, Figs.~\ref{fig:5_cont_vel_press}, \ref{fig:5_cont_press_stress} and \ref{fig:5_glyph} show contours of velocity norm, pressure contours in the fluid and of stress magnitude in the solid, and velocity vectors, respectively, of the solution obtained.

\section{General conclusions}\label{sec_conclusions}

This work first presents a displacement-stress-pressure formulation for a neo-Hookean solid using an updated Lagrangian formulation, followed by a new three-field FSI formulation stabilized by means of a VMS approach using time dependent sub-grid scales on both the fluid and the solid regions. Benchmarking was done for the solid formulation under static and dynamic scenarios; on itself the solid three-field formulation proves to be more accurate and less time step dependent than its irreducible counterpart.

The three-field formulation has proved to be robust and precise under all cases analyzed. Compared to the irreducible formulation, it preserves phase and amplitude in time with fewer elements and even with linear elements. This formulation proves to be accurate and efficient, as it is possible to use a coarser mesh that produces more accurate results than its standard irreducible counterpart with a finer mesh.

In conclusion, a new solid elasto-dynamic formulation has been benchmarked and coupled with the previously developed three-field fluid to produce a robust and accurate FSI formulation that is stable and efficient. 

\section*{Acknowledgements} 
A. Tello wants to acknowledge the doctoral scholarship received from the Colombian Government-Colciencias. R. Codina gratefully acknowledges the support received from the ICREA Acad\`emia Program, from the Catalan Government.


\section*{Appendix: Linearization of $J$ and $\ln(J)$}\label{app_lin}

The determinant $J$ of $F_{iJ}$ can be obtained as follows:
\begin{align}
    J &= \textrm{det}F_{iJ} = \textrm{det}\left(\tilde{F}_{iJ} + \dpar{\delta d_i}{X_J}\right), \nonumber\\
     &= \textrm{det}\left[\tilde{F}_{iK} \left( \delta_{KJ}+ \tilde{F}_{Kl}^{-1}\dpar{\delta d_l}{X_J}\right)\right],\nonumber
\end{align}
where by means of the multiplicative property of determinants we can express the previous expression as
\begin{align}
     J &= \tilde{J}\textrm{det}\left( \delta_{KJ}+ \tilde{F}_{Kl}^{-1}\dpar{\delta d_l}{X_J}\right),\nonumber
\end{align}
If we construct this determinant we notice that by discarding higher order terms we are left with:
\begin{align*}
     J &= \tilde{J}\left(1 + \tilde{F}_{Ji}^{-1}\dpar{\delta d_i}{X_J}\right),
\end{align*}
Finally, using the previous result we can linearize the logarithm of the determinant of the displacement gradient:
\begin{align}
     \textrm{ln}(J) &= \textrm{ln} \left[\tilde{J}\left(1 + \tilde{F}_{Ji}^{-1}\dpar{\delta d_i}{X_J}\right)\right], \nonumber \\
     &= \textrm{ln}(\tilde{J}) +  \textrm{ln}\left[\left(1 + \tilde{F}_{Ji}^{-1}\dpar{\delta d_i}{X_J}\right)\right], \nonumber
\end{align}
whereby using the fact that for small $x$ there holds $\textrm{ln}(1+x) = x + {\cal O}(x^2)$, we can express the previous result as:
\begin{align*}
     \textrm{ln}(J) &= \textrm{ln}(\tilde{J}) +  \tilde{F}_{Ji}^{-1}\dpar{\delta d_i}{X_J}.
\end{align*}


\bibliographystyle{abbrv}
\bibliography{trifil_vms.bib}

\end{document}